# Algebraic and Topological Persistence

by
**Luciano Melodia**

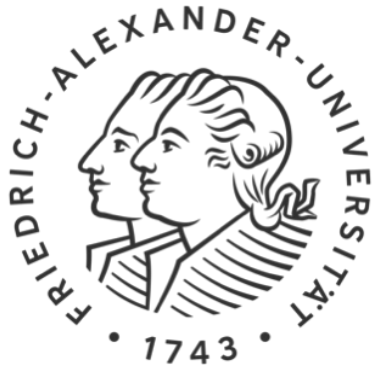

A bachelor thesis submitted in fulfillment
of the requirements for the degree of
Bachelor of Science (B.Sc.)
in Mathematics

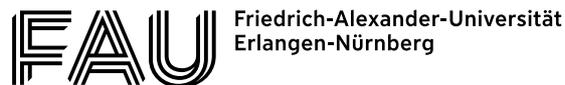

Friedrich-Alexander-Universität
Erlangen-Nürnberg

October 8, 2024

# Algebraic and Topological Persistence

## Luciano Melodia

## Abstract

This work will begin with a brief introduction, inspired by an applied context that has provided the foundation for this thesis. A substantial proportion of research in persistence theory has been conducted with the objective of applying persistent homology for the analysis of empirically-derived data. Nevertheless, we will adhere strictly to the theory of topological spaces and restrict ourselves to the provision of illustrative contexts for applications. This treatise addresses the theory of topological spaces (§2) and the foundations of persistence theory (§3). The subsequent discussion will address chain complexes (§2.3.8) and the associated simplicial homology groups (§2.2), as well as their relationship with singular homology theory (§2.3). Moreover, we present the fundamental concepts of algebraic topology, including exact and short exact sequences (§2.3.2) and relative homology groups derived from quotienting with subspaces of a topological space (§2.4). These tools are used to prove the Excision Theorem (§2.5) in algebraic topology. Subsequently, the theorem is applied to demonstrate the equivalence of simplicial and singular homology for triangulable topological spaces, i.e. those topological spaces which admit a simplicial structure (§2.6). This enables a more general theory of homology to be adopted in the study of filtrations of point clouds.

The chapter on homological persistence (§3) makes use of these tools throughout. We develop the theory of persistent homology (§3.2.2), the homology of filtrations of topological spaces (§3.1), and the corresponding dual concept of persistent cohomology (§3.2.3, §3.3). This work aims to provide mathematicians with a robust foundation for productive engagement with the aforementioned theories. The majority of the proofs have been rewritten to clarify the relationships between the techniques discussed. The novel aspect of this contribution is the canonical presentation of persistence theory and the associated ideas through a rigorous mathematical treatment for triangulable topological spaces and closing some gaps in the existing literature.

# Dedication

I dedicate this thesis to those closest to me, who have provided invaluable support throughout this journey. I could not have written these pages or undertaken this path without the indispensable support I received.

My academic career took me through German and Italian studies, then to a major in Information Science with a focus on Media Informatics, culminating in a Master's degree from Regensburg. I then did a three-year doctoral programme in Computer Science at Erlangen, after which I shifted my focus to theoretical Mathematics. Throughout, I was driven by a desire to understand the underlying structures of language and computing. I am proud of these experiences and the insights they have given me.

I am grateful to my mother, Beata, who supported my decisions wholeheartedly once they were made, despite initial reservations.

I would like to thank my late father, Domenico, for teaching me to value the past more than the future.

I'd also like to thank my sister Dominique and her husband Maxwell, who never let up on me about physics and political subjects. They made sure I kept up with these subjects whenever we had time to chat about them.

To Luciana, who brought joy and laughter back into my life during challenging times – thank you.

Also to Luciano and Arianna, who listened to me in difficult times and welcomed me to Sicily with delicious cassata.

And to all those unnamed, if these words strike a chord, I am aware that you have touched my life in significant ways.

I have learned a great deal from each of you.

# Declaration

I hereby declare that I
- alone wrote the submitted Bachelor's thesis without illicit or improper assistance.
- did not use any materials other than those listed in the bibliography and that all passages taken from these sources in part or in full have been marked as citations and their sources cited individually in the thesis. Citations include the version (edition and year of publication) and paragraphs of the used materials.
- have not submitted this Bachelor's thesis to another institution and that it has never been used for other purposes or to fulfill other requirements.

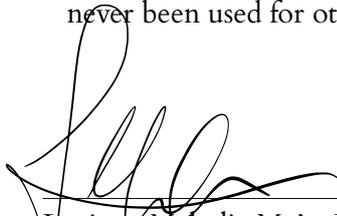

Luciano Melodia M.A., B.A.
Department of Mathematics
Friedrich–Alexander University Erlangen–Nürnberg

# Acknowledgements

I would like to express my gratitude to Prof. Dr. Catherine Meusburger, who not only taught me the rigours of the Bachelor's programme in Mathematics but also illuminated the more complex concepts of Linear Algebra I and II. Her assistance during our early sessions was crucial. We met in the study room to discuss exercise sheets and her guidance was invaluable. I am particularly grateful for the moments she joined us after a refreshing cup of coffee, skillfully translating challenging problems into clear solutions on the blackboard with swift strokes of chalk.

I would like to express my gratitude to my supervisor, Prof. Dr. Kang Li. His willingness to engage in discussions was instrumental in motivating me, particularly when I was overwhelmed by the vast sea of scientific papers and theorems. Prof. Li's expertise in simplifying complex research for my bachelor's thesis was priceless, and his guidance has been a driving force in my early academic journey in Mathematics.

I would be remiss if I did not acknowledge the invaluable contributions of my friends and colleagues, including Luciana Diaconescu, who helped proof-reading, Philipp Gäbelein, Michael Schleich, Aneek Butt, Sarah Wiesend, Adrian Hieber, Thomas Reichelt, Moritz Lanz, Sebastian Müller, Thomas Büttner, Anton Rechenauer, Tobias Simon, Michael Preeg, Constantin Winkler and Anton Hoof. I'd like to give my deepest cordial gratitude to Marie-Louise Isenberg. Without her, I probably wouldn't be here.

They all made my days better, supported me unwaveringly, and their presence was essential to my achievements. I know they will succeed wherever they are now, and I am always here for them.

# Contents

# Chapter 1

# Motivation

The development of persistence theory in topological data analysis (TDA) marks a significant advancement in the quantitative analysis of complex, high-dimensional data sets. Conventional statistical methods frequently prove inadequate for capturing the intricate geometric and topological structures of such data, particularly when these are obscured by noise or nonlinear relationships.

Persistence theory, and in particular persistent homology, provides a rigorous mathematical framework for the systematic identification and quantification of topologically relevant structures at varying levels of detail. Persistent homology is a mathematical concept that is employed for the analysis of homological features in data. These representations permit a precise description of the data structure and enable the distinction between significant features and noise. The persistent homology techniques are founded upon the principles of algebraic topology, wherein homology groups serve as algebraic invariants for the classification of topological spaces. In persistent homology, this concept is extended to a filtration where a real-valued function serves as the parameter of a nested sequence of topological spaces. The characteristics of this filtration are captured in persistence intervals, which indicate the degree of robustness associated with the features in question. The persistent homology is obtained by constructing simplicial complexes from the data through the efficient application of computational methods. In addition to their descriptive function, persistence diagrams and barcodes can be integrated into statistical and machine learning methods, enabling quantitative comparisons between data sets – as an example consider the work of Melodia et al. [21–23].

The theory of persistence is a valuable tool in TDA, as it remains stable in the analysis of data that has been subject to perturbations. This theory guarantees that minor alterations to the input data will only result in slight alterations to the persistence diagrams, thereby confirming the reliability of the topological features that have been extracted. This robustness is particularly beneficial when analysing data that is noisy or incomplete. Persistence theory serves to bridge the gap between abstract mathematical theory and the more practical field of data analysis. It has thus become a fundamental tool in this field.

## 1.1 History

The growth of persistence theory in TDA has been shaped by the seminal contributions of pioneering researchers, who have collectively established persistence as a core instrument for grasping the topological aspects of data.

**The field of size theory:** The origins of persistence theory can be traced back to the end of the 20th century, with the work of Patrizio Frosini and Massimo Ferri at the University of Bologna. Their research in set theory focused on the natural pseudodistance between functions defined on homeomorphic topological spaces [14, 16]. In particular, size theory, as applied to 0-dimensional homology, provided a framework for quantifying differences between topological spaces, and introduced a method for capturing persisting features across scales.

**Fractal geometry:** During her doctoral research at the University of Colorado Boulder, Vanessa Robins extended the application of persistence theory to fractal geometry. By means of alpha forms – a concept first proposed by Herbert Edelsbrunner et al. – Robins provided evidence of the efficacy of persistent homology in capturing the multiscale structure of fractal sets [13, 25].

**Algebraic foundations:** Herbert Edelsbrunner and his team at Duke University have made significant progress in formalising the algebra of persistence theory. They introduced fundamental concepts, including simplicial filtrations, which systematically construct topological spaces from data points by adding simplices in a hierarchical manner [12]. This filtration process tracks topological features, such as connected components, loops and voids, across different scales, leading to the birth and death of these features. The distinction between positive and negative simplices, introduced by Edelsbrunner's group, is crucial for understanding the emergence and decay of homological features within a filtration of point clouds.

## 1.2 Computational Aspects

The computational efficiency of persistence theory has been a significant factor in its extensive adoption across a range of scientific disciplines. This efficiency, underpinned by robust algorithmic principles, has established persistence theory as a key analytical tool for complex data sets.

**Algorithmic developments:** The algorithms for computing persistent homology, developed by Herbert Edelsbrunner and colleagues, are founded upon rigorous mathematical theory and optimised for practical application. These algorithms entail the construction and reduction of boundary matrices for the efficient computation of homology groups across a filtration of simplicial complexes. Among these algorithms, the matrix reduction algorithm plays a pivotal role in the computation of persistence intervals. The practical implementation of these algorithms has led to the development of standard software packages, including `GUDHI`, `Ripser`, and `Dionysus` [24].

**Applications in life sciences:** In the field of life sciences, persistent homology has emerged as a powerful tool for the study of the structure and function of biological molecules, particularly proteins. Proteins are complex macromolecules, the function of which is intricately linked to their three-dimensional structure. Persistent homology is a method for identifying patterns within a protein's structure that persist across different scales. These patterns frequently correspond to critical functional regions, such as binding pockets or active centers, and provide insights into protein interactions and stability under diverse conditions [20]. In the field of neuroscience, persistence theory has been employed to examine brain networks, thereby providing a novel methodology for analysing the intricate connectivity patterns that underpin brain function. By constructing simplicial complexes from neural data, researchers utilise persistent homology to monitor alterations in brain region connectivity over time. This approach has yielded new insights into the manner in which brain networks reorganise during cognitive processes or in response to neurological diseases [17, 26].

**Applications in machine learning:** In the field of machine learning, persistent homology is a valuable tool for uncovering the underlying topological structure of data distributions. This allows for the effective completion of tasks such as clustering, classification, and anomaly detection. The topological features of data provide essential insights into the organisation of the data space. For example, persistent homology can identify clusters of data points that form distinct topological features, such as loops or voids, which correspond to different classes or subpopulations within the data set. This has been demonstrated in studies referenced in the literature [19, 21–23]. The incorporation of topological data into machine learning models improves their resilience and generalisation capabilities, particularly in the context of noisy or incomplete data.

**Robustness to noise:** A defining feature of persistence theory is its resilience to noise, which represents a substantial advantage in the context of real-world data analysis. Persistence theory is concerned with the topological features that remain consistent across different scales, thereby filtering out noise-induced artefacts. Robustness is grounded in the stability theorems, which guarantee that minor changes are induced in persistence by small perturbations of the data [8, §3.1].

## 1.3   Advanced Extensions

The adaptability and extensibility of persistence theory are exemplified by its sophisticated developments. Among these extensions, discrete morse theory, multi-parameter persistence, and zigzag persistence represent notable ones.

**Discrete morse theory:** The discrete morse theory, initially proposed by Robin Forman, represents an extension of the classical morse theory to discrete spaces, such as simplicial complexes. This extension is of particular value in TDA, as it facilitates the simplification of complex spaces while ensuring the preservation of essential topological characteristics. The construction of discrete Morse functions on simplicial complexes enables the reduction of the number of critical simplices, thereby facilitating more efficient persistent homology computations [15]. The critical simplices correspond to significant topological features, and focusing on these reduces the computational complexity of calculating persistence intervals. This approach is especially advantageous in large data sets, where the number of simplices can render traditional homology computations impractical.

**Multiparameter persistence:** Multiparameter persistence represents an extension of the traditional single-parameter filtration approach, whereby filtrations are considered that are indexed by multiple parameters simultaneously. This allows for a more comprehensive examination of topological characteristics, as distinct parameters facilitate the capture of diverse aspects of the data's structural elements. To illustrate, one might filter a data set by both scale and density, thereby uncovering topological features that would otherwise remain invisible under a one-parameter filtration [4]. The use of multiple parameters results in a more complex algebraic structure, which in turn gives rise to computational and visualisation challenges associated with multiparameter persistence modules. In contrast to single-parameter persistence, where persistence diagrams or barcodes offer a comprehensive invariant, multiparameter persistence lacks a straightforward representation. In contrast, more complex invariants, such as generalized persistence diagrams or rank invariants, are used to capture relationships between parameters.

**Zigzag persistence:** Zigzag persistence represents an extension of the traditional persistent homology framework to allow for the analysis of dynamic data sets, whereby the data may undergo change over time. In contrast to the standard persistence approach, which requires a monotonically increasing sequence of spaces, zigzag persistence permits both forward and backward inclusions throughout the filtration process [3]. This flexibility permits the analysis of topological features in settings where data is not static, such as time-varying networks or data sets undergoing changes due to external factors. Zigzag persistence is a methodology that captures the evolution of topological features, including their appearance, disappearance, and reappearance as the data changes. The algebraic structure underlying zigzag persistence is more intricate than that of standard persistence. It involves directed graphs and a more complex form of homological algebra.

## 1.4   Contribution

This work presents a comprehensive summary of persistent homology theory, harmonising various proofs to provide a unified examination of persistence modules and persistent (co)homology from an algebraic-topological perspective. The objective of this study is to provide a mathematically rigorous explanation for the remarkable success of persistent homology in data analytics. Our approach is based on the assumption that real-world data can be represented on a triangulable topological space and that simplicial structures allow us to derive invariants of this space. We provide a rigorous proof of the merits of this approach and demonstrate its clear and well-defined foundations. Moreover, we shed light on the interplay between relative and absolute (co)homology on persistence intervals and fill some gaps within the literature.

# Chapter 2

# Topological Spaces and Groups

Topological persistence is deeply rooted in algebraic topology. This is the study of topological spaces and functions based on algebraic objects and their properties. Examples include homotopy and homology theories, which are essential for understanding the construction and connectedness of surfaces, a central aspect of data analysis. Persistent homology, the central tool of topological persistence, extends classical homology to identify features across multiple scales. Introduced by Edelsbrunner et al. in their seminal paper [10], persistent homology examines multi-scale topological features through a filtration process - an indexed family of nested spaces that starts with the empty set and progressively covers the entire space under study. Each filtration stage represents a snapshot of some topological space at a certain resolution, capturing the appearance and disappearance of multi-dimensional homology groups encoding topological properties such as connected components, holes and voids. Mathematically, the persistence of homological features is visualised using diagrams or barcode representations. These visual aids represent the emergence (birth) and disappearance (death) of topological features as the filter parameter changes. The duration of a feature's presence, represented by the length of its interval in the barcode, indicates its significance, with longer intervals suggesting features that represent the likely true characteristics of the underlying data rather than mere noise. The robustness of persistent homology, in particular its resistance to small perturbations in the data, is captured in the Stability Theorem [8, §3.1]. This theorem, proved by Cohen-Steiner, Edelsbrunner and Harer [1, §3], states that small variations in the input data lead to small changes in the persistence barcodes. This property is crucial for practical applications, as it guarantees that the topological summaries are both reliable and meaningful for the actual underlying structures.

We start with basic concepts such as topological spaces and groups, which are crucial for understanding and encoding connectedness and other invariants. The discussion extends to simplicial complexes, which are essential for modelling data structures in topological data analysis. We explore simplicial and singular homology groups to accurately quantify topological features, and dwell on singular chain complexes and exact sequences to deepen the algebraic aspects of persistence theory. This will turn useful with regard to the algebra of persistence modules.

## 2.1 Simplicial Complexes

We note that a set of points $X = \{x_0, x_1, \ldots, x_d\}$ in $\mathbb{R}^n$ is affinely independent if no affine subspace of dimension less than $d$ contains all the points in $X$. Such a set of points is commonly called point cloud.

**Definition 2.1.1** ($d$-simplex). *[2, Definition 2.1] A $d$-dimensional simplex $\sigma^{(d)}$, or $d$-simplex, is the set of all convex combinations of $X = \{x_0, x_1, \ldots, x_d\} \subset \mathbb{R}^n$, where $X$ consists of $d + 1$ affinely independent points.*

*Formally, $\sigma^{(d)}$ is defined by:*

$$\sigma^{(d)} := [x_0, \ldots, x_d] = \left\{ \sum_{i=0}^{d} \lambda_i x_i \ \middle| \ \sum_{i=0}^{d} \lambda_i = 1, \ \lambda_i \geq 0 \right\}. \tag{2.1}$$

As a convention, the empty set is considered a face, corresponding to the simplex formed by the empty subset of vertices. Specifically, a $0$-simplex corresponds to a single point, a $1$-simplex to a line segment between two points, a $2$-simplex to a triangle, and a $3$-simplex to a tetrahedron. Notably, a $d$-simplex is homeomorphic to the $d$-dimensional disk $D^d$.

**Lemma 2.1.2** ($d$-simplex to $d$-disk). *The $d$-simplex $\sigma^{(d)}$ is homeomorphic to the $d$-dimensional disk $D^d$.*

*Proof.* Define the standard $d$-simplex $\sigma^{(d)}$ as

$$\sigma^{(d)} := \left\{ (x_0, \ldots, x_d) \in \mathbb{R}^{d+1} \ \middle| \ \sum_{i=0}^{d} x_i = 1, \ x_i \geq 0 \right\}, \tag{2.2}$$

and the $d$-dimensional disk $D^d$ as

$$D^d := \left\{ (x_0, \ldots, x_{d-1}) \in \mathbb{R}^d \ \middle| \ \sum_{i=0}^{d-1} x_i^2 \leq 1 \right\}. \tag{2.3}$$

We construct a homeomorphism $f : \sigma^{(d)} \to D^d$ by

$$f(x_0, \ldots, x_d) = (\sqrt{x_0}, \ldots, \sqrt{x_{d-1}}), \tag{2.4}$$

where $x_d = 1 - \sum_{i=0}^{d-1} x_i$. This map is well-defined since

$$\sum_{i=0}^{d-1} (\sqrt{x_i})^2 = \sum_{i=0}^{d-1} x_i \leq 1. \tag{2.5}$$

The inverse $g : D^d \to \sigma^{(d)}$ is given by

$$g(y_0, \ldots, y_{d-1}) = (y_0^2, \ldots, y_{d-1}^2, 1 - \sum_{i=0}^{d-1} y_i^2), \tag{2.6}$$

ensuring that $g$ is well-defined because $\sum_{i=0}^{d-1} y_i^2 \leq 1$ implies $1 - \sum_{i=0}^{d-1} y_i^2 \geq 0$. Both $f$ and $g$ are continuous and are inverses of each other, as shown by $f(g(y)) = y$ for all $y \in D^d$ and $g(f(x)) = x$ for all $x \in \sigma^{(d)}$. □

Furthermore, it is important to note that $\sigma^{(d)}$ represents the convex hull of the points $X = \{x_0, x_1, \ldots, x_d\}$, defined as the smallest convex subset of $\mathbb{R}^n$ that contains all of these points. The faces of the simplex $\sigma^{(d)}$, with vertex set $X$, are formed by the simplices corresponding to subsets of $X$. A $d$-face of a simplex consists of a subset of the vertices with cardinality $d + 1$. The faces of a $d$-simplex with dimension less than $d$ are known as its proper faces. Two simplices are considered properly situated if their intersection is either empty or a face of both simplices. By identifying simplices along entire faces, we can construct the corresponding simplicial complexes.

**Definition 2.1.3** (Simplicial complex). *[2, Definition 2.2] A simplicial complex $K$ is a finite collection of simplices that satisfies the following properties:*

1. *For every simplex $\sigma^{(k)}$ in $K$, and every face $\tau^{(k)}$ of $\sigma^{(d)}$ with $k < d$, it follows that $\tau^{(k)}$ is also in $K$.*

2. *Any two simplices $\sigma^{(d)}$ and $\tau^{(k)}$ in $K$ are properly situated; that is, their intersection is either empty or a face of both simplices.*

The dimension of a simplicial complex $K$ is defined as the highest dimension among its simplices. For a simplicial complex $K$ in $\mathbb{R}^n$, the underlying space $|K|$ is the union of all the simplices in $K$. The topology of $K$ is determined by the topology induced on $|K|$ by $\mathbb{R}^n$'s standard topology. Notably, when the vertex set $V(K)$ is specified, a simplicial complex in $\mathbb{R}^n$ can be fully characterized by listing its simplices. Thus, it can be described purely in terms of combinatorics using abstract simplicial complexes.

**Definition 2.1.4** (Abstract simplicial complex). *[2, Definition 2.3] Consider a finite set $V(K) = \{v_0, \ldots, v_d\}$ for some simplicial complex $K$. An abstract simplicial complex $\tilde{K}$ with vertex set $V(K)$ is a collection of finite subsets of $V(K)$ that satisfies the conditions:*

1. *Every singleton set $\{v_i\}$, where $v_i \in V(K)$, is included in $\tilde{K}$.*
2. *If a set $A$ is in $\tilde{K}$ and $B$ is a subset of $A$, then $B$ must also be in $\tilde{K}$.*

We have a correspondence of $A \subset \tilde{K}$ if and only if $\sigma^{(d)} \in K$ and $B \subset \tilde{K}$ if and only if $\tau^{(k)} \in K$. The abstract simplicial complex $\tilde{K}$ associated with a simplicial complex $K$ is commonly referred to as its vertex scheme. Conversely, if an abstract complex $\tilde{K}$ serves as the vertex scheme for a complex $K$ in $\mathbb{R}^n$, then $K$ is known as a geometric realization of $\tilde{K}$.

**Proposition 2.1.5.** *Every finite abstract simplicial complex $\tilde{K}$ can be geometrically realized in Euclidean space.*

*Proof.* Let $\{v_0, v_1, \ldots, v_d\}$ denote the vertex set of $\tilde{K}$, with $0 \leq d < n$ representing the number of vertices in $\tilde{K}$. Consider $\sigma^{(d-1)} \subset \mathbb{R}^n$, the simplex $[e_1, e_2, \ldots, e_d]$, where $e_i$ represents the $i$-th unit vector. In this context, $K$ refers to the subcomplex of $\sigma^{(d-1)}$ such that $[e_{i_0}, \ldots, e_{i_l}]$ is a $l$-simplex of $K$ with $0 \leq l \leq d$ if and only if $\{v_{i_0}, \ldots, v_{i_l}\}$ is a subset of $\tilde{K}$. □

Pay attention, that this result is in particular interesting for data analysis, as computer aided methods deal with finite point sets. All realizations of an abstract simplicial complex are homeomorphic to each other. The specific realization mentioned above is referred to as the natural realization.

**Proposition 2.1.6.** *Let $\tilde{K}$ be an abstract simplicial complex. Any two geometric realizations $K$ and $K'$ of $\tilde{K}$ are homeomorphic.*

*Proof.* Let $K$ and $K'$ be two geometric realizations of the abstract simplicial complex $\tilde{K}$ in $\mathbb{R}^n$. Let the vertices of $\tilde{K}$ be $\{v_1, \ldots, v_d\}$. Each vertex $v_i$ is mapped to a point $p_i$ in $K$ and to a point $q_i$ in $K'$. There is a bijection of vertices from $K$ to $K'$:

$$\phi : \{p_1, \ldots, p_d\} \to \{q_1, \ldots, q_d\}, \quad \phi(p_i) = q_i \text{ for each } i. \tag{2.7}$$

Extend $\phi$ to a map $f : K \to K'$ on simplices. Let $\sigma^{(k)} = [p_{i_0}, \ldots, p_{i_k}]$ be a $k$–simplex in $K$, corresponding to the simplex $\tau^{(k)} = [q_{i_0}, \ldots, q_{i_k}]$ in $K'$. For $x \in \sigma^{(k)}$, write $x = \sum_{j=0}^{k} \lambda_j p_{i_j}$ with $\sum_{j=0}^{k} \lambda_j = 1$, and define $f(x) = \sum_{j=0}^{k} \lambda_j q_{i_j}$. This defines $f$ linearly on simplices. Since $f$ respects face relations and is continuous on each simplex, $f$ is globally continuous. The map $f$ is bijective, as it maps each vertex $p_i$ to the corresponding $q_i$, and the linear extension preserves this correspondence on simplices. Hence, $f$ is a homeomorphism. Similarly, the inverse map $f^{-1} : K' \to K$ is continuous by the same construction. $\qquad \square$

Furthermore, it has been proven that any finite abstract simplicial complex of dimension $d$ can be realized as a simplicial complex in $\mathbb{R}^{2d+1}$.

**Theorem 2.1.7.** *Any finite abstract simplicial complex of dimension $d$ can be realized as a simplicial complex in $\mathbb{R}^{2d+1}$.*

*Proof.* Let $\tilde{K}$ be a finite abstract simplicial complex of dimension $d$. We construct an injective geometric realization $f : \tilde{K} \to \mathbb{R}^{2d+1}$, $\text{im}(f) = K$. Let $V(K)$ be the vertex set of $\tilde{K}$. First, define an injective map $\tilde{f} : V(K) \to \mathbb{R}^{2d+1}$, which is possible because $V(K)$ is finite and $\mathbb{R}^{2d+1}$ has sufficient dimensionality. If necessary, we adjust $\tilde{f}$ by slight perturbations to ensure that the images of the vertices of each $k$–simplex $\sigma^{(k)} \in \tilde{K}$ are affinely independent in $\mathbb{R}^{2d+1}$. Next, extend $\tilde{f}$ to a map $f : \tilde{K} \to \mathbb{R}^{2d+1}$ by defining it on each subset $\{v_0, \ldots, v_k\}$ through the unique affine map such that $f(v_i) = \tilde{f}(v_i)$. The injectivity of $\tilde{f}$ on the vertices and the affine independence of the vertex images guarantee that $f$ is injective on each simplex and preserves the simplicial structure. In particular, for any $k$–simplices $\sigma^{(k)}, \tau^{(k)} \in K$, we have $f(V(\sigma^{(k)}) \cap V(\tau^{(k)})) = f(V(\sigma^{(k)})) \cap f(V(\tau^{(k)}))$. This ensures that $f$ is a well-defined injective map on the geometric realization of $\tilde{K}$. $\qquad \square$

## 2.2 Simplicial Homology

Given a set $V(\sigma^{(d)})$ representing the vertices of a $d$–simplex $\sigma^{(d)}$, we can establish an orientation for the simplex by selecting a specific ordering of the vertices. If the vertex ordering differs from our chosen order by an odd permutation, the orientation is considered reversed, while even permutations preserve the orientation. Thus, a simplex can have only two possible orientations. To denote the orientation, we use round brackets $(\cdot)$ instead of square brackets $[\cdot]$ for simplices. Moreover, the orientation of a $d$-simplex induces an orientation on its $(d-1)$-faces. Specifically, if $\sigma^{(d)} := (v_0, v_1, \ldots, v_d)$ represents an oriented $d$-simplex, then the orientation of the $(d-1)$-face $\tau_i^{(d-1)}$ of $\sigma^{(d)}$, omitting the vertex $v_i$, is given by

$$\tau_i^{(d-1)} = (-1)^i (v_0, \ldots, v_{i-1}, v_{i+1}, \ldots, v_d). \tag{2.8}$$

**Definition 2.2.1** ($d$-chain). *[28, §2.3] Given a set $\{\sigma_0^{(d)}, \ldots, \sigma_k^{(d)}\}$ of arbitrarily oriented $d$-simplices in a complex $K$ and an abelian group $(G, +)$, a $d$-chain $c$ with coefficients $g_i \in G$ is defined as a formal sum:*

$$c := g_0 \sigma_0^{(d)} + g_1 \sigma_1^{(d)} + \ldots + g_k \sigma_k^{(d)} = \sum_{i=0}^{k} g_i \sigma_i^{(d)}. \tag{2.9}$$

Henceforth, we will assume that $(G, +) = (\mathbb{Z}, +)$.

**Lemma 2.2.2.** *The set of simplicial $d$-chains $C_d^{\triangle}$ forms an abelian group $(C_d^{\triangle}, +)$.*

*Proof.* The identity element of the group is the empty chain, given by:

$$e_{C_d^{\triangle}} = \sum_{i \in \varnothing} g_i \sigma_i^{(d)} := 0. \tag{2.10}$$

The sum of two chains is defined as:

$$c + c' = \sum_{i=0}^{k} g_i \sigma_i^{(d)} + \sum_{j=0}^{l} g_j' \sigma_j^{(d)} \tag{2.11}$$

$$= \sum_{i=0}^{\min(k,l)} (g_i + g_i') \sigma_i^{(d)} + \begin{cases} \sum_{j=l+1}^{k} g_j \sigma_j^{(d)} & \text{if } k > l, \\ 0 & \text{if } k = l, \\ \sum_{j=k+1}^{l} g_j' \sigma_j^{(d)} & \text{if } k < l. \end{cases} \tag{2.12}$$

Hence, $c + c' \in C_d^{\triangle}$. The associativity of the group operation in $C_d^{\triangle}$ follows directly from the associativity of the group operation in $G$. The inverse is

$$c + (-c) = \sum_{i=0}^{k} g_i \sigma_i^{(d)} + \sum_{i=0}^{k} (-g_i) \sigma_i^{(d)} = \sum_{i=0}^{k} (g_i - g_i) \sigma_i^{(d)} = e_{C_d^{\triangle}} = 0. \tag{2.13}$$

$\square$

**Definition 2.2.3** (Boundary). *[18, §2, p.106] Let $\sigma^{(d)}$ be an oriented $d$-simplex in a complex $K$. The boundary of $\sigma^{(d)}$ is defined as the simplicial $(d-1)$-chain of $K$ with coefficients in the abelian group $G$, given by*

$$\partial_d(\sigma^{(d)}) = \sum_{i=0}^{d} (-1)^i \sigma_i^{(d-1)}, \tag{2.14}$$

where $\sigma_i^{(d-1)}$ is a $(d-1)$-face of $\sigma^{(d)}$. If $d = 0$, we define $\partial_0(\sigma^{(0)}) = 0$.

Since $\sigma^{(d)}$ is an oriented simplex, the $\sigma_i^{(d-1)}$ faces also have associated orientations. We extend the definition of the boundary linearly to elements of $C_d^\triangle$.

**Lemma 2.2.4.** *The boundary operator is a group homomorphism*

$$\partial_d : C_d^\triangle \to C_{d-1}^\triangle. \tag{2.15}$$

*Proof.* We define the boundary operator for a $d$-chain $c = \sum_{i=0}^k g_i \sigma_i^{(d)}$:

$$\begin{aligned}
\partial_d(c) &= \sum_{i=0}^k g_i \partial_d(\sigma_i^{(d)}) \\
&= \sum_{i=0}^k g_i \sum_{j=0}^d (-1)^j \sigma_{ij}^{(d-1)} \\
&= \sum_{i=0}^k \sum_{j=0}^d g_i (-1)^j \sigma_{ij}^{(d-1)},
\end{aligned} \tag{2.16}$$

which is an element of $C_{d-1}^\triangle$, where $\sigma_{ij}^{(d-1)}$ are the $(d-1)$-faces of the $d$-simplices $\sigma_i^{(d)}$ in $K$. To verify that $\partial_d$ is a group homomorphism, consider two $d$-chains $c = \sum_{i=0}^k g_i \sigma_i^{(d)}$ and $c' = \sum_{j=0}^l g_j' \sigma_j^{(d)}$. We compute for $l > k$ w.l.o.g.:

$$\begin{aligned}
\partial_d(c + c') &= \partial_d \left( \sum_{i=0}^k g_i \sigma_i^{(d)} + \sum_{j=0}^l g_j' \sigma_j^{(d)} \right) \\
&= \partial_d \left( \sum_{i=0}^k (g_i + g_i') \sigma_i^{(d)} + \sum_{j=k+1}^l g_j' \sigma_j^{(d)} \right) \\
&= \sum_{i=0}^k (g_i + g_i') \partial_d(\sigma_i^{(d)}) + \sum_{j=k+1}^l g_j' \partial_d(\sigma_j^{(d)}) \\
&= \sum_{i=0}^k g_i \partial_d(\sigma_i^{(d)}) + \sum_{j=0}^l g_j' \partial_d(\sigma_j^{(d)}) \\
&= \partial_d(c) + \partial_d(c'). \tag{2.17}
\end{aligned}$$

$\square$

**Remark 2.2.5.** *The compatibility with addition can in principle be seen by the definition of the boundary operator, already.*

**Example 2.2.6.** *Let's consider the 2-simplex $\sigma^{(2)}$ with vertices $v_0, v_1,$ and $v_2$. The 1-faces of this simplex are:*

$$e_0 = (v_1, v_2), \quad \text{connecting } v_1 \text{ and } v_2,$$

$$e_1 = (v_2, v_0), \quad \textit{connecting } v_2 \textit{ and } v_0,$$
$$e_2 = (v_0, v_1), \quad \textit{connecting } v_0 \textit{ and } v_1. \tag{2.18}$$

*Now, let's proceed with the computation:*

$$
\begin{aligned}
\partial_1\big(\partial_2(\sigma^{(2)})\big) &= \partial_1\big(\partial_2(e_0 + e_1 + e_2)\big) \\
&= \partial_1\big(\partial_2(e_0)\big) + \partial_1\big(\partial_2(e_1)\big) + \partial_1\big(\partial_2(e_2)\big) \\
&= \partial_1(v_1, v_2) + \partial_1(v_2, v_0) + \partial_1(v_0, v_1) \\
&= [(v_2) - (v_1)] + [(v_0) - (v_2)] + [(v_1) - (v_0)] \\
&= 0.
\end{aligned} \tag{2.19}
$$

*We observe that $C_0^{\triangle}$ is an abelian group and that oppositely oriented simplices cancel each other out, resulting in:*

$$\partial_1\big(\partial_2(\sigma^{(2)})\big) = 0. \tag{2.20}$$

*This property can be generalized to higher dimensions through induction. Therefore, since $\partial_d$ is a group homomorphism and the chain $c$ is a sum of $d$-simplices, we conclude for $d \geq 1$:*

$$\partial_d^2(c) := \partial_{d-1}\big(\partial_d(c)\big) = 0 \quad \textit{for any } d\textit{-chain } c \textit{ in } C_d^{\triangle}. \tag{2.21}$$

*Consequently, the boundary of the boundary is zero. Moreover, if the boundary of a simplex is zero, it is referred to as a cycle. By this definition, we can deduce that the boundary of any simplex is a cycle.*

**Definition 2.2.7** ($d$-cycle). *[18, §2, p.106] A $d$-chain is called a $d$-cycle if its boundary is equal to zero.*

We denote the set of $d$-cycles of a complex $K$ over an abelian group $G$ – and in particular in this thesis over the group $\mathbb{Z}$ – as $Z_d^{\triangle}$, the simplicial cycle group. It is important to note that $Z_d^{\triangle}$ is a subgroup of $C_d^{\triangle}$ and can also be expressed as:

$$Z_d^{\triangle} := \ker(\partial_d). \tag{2.22}$$

A $d$-cycle of a $k$-complex $K$ is said to be homologous to zero if it can be expressed as the boundary of a $(d+1)$-chain in $K$, where $d = 0, 1, \dots, k-1$. In other words, a cycle is considered a boundary if it can be 'filled in' by a higher-dimensional chain. This equivalence relation is denoted as $c \sim 0$.

**Definition 2.2.8** (Boundary group). *[28, §2.3] The subgroup of $Z_d^{\triangle}$ consisting of boundaries is referred to as the simplicial boundary group $B_d^{\triangle}$.*

It is worth noting that $B_d^{\triangle}$ is equal to the image of the boundary operator $\partial_{d+1}$. Since $B_d^{\triangle}$ is a subgroup of $Z_d^{\triangle}$ and $Z_d^{\triangle}$ is an abelian group, every subgroup of $Z_d^{\triangle}$ is normal. Therefore, we can construct the group quotient $H_d^{\triangle} := Z_d^{\triangle} / B_d^{\triangle}$.

**Definition 2.2.9** (Simplicial homology group). *[18, §2, p.106] The group $H_d^\triangle$ represents the $d$-dimensional simplicial homology group of the complex $K$ over $\mathbb{Z}$. It is expressed as the group quotient:*

$$H_d^\triangle := Z_d^\triangle / B_d^\triangle = \ker(\partial_d)/\operatorname{im}(\partial_{d+1}). \tag{2.23}$$

Next, we want to examine the structure of this homology group by shedding light on its relationship to the connected components of a simplicial complex. We will find that the homology groups of the connected components of the complex, which in turn form a complex themselves, yield the direct sum of the homology group of the entire complex.

**Definition 2.2.10.** *A subcomplex is a subset $S$ of the simplices belonging to a complex $K$, where $S$ itself forms a complex.*

**Definition 2.2.11.** *The collection of all simplices in a complex $K$ with dimension less than or equal to $d$ is referred to as the $d$-skeleton of $K$.*

By definition, the $d$-skeleton forms a subcomplex.

**Definition 2.2.12.** *A complex $K$ is considered connected if it cannot be expressed as the disjoint union of two or more non-empty subcomplexes. A geometric complex is path-connected if there exists a path of 1-simplices connecting any vertex to any other one.*

**Lemma 2.2.13.** *A geometric complex is path-connected if and only if it is connected.*

*Proof.* Suppose $K$ is path-connected. Assume for contradiction that $K$ can be expressed as the disjoint union of two non-empty subcomplexes $L$ and $M$. Since $K$ is path-connected, there exists a path of 1-simplices between any two vertices in $K$. Let $l \in L$ and $m \in M$ be any two vertices. By path-connectedness, there is a path from $l$ to $m$, which contradicts the assumption that $L$ and $M$ are disjoint. Therefore, $K$ is connected.

Suppose $K$ is connected. Pick any vertex $v \in K$. Let $L$ denote the subcomplex of $K$ containing all vertices reachable from $v$ via paths of 1-simplices. If $L \neq K$, then $L$ and $K \setminus L$ form a disjoint union of two non-empty subcomplexes, contradicting the connectedness of $K$. Hence, $L = K$. $\square$

**Proposition 2.2.14.** *[18, Proposition 2.6] Let $K_1, \ldots, K_p$ be the collection of all connected components of a complex $K$. Furthermore, let $H_{d_i}^\triangle$ represent the $d$-th simplicial homology group of $K_i$, and $H_d^\triangle$ denote the $d$-th simplicial homology group of $K$. Then, $H_d^\triangle$ is isomorphic to the direct sum $H_{d_1}^\triangle \oplus \cdots \oplus H_{d_p}^\triangle$.*

*Proof.* Let $C_d^\triangle$ represent the group of simplicial $d$-chains of $K$, and let $K_i$ denote the $i$-th component of $K$. Define $C_{d_i}^\triangle$ as the group of simplicial $d$-chains of $K_i$. It is evident that $C_{d_i}^\triangle$ is a subgroup of $C_d^\triangle$. Furthermore, we observe that $C_d^\triangle$ can be expressed as the direct sum of $C_{d_1}^\triangle, \ldots, C_{d_p}^\triangle$:

$$C_d^\triangle = C_{d_1}^\triangle \oplus \cdots \oplus C_{d_p}^\triangle. \tag{2.24}$$

Our goal is to demonstrate that a similar decomposition applies to the groups $B_d^\triangle$ and $Z_d^\triangle$. By considering $B_{d_i}^\triangle$ as the image of $\partial_{d+1}$ restricted to the subgroup $C_{(d+1)_i}^\triangle$, we can represent the group $B_d^\triangle$ as the direct sum of these restrictions:

$$B_d^\triangle = B_{d_1}^\triangle \oplus \cdots \oplus B_{d_p}^\triangle. \tag{2.25}$$

Thus, for any element $c \in C_{d+1}^\triangle$, we have:

$$c = c_1 + \cdots + c_p, \quad \partial_{d+1}(c) = \partial_{d+1}(c_1) + \cdots + \partial_{d+1}(c_p) \in B_d^\triangle, \tag{2.26}$$

where $c_i \in C_{(d+1)_i}^\triangle$. Let us define $Z_{d_i}^\triangle$ as the intersection of the kernel of $\partial_d$ and $C_{d_i}^\triangle$. It follows that $Z_d^\triangle$ can be expressed as the direct sum of $Z_{d_1}^\triangle, \ldots, Z_{d_p}^\triangle$:

$$Z_d^\triangle = Z_{d_1}^\triangle \oplus \cdots \oplus Z_{d_p}^\triangle. \tag{2.27}$$

To verify this, consider an element $c \in C_d^\triangle$ that belongs to $Z_d^\triangle$. We require $\partial_d(c) = 0$. However, we can express $\partial_d(c)$ as $\partial_d(c_1) + \cdots + \partial_d(c_p)$. Therefore, for $\partial_d(c) = 0$, it must be that $\partial_d(c_i) = 0$, indicating that $c_i \in Z_{d_i}^\triangle$.

Since both $Z_d^\triangle$ and $B_d^\triangle$ can be decomposed componentwise, we conclude that:

$$Z_d^\triangle / B_d^\triangle = Z_{d_1}^\triangle / B_{d_1}^\triangle \oplus \cdots \oplus Z_{d_p}^\triangle / B_{d_p}^\triangle, \tag{2.28}$$

$$H_d^\triangle = H_{d_1}^\triangle \oplus \cdots \oplus H_{d_p}^\triangle. \tag{2.29}$$

$\square$

**Definition 2.2.15** (Index). *The index of a chain* $c = \sum_{i=0}^k g_i \sigma_i^{(d)}$ *is* $I(c) = \sum_{i=0}^k g_i$.

**Remark 2.2.16.** *Sometimes, when we want to emphasise the ground group or the complex used, we write* $H_k^\triangle(K; \mathbb{Z})$ *for* $H_k^\triangle$, $Z_k^\triangle(K; \mathbb{Z})$ *instead of* $Z_k^\triangle$ *or* $B_k^\triangle(K; \mathbb{Z})$ *instead of* $B_k^\triangle$. *We do this in particular when we establish isomorphisms into the ground group or the ground field. However, the coefficients used should be obvious from the context.*

**Proposition 2.2.17.** *[18, Proposition 2.7] If $K$ is a connected complex and $c$ is a 0-chain with $I(c) = 0$, then $I(c) = 0$ is equivalent to $c \sim 0$, where $\sim$ denotes homology equivalence. Furthermore, in this case,* $H_0^\triangle(K; \mathbb{Z}) \cong \mathbb{Z}$.

*Proof.* Let $\sigma^{(1)} = (v_0, v_1)$ be a 1-simplex. For the chain $c = \partial_1(g\sigma^{(1)}) = gv_1 - gv_0$, we have $c \sim 0$, and it is clear that $I(c) = g - g = 0$. Since $I(c + c') = I(c) + I(c')$, $I$ is a group homomorphism. For any 1-chain $c \in C_1^\triangle$ of the form $\sum_{i=0}^k g_i \sigma_i^{(1)}$, where $\sigma_i^{(1)} = (v_i, v_{i+1})$, we have:

$$c = \partial_1(c) \implies c \sim 0 \implies I(c) = I(\partial_1(c)) = 0. \tag{2.30}$$

Consider two vertices $v$ and $w$ in $K$. Since $K$ is connected, there exists a path between $v$ and $w$ consisting of 1-simplices $\sigma_i^{(1)} = (v_i, v_{i+1})$, $i = 0, \ldots, k-1$, where

$v_0 = v$ and $v_k = w$. The boundary of the chain $c = \sum_{i=0}^{k} g\sigma_i^{(1)}$ is given by:

$$\partial_1(c) = \sum_{i=0}^{k} g\partial_1(\sigma_i^{(1)}) = \sum_{i=0}^{k} g\left((v_{i+1}) - (v_i)\right) = gw - gv. \tag{2.31}$$

Since $\partial_1(c)$ is a boundary, we have $c = \partial_1(c) \sim 0$. Thus, $(gw - gv) \sim 0$, implying $gw \sim gv$. Therefore, any 0-chain $c$ in $K$ is homologous to the chain $gv$. We observe that homologous chains have equal indices, i.e., $I(c) = I(gv) = g$. Thus, $c \sim gv \implies c \sim I(c)v$. This shows that if $I(c) = 0$, then $c \sim 0$. Since $I$ is a homomorphism from $C_0^{\triangle} = Z_0^{\triangle}$ to $\mathbb{Z}$, for any 0-simplex $c$ and $g \in \mathbb{Z}$, the chain $gc \in C_0^{\triangle}$ is a cycle with $I(gc) = g$. Therefore, $I(Z_0^{\triangle}) = \mathbb{Z}$. Since $I(c) = 0$ is equivalent to $c \sim 0$, we have $B_0^{\triangle} = \ker(I)$. Thus,

$$H_0^{\triangle} = Z_0^{\triangle}/B_0^{\triangle} = C_0^{\triangle}/\ker(I) \cong \mathbb{Z}. \tag{2.32}$$

$\square$

**Corollary 2.2.18.** *The group $H_d^{\triangle}(K;\mathbb{Z})$ can be represented as $\mathbb{Z}^p = \bigoplus_{i=1}^{p} \mathbb{Z}$, where $p$ denotes the number of connected components present in $K$.*

*Proof.* Proposition 2.2.17 says that if $K$ is connected, then $H_0^{\triangle}(K;\mathbb{Z}) \cong \mathbb{Z}$. Consider a complex $K$ that consists of $p$ connected components $K_1, K_2, \ldots, K_p$. The 0-dimensional simplicial homology group $H_0^{\triangle}(K;\mathbb{Z})$ can be expressed as

$$H_0^{\triangle}(K;\mathbb{Z}) = H_0^{\triangle}(K_1;\mathbb{Z}) \oplus H_0^{\triangle}(K_2;\mathbb{Z}) \oplus \cdots \oplus H_0^{\triangle}(K_p;\mathbb{Z}). \tag{2.33}$$

Since each $K_i$ is a connected component, by Proposition 2.2.17, each $H_0^{\triangle}(K_i;\mathbb{Z})$ is isomorphic to $\mathbb{Z}$. Therefore, we have:

$$H_0^{\triangle}(K;\mathbb{Z}) \cong \bigoplus_{i=1}^{p} \mathbb{Z} = \mathbb{Z}^p. \tag{2.34}$$

$\square$

**Example 2.2.19.**

1. *The zeroth homology group of the circle is isomorphic to $\mathbb{Z}$. Consider a simplicial representation of the circle using four 1-simplices: $\sigma_0^{(1)} = (v_0, v_1)$, $\sigma_1^{(1)} = (v_1, v_2)$, $\sigma_2^{(1)} = (v_2, v_3)$, and $\sigma_3^{(1)} = (v_3, v_0)$. The group $Z_0^{\triangle}$ consists of sums over the four zero-simplices $v_0, v_1, v_2$, and $v_3$ with coefficients in $\mathbb{Z}$. Let $c$ be a zero-chain with non-zero coefficients given by:*

$$c = g_0 v_0 + g_1 v_1 + g_2 v_2 + g_3 v_3. \tag{2.35}$$

*To reduce it to an element of $H_0^{\triangle}$, subtract the chain $c' = g_3 v_3 - g_3 v_2 \sim 0$:*

$$c - c' = g_0 v_0 + g_1 v_1 + (g_2 + g_3) v_2. \tag{2.36}$$

*By repeating this process, we obtain a new chain:*

$$c'' = (g_0 + g_1 + g_2 + g_3)v_0. \tag{2.37}$$

*Since $c'' \sim c$, it represents an element of $H_0^{\triangle}$. Moreover, since $g_i \in \mathbb{Z}$, we can write $g_0 + g_1 + g_2 + g_3 = g$ as $c'' = gv_0$, where $g$ is an element of $\mathbb{Z}$. Therefore, we can choose any $g$, implying that $H_0^{\triangle} \cong \mathbb{Z}$.*

2. *We will demonstrate that $H_d^{\triangle}(S^d) \cong \mathbb{Z}$. The $d$-simplex $\sigma^{(d)}$ and the $d$-ball are homeomorphic, and their boundaries, which consist of $(d-1)$-simplices, are homeomorphic to the $d$-sphere. Thus, the appropriate simplicial structure to impose on $S^d$ is that of the boundary of the $(d+1)$-simplex $\sigma^{(d+1)}$. Let $\{v_0, \ldots, v_d\}$ denote the set of vertices of $\sigma^{(d+1)}$. This set is not oriented, and the orientations of the $(d-1)$-simplices can be arbitrarily determined. We'll utilize their numbering to establish orientations. Consequently, all $d$-chains on this structure are*

$$c = \sum_{i=0}^{d+1} g_i(v_0, \ldots, v_{i-1}, v_{i+1}, \ldots, v_d), \tag{2.38}$$

*where $g_i \in \mathbb{Z}$. Since $\sigma^{(d+1)}$ itself is not part of the structure, there are no boundaries in $Z_d^{\triangle}$, the group of simplicial cycles. Thus, $H_d^{\triangle} = Z_d^{\triangle}/B_d^{\triangle}$ represents the group of simplicial cycles. If $c \in Z_d^{\triangle}$, then $\partial_{d+1}(c) = 0$. Using Eq. 2.38, we have:*

$$\partial_{d+1}(c) = \partial_{d+1}\left(\sum_{i=0}^{d+1} g_i(v_0, \ldots, v_{i-1}, v_{i+1}, \ldots, v_d)\right)$$
$$= \sum_{i=0}^{d+1} g_i \left(\sum_{j=0}^{d+1}(-1)^j(v_0, \ldots, v_{i-1}, v_{i+1}, \ldots, v_{j-1}, v_{j+1}, \ldots, v_d)\right).$$

*By rearranging this sum, we obtain terms of the form:*

$$(g_k - g_l)(v_0, \ldots, v_{j-1}, v_{j+1}, \ldots, v_{i-1}, v_{i+1}, \ldots, v_d), \tag{2.39}$$

*where $k, l = 0, \ldots, d+1$ for all $i, j = 0, \ldots, d$. Each pair of $d$-simplices of $\sigma^{(d+1)}$ intersects along a $(d-1)$-face. Therefore, we obtain terms of the form given in Eq. 2.39 for each of these faces. From this, we deduce that if $\partial_d(c) = 0$, we must have $g_k = g_l$ for all $k, l = 0, \ldots, d+1$. In other words, $g_0 = g_1 = \cdots = g_{d+1}$. Consequently, our original $d$-chain is:*

$$c = \sum_{i=0}^{d+1} g_0(v_0, \ldots, v_{i-1}, v_{i+1}, \ldots, v_d), \tag{2.40}$$

*allowing us to choose $g_0$ from $\mathbb{Z}$. Thus, we conclude that $H_d^{\triangle}(S^d) \cong \mathbb{Z}$.*

3. *We demonstrate that $H_d^{\triangle}(D^d) = 0$. The simplest simplicial structure for $D^d$ is that of the $d$-simplex $\sigma^{(d)}$. Consequently, all $d$-chains can be expressed as $c = g\sigma^{(d)}$, where $g \in \mathbb{Z}$. This form is never a boundary, implying that $H_d^{\triangle} = Z_d^{\triangle}$. However,* $\partial_d(c) = 0$ *is generally only true when* $g = 0$. *Thus,* $H_d^{\triangle}(D^d) \cong 0$.

## 2.3 Singular Homology

In the context of lower dimensions, there is an intuitive understanding of when two topological spaces are fundamentally 'equivalent'. To formalise and strengthen this intuition, various methods have been developed, one of which is the concept of homeomorphism. It would be highly desirable to establish a relation between the homology groups of homeomorphic spaces. Interestingly, it has been found that the homology groups of two topological spaces are isomorphic if they are homeomorphic. This fact requires verification.

In order to achieve this objective, it is necessary to develop a methodology for comparing homology groups. However, it is not immediately evident how this can be accomplished with the tools that have been developed thus far. Indeed, this represents a significant challenge. To address this issue, we propose the introduction of the concept of singular homology. The fundamental principles underlying this concept are analogous to those that have been previously explored.

**Definition 2.3.1** (Singular $d$-simplex). *A singular $d$-simplex in a topological space $X$ is a continuous map* $\tilde{\sigma}^{(d)} : \sigma^{(d)} \to X$.

We define the boundary map $\partial_d$ in a similar manner as before. The singular boundary map, also denoted as $\partial_d$ as we won't need any distinction, is a function that operates on the chain group $C_d(X)$ and maps it to the chain group $C_{d-1}(X)$. It is defined as follows: For any singular $d$-simplex $\tilde{\sigma}^{(d)}$ in $X$, the boundary map $\partial_d(\tilde{\sigma}^{(d)})$ is obtained by summing over all the $(d-1)$-simplices that are obtained by removing one vertex from $\tilde{\sigma}^{(d)}$. Each term in the sum is multiplied by $(-1)^i$, where $i$ represents the index of the removed vertex. In other words, if $v_i$ represents the 0-simplex (vertex) of $\tilde{\sigma}^{(d)}$, then the boundary map can be expressed as:

$$\partial_d(\tilde{\sigma}^{(d)}) = \sum_{i=0}^{d} (-1)^i \tilde{\sigma}^{(d)}|_{[v_0,\ldots,v_{i-1},v_{i+1},\ldots,v_d]}. \tag{2.41}$$

This sum, however, can be a series, as we can take limits on simplices, see Proposition 2.5.3. This does not cause any trouble considering simplicial complexes, as they are considered to be finite, but for CW-complexes we have more flexible definitions, which will be discussed later on. However, the convergence is guaranteed by the definition of simplices. Here, $v_i$ is a map that takes the 0-simplex $\sigma^{(0)}$ to the corresponding vertex in $X$, such that $v_i : \sigma^{(0)} \to X$ is continuous.

As mentioned earlier, when we apply the boundary map twice to a $d$-chain $c$, denoted as $\partial_d^2(c)$ or $\partial_{d-1}(\partial_d(c))$, the result is always zero. This observation leads us to the idea of defining the singular homology groups in a similar way to the simplicial homology groups.

**Definition 2.3.2** (Singular homology). *The singular homology group $H_d(X)$ is*

$$H_d(X) = \ker(\partial_d)/\operatorname{im}(\partial_{d+1}). \tag{2.42}$$

In the following section, we will explore how this definition of homology allows us to establish a simple relationship between homeomorphic spaces and their corresponding homology groups. This relationship becomes apparent when we consider the fact that the definitions of $H_d$ and $H_d^\triangle$ are analogous. Intuitively, we would expect these two groups to be the same for triangulable topological spaces. However, this is not immediately obvious. One reason for this is that $H_d^\triangle$ is finitely generated, while the chain group $C_d(X)$, from which we derived $H_d$, is uncountable. Interestingly, for spaces where both simplicial and singular homology groups can be calculated, these two groups are indeed equivalent. We will provide a proof for this later on. But before we do, let us present some facts about singular homology that support the intuition that $H_d$ is isomorphic to $H_d^\triangle$.

**Remark 2.3.3.** *If we omit information in the notation, e.g. omit the ground group and write $H_k(X)$ instead of $H_k(X; \mathbb{Z})$ or omit the dimension $H(X)$ instead of $H_k(X)$, all results apply to the general case, i.e. to all omitted information.*

**Proposition 2.3.4.** *[18, Proposition 2.6] In the context of a topological space $X$, $H_d(X)$ is isomorphic to the direct sum $H_d(X_1) \oplus \cdots \oplus H_d(X_p)$, where $X_i$ represents the path-connected components of $X$.*

This equivalence serves as the counterpart to Proposition 2.2.17.

*Proof.* As the maps $\tilde{\sigma}^{(d)}$ are continuous, a singular simplex always has a path-connected image in $X$. Consequently, $C_d(X)$ can be expressed as the direct sum of subgroups $C_d(X_1) \oplus \cdots \oplus C_d(X_p)$. The boundary map $\partial_d$ functions as a homomorphism, thereby preserving this decomposition. Therefore, $\ker(\partial_d)$ and $\mathrm{im}(\partial_{d+1})$ also decompose:

$$
\begin{aligned}
H_d(X) &\cong \frac{\ker(\partial_d|_{X_1})}{\mathrm{im}(\partial_{d+1}|_{X_1})} \oplus \frac{\ker(\partial_d|_{X_2})}{\mathrm{im}(\partial_{d+1}|_{X_2})} \oplus \cdots \oplus \frac{\ker(\partial_d|_{X_p})}{\mathrm{im}(\partial_{d+1}|_{X_p})} \\
&= H_d(X_1) \oplus H_d(X_2) \oplus \cdots \oplus H_d(X_p).
\end{aligned}
\tag{2.43}
$$

$\square$

**Definition 2.3.5** (Index). *The index is defined as a map $I : C_d(X) \to \mathbb{Z}$ through*

$$
I(c) = \sum_i g_i \quad \text{for } c = \sum_i g_i \tilde{\sigma}^{(d)} \in C_d(X).
\tag{2.44}
$$

**Proposition 2.3.6.** *The 0-dimensional homology group of a space $X$ is the direct sum of $\mathbb{Z}$ copies, with each copy corresponding to a distinct path-component of $X$.*

*Proof.* To establish the isomorphism $H_0(X) \cong \mathbb{Z}$, it is sufficient to consider the case where $X$ is path-connected. For a 0-chain $c$, the boundary operator $\partial_0(c)$ is always zero since the boundary of any 0-simplex vanishes. Consequently, $\ker(\partial_0) = C_0(X)$, which implies that $H_0(X) = C_0(X)/\mathrm{im}(\partial_1)$. Our goal is to demonstrate that $\ker(I) = \mathrm{im}(\partial_1)$. That is, for any 0-chain $c$, $I(c) = 0$ if and only if $c \sim 0$. Suppose $I(c) = 0$ for $c = \sum_i g_i \tilde{\sigma}^{(0)}$. Since $X$ is path-connected, we can express $c$ as the boundary of a 1-chain $c = \partial_1(\sum_j g_j \tilde{\sigma}^{(1)})$ implying $c \sim 0$. If $c \sim 0$, then $c = \partial_1(\sum_j g_j \tilde{\sigma}^{(1)})$ for some 1-chain $\sum_j g_j \tilde{\sigma}^{(1)}$. By linearity of $I$, we have $I(c) = I(\partial_1(\sum_j g_j \tilde{\sigma}^{(1)})) = 0$. Thus, $I$ induces an isomorphism between $H_0(X)$ and $\mathbb{Z}$ when $X$ is path-connected. For the general case where $X$ has $p$ path-components $X_1, \ldots, X_p$, we apply the above argument to each component. Therefore, $H_0(X)$ is isomorphic to the direct sum $H_0(X) \cong \bigoplus_{i=1}^{p} \mathbb{Z} = \mathbb{Z}^p$. $\qquad\square$

This correspondence serves as the parallel to Corollary 2.2.18.

**Example 2.3.7.** *For any $d \geq 0$, the $d$-th homology group of the $d$-dimensional sphere is isomorphic to $\mathbb{Z}$:*

$$H_d(S^d) \cong \mathbb{Z}.$$

*Proof.* First, we recall the general homology of spheres. The singular homology groups $H_n(S^d)$ of the $d$-dimensional sphere are well-known and are given by:

$$H_n(S^d) \cong \begin{cases} \mathbb{Z}, & n = 0, \\ 0, & 0 < n < d, \\ \mathbb{Z}, & n = d, \\ 0, & n > d. \end{cases}$$

Specifically, the $d$-th homology group $H_d(S^d)$ corresponds to the free abelian group generated by the fundamental class $[\tilde{\sigma}^{(d)}]$, where $\tilde{\sigma}^{(d)}$ is a singular $d$-simplex that maps onto the entire sphere $S^d$. This $d$-simplex represents the whole $d$-sphere and is a cycle because the boundary of $S^d$ is empty. In other words, $\partial_d(\tilde{\sigma}^{(d)}) = 0$, meaning that $[\tilde{\sigma}^{(d)}]$ is a $d$-cycle. Since $S^d$ has no $(d+1)$-dimensional simplices, there are no non-trivial $d$-boundaries, i.e., there are no $(d+1)$-chains to bound our $d$-cycle. Therefore, the class $[\tilde{\sigma}^{(d)}]$ is not a boundary and generates $H_d(S^d)$. $\qquad\square$

## 2.3.1 Singular Chain Complexes

In order to demonstrate the equivalence of $H_d^{\triangle}$ and $H_d$, we will introduce a few concepts that will assist us in our proof.

**Definition 2.3.8** (Chain complex)**.** *A chain complex $\mathcal{C}$ is a sequence of abelian groups $\{C_d\}_{d \in \mathbb{Z}}$ connected by homomorphisms $\partial_d : C_d \to C_{d-1}$ (called boundary operators), such that the composition of any two consecutive maps is zero, i.e., $\partial_{d-1} \circ \partial_d = 0$ for all $d$. Formally, a chain complex is:*

$$\mathcal{C}: \quad \cdots \xrightarrow{\partial_{d+2}} C_{d+1} \xrightarrow{\partial_{d+1}} C_d \xrightarrow{\partial_d} C_{d-1} \xrightarrow{\partial_{d-1}} \cdots \tag{2.45}$$

*where $\partial_{d-1} \circ \partial_d = 0$.*

**Example 2.3.9.** *The groups $C_d(X)$ represent the collection of singular $d$-chains that form a part of a chain complex, where the boundary operator $\partial_d$ is the map between these groups in the respective dimensions:*

$$\mathcal{C}(X): \quad \cdots \xrightarrow{\partial_{d+2}} C_{d+1} \xrightarrow{\partial_{d+1}} C_d \xrightarrow{\partial_d} C_{d-1} \xrightarrow{\partial_{d-1}} \cdots \xrightarrow{\partial_2} C_1 \xrightarrow{\partial_1} C_0 \xrightarrow{\partial_0} 0. \tag{2.46}$$

**Definition 2.3.10.** *A chain map $f$ between two chain complexes $(A, \partial^{(A)})$ and $(B, \partial^{(B)})$ is a collection of maps $f_d : A_d \to B_d$ such that for each $d$, the following condition holds:*

$$\partial_d^{(B)} \circ f_d = f_{d-1} \circ \partial_d^{(A)}. \tag{2.47}$$

*This can be written in the form of a commuting diagram:*

$$
\begin{array}{ccccccccc}
\cdots \xrightarrow{\partial_{d+2}^{(A)}} & A_{d+1} & \xrightarrow{\partial_{d+1}^{(A)}} & A_d & \xrightarrow{\partial_d^{(A)}} & A_{d-1} & \xrightarrow{\partial_{d-1}^{(A)}} & \cdots \\
& \downarrow{f_{d+1}} & & \downarrow{f_d} & & \downarrow{f_{d-1}} & & \\
\cdots \xrightarrow{\partial_{d+2}^{(B)}} & B_{d+1} & \xrightarrow{\partial_{d+1}^{(B)}} & B_d & \xrightarrow{\partial_d^{(B)}} & B_{d-1} & \xrightarrow{\partial_{d-1}^{(B)}} & \cdots .
\end{array}
\tag{2.48}
$$

**Theorem 2.3.11** (Chain map homomorphism). *[27, Theorem 1.3.1] A chain map $f$ between two chain complexes $(A, \partial^{(A)})$ and $(B, \partial^{(B)})$ induces a homomorphism between their respective homology groups $f_\star : H_d(A) \to H_d(B)$.*

*Proof.* By definition of a chain map, we have $f_{d-1} \circ \partial_d^{(A)} = \partial_d^{(B)} \circ f_d$. Let $c \in Z_d(A)$ be a cycle in $A$, i.e., $\partial_d^{(A)}(c) = 0$. Applying $f$ to this equation, we obtain:

$$f_{d-1}(\partial_d^{(A)}(c)) = \partial_d^{(B)}(f_d(c)). \tag{2.49}$$

Since $\partial_d^{(A)}(c) = 0$, we have:

$$f_{d-1}(0) = \partial_d^{(B)}(f_d(c)), \tag{2.50}$$

which implies that $\partial_d^{(B)}(f_d(c)) = 0$. Thus, $f_d(c)$ is a cycle in $B$. Now, let $b \in B_d(A)$ be a boundary in $A$, i.e., there exists $a \in A_{d+1}$ such that $\partial_{d+1}^{(A)}(a) = b$. Applying $f$ to both sides, we obtain:

$$f_d(\partial_{d+1}^{(A)}(a)) = \partial_{d+1}^{(B)}(f_{d+1}(a)). \tag{2.51}$$

Since $\partial_{d+1}^{(A)}(a) = b$, we have:

$$f_d(b) = \partial_{d+1}^{(B)}(f_{d+1}(a)). \tag{2.52}$$

Therefore, $f_d(b)$ is a boundary in $B$. From the above, we see that $f$ maps cycles in $A$ to cycles in $B$ and boundaries in $A$ to boundaries in $B$. Hence, $f$ induces a well-defined map $f_\star : H_d(A) \to H_d(B)$. To show that $f_\star$ is a homomorphism, let $[c_1], [c_2] \in H_d(A)$ be two homology classes. We want to show that $f_\star([c_1] + [c_2]) = f_\star([c_1]) + f_\star([c_2])$. Let $c_1$ and $c_2$ be representatives of $[c_1]$ and $[c_2]$, respectively. Then, $[c_1] + [c_2]$ is represented by $c_1 + c_2$. Applying $f$ to both sides, we have:

$$f_d(c_1 + c_2) = f_d(c_1) + f_d(c_2). \tag{2.53}$$

Since $f_d(c_1)$ and $f_d(c_2)$ are cycles in $B$, we have:

$$[f_d(c_1 + c_2)] = [f_d(c_1)] + [f_d(c_2)]. \tag{2.54}$$

Therefore, $f_\star([c_1] + [c_2]) = f_\star([c_1]) + f_\star([c_2])$. $\qquad\qquad\qquad\square$

### 2.3.2 Exact and Short Exact Sequences

We can apply Theorem 2.3.11 to the case of singular homology. Consider two topological spaces $X$ and $Y$. For any map $f : X \to Y$, we can define an induced homomorphism $f_\star : C_d(X) \to C_d(Y)$ by composing singular $d$-simplices $\tilde{\sigma}^{(d)} : \sigma^{(d)} \to X$ with $f$. Specifically, we have

$$f_\star := f \circ \tilde{\sigma}^{(d)} : \sigma^{(d)} \to Y. \tag{2.55}$$

We can extend this definition by applying $f_\star$ to $d$-chains in $C_d(X)$. This gives us the following commutative diagram:

$$\cdots \xrightarrow{\partial_{d+2}} C_{d+1}(X) \xrightarrow{\partial_{d+1}} C_d(X) \xrightarrow{\partial_d} C_{d-1}(X) \xrightarrow{\partial_{d-1}} \cdots$$
$$\downarrow{f_{d+1}} \qquad\quad \downarrow{f_d} \qquad\quad \downarrow{f_{d-1}} \qquad\qquad \tag{2.56}$$
$$\cdots \xrightarrow{\partial_{d+2}} C_{d+1}(Y) \xrightarrow{\partial_{d+1}} C_d(Y) \xrightarrow{\partial_d} C_{d-1}(Y) \xrightarrow{\partial_{d-1}} \cdots.$$

The chain map $f_d$ gives rise to a homomorphism $f_\star : H_d(X) \to H_d(Y)$. It becomes evident that if $X$ and $Y$ are homeomorphic, then the induced map $f_\star$ is an isomorphism.

Next, we want to formalize the relationships between the homology groups of a topological space $X$, a subset $A \subset X$, and the quotient space $X/A$. To do so, we introduce the concept of exact sequences.

**Definition 2.3.12** (Exact sequences). *An arrangement of elements in the form*

$$\cdots \to A_{d+1} \xrightarrow{\alpha_{d+1}} A_d \xrightarrow{\alpha_d} A_{d-1} \to \cdots \tag{2.57}$$

*is referred to as an exact sequence when the $A_i$ are abelian groups and the $\alpha_i$ are homomorphisms, and it satisfies the condition that $\mathrm{im}(\alpha_{d+1}) = \ker(\alpha_d)$ for all $d$.*

**Remark 2.3.13.**
- *The condition $\mathrm{im}(\alpha_{d+1}) = \ker(\alpha_d)$ implies that $\mathrm{im}(\alpha_{d+1})$ is a subset of $\ker(\alpha_d)$, which is equivalent to $\alpha_d \circ \alpha_{d+1} = 0$. Thus, an exact sequence is a chain complex.*
- *Since $\ker(\alpha_d) \subseteq \mathrm{im}(\alpha_{d+1})$, the homology groups of an exact sequence are trivial.*

**Proposition 2.3.14.** *We can establish the following equivalences:*
1. *$0 \to A \xrightarrow{a} B$ is exact $\iff \ker(a) = 0$, or $a$ is injective.*
2. *$A \xrightarrow{a} B \to 0$ is exact $\iff \mathrm{im}(a) = B$, or $a$ is surjective.*
3. *$0 \to A \xrightarrow{a} B \to 0$ is exact if and only if $a$ is an isomorphism.*

4. *A sequence of the form*

$$0 \to A \xrightarrow{a} B \xrightarrow{b} 0 \tag{2.58}$$

*is said to be exact if and only if the following conditions hold:*
- *The map $a : A \to B$ is injective, meaning that $\ker(a) = 0$.*
- *The map $b : B \to 0$ is surjective, meaning that $\text{im}(b) = 0$.*
- *The kernel of $b$ is equal to the image of $a$, i.e., $\ker(b) = \text{im}(a)$.*
- *If $a : A \hookrightarrow B$ is an inclusion, then $B/\text{im}(a) \cong B/A$.*

*These exact sequences are called short exact sequences.*

*Proof.*

1. Assume $0 \to A \xrightarrow{a} B$ is exact. This means $\text{im}(0) = \ker(a)$, which implies $\ker(a) = 0$ since the image of the zero map is always the trivial group. Therefore, $a$ is injective. Conversely, suppose $\ker(a) = 0$. We need to show $\text{im}(0) = \ker(a)$. Since $\ker(a) = 0$, the only element mapped to the identity in $B$ is the zero element of $A$. Thus, the sequence $0 \to A \xrightarrow{a} B$ is exact.

2. Assume $A \xrightarrow{a} B \to 0$ is exact. This means $\text{im}(a) = \ker(0)$, which implies $\text{im}(a) = B$ since the kernel of the zero map is always the entire group. Therefore, $a$ is surjective. Conversely, suppose $\text{im}(a) = B$. We need to show $\text{im}(a) = \ker(0)$. Since $\text{im}(a) = B$, every element in $B$ has a preimage in $A$ under the map $a$. Thus, the sequence $A \xrightarrow{a} B \to 0$ is exact.

3. Assume $0 \to A \xrightarrow{a} B \to 0$ is exact. From the exactness, we have $\ker(a) = \text{im}(0) = 0$ and $\text{im}(a) = \ker(0) = B$. Therefore, $a$ is both injective and surjective, hence an isomorphism. Conversely, suppose $a$ is an isomorphism, meaning $a$ is both injective (no kernel) and surjective (maps onto $B$). Thus, $0 \to A \xrightarrow{a} B \to 0$ is exact, satisfying $\text{im}(0) = \ker(a)$ and $\text{im}(a) = \ker(0)$.

4. Assume the sequence $0 \to A \xrightarrow{a} B \xrightarrow{b} 0$ is exact. $a$ is injective, and $\text{im}(a) = \ker(b)$. As $b$ is the zero map and surjective, $B/\text{im}(a) = 0$, implying $B = \text{im}(a)$. Therefore, $a$ is an isomorphism, thus $B \cong A$. Conversely, if $a$ is an isomorphism, $\ker(a) = 0$ and $\text{im}(a) = B$.

$\square$

## 2.4 Relative Homology

The concept we will now discuss is that of relative homology groups. Let $X$ be a topological space and $A$ a subspace of $X$. We define $C_d(X, A)$ as the quotient group $C_d(X)/C_d(A)$. This means that chains in $A$ are considered equivalent to the trivial chains in $C_d(X)$. Since the boundary operator $\partial_d : C_d(X) \to C_{d-1}(X)$ also maps $C_d(A)$ to $C_{d-1}(A)$, a natural boundary map on the quotient group $\partial_d : C_d(X, A) \to C_{d-1}(X, A)$ is obtained. This gives rise to the following sequence:

$$\mathcal{C}(X, A) : \quad \cdots \to C_{d+1}(X, A) \xrightarrow{\partial_{d+1}} C_d(X, A) \xrightarrow{\partial_d} C_{d-1}(X, A) \to \cdots . \tag{2.59}$$

This sequence forms a chain complex because $\partial_d \circ \partial_{d+1} = 0$. We can then define the relative homology groups $H_d(X, A)$ as the homology groups of the chain groups within the chain complex and $H_\bullet(\mathcal{C}(X, A))$ as the homology of the entire chain complex, this becomes important in Chapter 3 on persistent homology.

We propose two important facts about $H_d(X, A)$:

**Proposition 2.4.1.** *[18, §2.1, p.115]*
1. *Elements in $H_d(X, A)$ are represented by relative cycles, which are $d$–chains $c$ in $C_d(X)$ such that $\partial_d(c) \in C_{d-1}(A)$.*
2. *A relative cycle $c$ is trivial in $H_d(X, A)$ if and only if it is a relative boundary, i.e., $c$ is the sum of a chain in $C_d(A)$ and the boundary of a chain in $C_{d+1}(X)$.*

*Proof.*
1. Assume $[c] \in H_d(X, A)$ represents a homology class of a chain $c \in C_d(X)$. Since $[c]$ is a class in the relative homology group, $\partial_d(c) \in C_{d-1}(A)$, implying that $c$ is a relative cycle because its boundary maps into $A$. Conversely, if $\partial_d(c) \in C_{d-1}(A)$, then $c$ qualifies as a relative cycle by definition, and any chain homologous to $c$ in $C_d(X)$ that differs from $c$ by a boundary in $C_d(A)$ will also have its boundary in $C_{d-1}(A)$, confirming $[c] \in H_d(X, A)$.

2. For a chain $c$ in $C_d(X)$ to be a relative boundary, it must be expressible as $c = a + \partial_{d+1}(b)$ where $a \in C_d(A)$ and $b \in C_{d+1}(X)$. Applying the boundary operator, we have:

$$\partial_d(c) = \partial_d(a) + \partial_d(\partial_{d+1}(b)) = \partial_d(a) + 0 = \partial_d(a). \tag{2.60}$$

Since $\partial_d(a) \in C_{d-1}(A)$ and $\partial_d \partial_{d+1} = 0$, $c$ is a cycle relative to $A$, making it trivial in $H_d(X, A)$. Conversely, if a relative cycle $c$ is trivial in $H_d(X, A)$, then it must be homologous to a boundary in $C_d(A)$, meaning there exists $a \in C_d(A)$ and $b \in C_{d+1}(X)$ such that $c = a + \partial_{d+1}(b)$. □

Before we further investigate the decomposition of the relative homology groups into exact sequences, we need an intermediate result for commutative diagrams, which facilitates further argumentation. This intermediate result is also known as the Snake Lemma.

**Definition 2.4.2** (Cokernel and coset). *Let $A$ and $B$ be abelian groups (or modules), and let $a : A \to B$ be a homomorphism. The cokernel of $a$, denoted by $coker(a)$, is defined as the quotient group*

$$coker(a) = B/im(a). \tag{2.61}$$

*For any element $y \in B$, $cl(y)$ represents the coset of $y$ in $coker(a)$, defined by*

$$cl(y) = y + im(a). \tag{2.62}$$

This is the set of all elements in $B$ that differ from $y$ by an element of $im(a)$. Or, two elements $y_1, y_2 \in B$ satisfy $cl(y_1) = cl(y_2)$ if and only if $y_1 - y_2 \in im(a)$.

**Lemma 2.4.3** (The Snake Lemma). *[27, Snake Lemma 1.3.2] Consider the following*

*commutative diagram with exact rows:*

$$0 \longrightarrow A' \xrightarrow{f'} A \xrightarrow{f} A'' \longrightarrow 0$$
$$\downarrow{a'} \quad \downarrow{a} \quad \downarrow{a''} \tag{2.63}$$
$$0 \longrightarrow B' \xrightarrow{g'} B \xrightarrow{g} B'' \longrightarrow 0.$$

*There exists an exact sequence:*

$$\ker(a') \xrightarrow{\delta} \ker(a) \xrightarrow{\delta} \ker(a'') \xrightarrow{\delta} \cdots$$
$$\cdots \xrightarrow{\delta} \operatorname{coker}(a') \xrightarrow{\delta} \operatorname{coker}(a) \xrightarrow{\delta} \operatorname{coker}(a'') \tag{2.64}$$

*with the connecting homomorphism $\delta$.*

*Proof.* Let $x'' \in \ker(a'') \subseteq A''$. Since $f$ is surjective, there exists $x \in A$ such that $f(x) = x''$. Since $a''(x'') = 0$, we have $a''(f(x)) = 0$. By the commutativity of the diagram, we get $g(a(x)) = 0$, so $a(x) \in \ker(g) = \operatorname{im}(g')$. Therefore, there exists $y' \in B'$ such that $g'(y') = a(x)$. Define the map $\delta(x'') = \operatorname{cl}(y') \in \operatorname{coker}(a')$, where $\operatorname{cl}(y')$ denotes the equivalence class of $y'$ in $\operatorname{coker}(a')$. We now need to show that $\delta(x'')$ is independent of the choice of $x$. Suppose there is another $x_1 \in A$ such that $f(x_1) = x''$. Then $x - x_1 \in \ker(f)$. Since $f(x - x_1) = 0$, and by commutativity $g(a(x - x_1)) = 0$, we have $a(x - x_1) \in \operatorname{im}(g')$. Thus, there exists $z' \in B'$ such that $g'(z') = a(x - x_1)$. This shows that the difference in the image under $a$ is in the image of $g'$, so $\operatorname{cl}(y') = \operatorname{cl}(y_1')$. Therefore, $\delta(x'')$ is independent of the choice of $x$. Since all maps involved (i.e., $f$, $a$, $g$, and $g'$) are linear, the map $\delta$ is also linear. Now we show that $\ker(\delta) = \operatorname{im}(\ker(a) \to \ker(a''))$. Let $y'' \in \ker(\delta)$, meaning $y'' \in \ker(a'')$ and $\delta(y'') = 0$. By definition of $\delta$, this means that $y' \in B'$ such that $g'(y') = a(x)$ must satisfy $\operatorname{cl}(y') = 0$, i.e., $y' \in \operatorname{im}(a')$. Thus, there exists $w' \in A'$ such that $y' = a'(w')$. This gives $a(x) = g'(a'(w')) = a(g(w'))$, so $x - g(w') \in \ker(a)$, and hence $y'' = f(x - g(w')) \in \operatorname{im}(f)$. Thus, $y'' \in \operatorname{im}(\ker(a) \to \ker(a''))$. This shows that $\ker(\delta) = \operatorname{im}(\ker(a) \to \ker(a''))$. Finally, we consider $y \in \operatorname{coker}(a)$ such that $y \in \ker(\operatorname{coker}(a) \to \operatorname{coker}(a''))$. By definition, this means $y = \operatorname{cl}(z)$ for some $z \in B$ with $a(z) = g'(w')$ for some $w' \in B'$. Since $z \in \ker(g)$, there exists $x \in A$ such that $z = g(x)$. Thus, $y = \operatorname{cl}(g(x)) = 0$. This shows that $\ker(\operatorname{coker}(a) \to \operatorname{coker}(a'')) = 0$. $\qquad\square$

**Theorem 2.4.4.** *[18, p.115f.] The relative homology groups $H_d(X, A)$ are part of the exact sequence:*

$$\cdots \to H_d(A) \to H_d(X) \to H_d(X, A) \to H_{d-1}(A) \to \cdots$$
$$\cdots \to H_0(X, A) \to 0. \tag{2.65}$$

*Proof.* Consider the following diagram:

$$
\begin{array}{ccccccccc}
0 & \longrightarrow & C_d(A) & \overset{i}{\hooklongrightarrow} & C_d(X) & \overset{j}{\twoheadrightarrow} & C_d(X,A) & \longrightarrow & 0 \\
& & \downarrow{\scriptstyle\partial_d} & & \downarrow{\scriptstyle\partial_d} & & \downarrow{\scriptstyle\partial_d} & & \\
0 & \longrightarrow & C_{d-1}(A) & \overset{i}{\hooklongrightarrow} & C_{d-1}(X) & \overset{j}{\twoheadrightarrow} & C_{d-1}(X,A) & \longrightarrow & 0.
\end{array}
\tag{2.66}
$$

Here, $i$ is the inclusion map $C_d(A) \hookrightarrow C_d(X)$, and $j$ is the quotient map $C_d(X) \twoheadrightarrow C_d(X,A)$. Both rows are exact, and the diagram commutes. This gives rise to a long exact sequence of homology groups by the Snake Lemma 2.4.3.

To clarify, consider the chain complexes

$$
\begin{aligned}
\mathcal{A} &= \{C_d(A), \partial_d\}_{d \in \mathbb{N}_0}, \\
\mathcal{B} &= \{C_d(X), \partial_d\}_{d \in \mathbb{N}_0}, \\
\mathcal{C} &= \{C_d(X,A), \partial_d\}_{d \in \mathbb{N}_0}.
\end{aligned}
\tag{2.67}
$$

The horizontal maps $i$ and $j$ are chain maps and thus induce homomorphisms on homology groups $i_\star : H_d(A) \to H_d(X)$ and $j_\star : H_d(X) \to H_d(X,A)$. The exactness of the rows implies that $\ker(j_\star) = \operatorname{im}(i_\star)$. The connecting homomorphism $\partial_d : H_d(X,A) \to H_{d-1}(A)$ can be defined as follows: for any class $[c] \in H_d(X,A)$ represented by a cycle $c \in C_d(X)$ with $\partial_d(c) \in C_{d-1}(A)$, choose a chain $b \in C_d(X)$ such that $j(b) = c$. Since $c$ is a cycle, $\partial_d(b) \in \ker(j) = \operatorname{im}(i)$, so there exists a chain $a \in C_{d-1}(A)$ such that $\partial_d(b) = i(a)$. Define $\partial([c]) = [a] \in H_{d-1}(A)$. We need to verify that the map $\partial_d : H_d(X,A) \to H_{d-1}(A)$ is well-defined: If $b$ and $b'$ both map to $c$ via $j$, then $b - b' \in \ker(j) = \operatorname{im}(i)$. Thus, there exists $a' \in C_{d-1}(A)$ such that $b - b' = i(a')$. Hence,

$$
\partial_d(b) = \partial_d(b') + \partial_d(i(a')) = \partial_d(b') + i(\partial_{d-1}(a')).
\tag{2.68}
$$

So, $[a] = [a']$ in $H_{d-1}(A)$, ensuring the uniqueness of $\partial([c])$. If $c$ and $c'$ are homologous in $H_d(X,A)$, then $c - c' = \partial_d(c'')$ for some $c'' \in C_{d+1}(X,A)$. Thus, $b - b' = \partial_{d+1}(c'')$ in $C_d(X)$, and the same argument as above applies. Thus,

$$
\begin{aligned}
&\cdots \to H_d(A) \to H_d(X) \to H_d(X,A) \to H_{d-1}(A) \to \cdots \\
&\cdots \to H_0(X,A) \to 0
\end{aligned}
\tag{2.69}
$$

is a long exact sequence. $\qquad\square$

We prove now a general result for arbitrary chain complexes, which will apply later on to our situation. Consider the chain complexes $\mathcal{A}, \mathcal{B}, \mathcal{C}$, corresponding to the rows of the following commuting diagram, such that the columns are exakt:

$$
\begin{array}{ccccccccc}
& & & 0 & & 0 & & 0 & \\
& & & \downarrow & & \downarrow & & \downarrow & \\
\mathcal{A}: & \cdots \longrightarrow & A_{d+1} & \xrightarrow{\partial_{d+1}} & A_d & \xrightarrow{\partial_d} & A_{d-1} & \xrightarrow{\partial_{d-1}} & \cdots \\
& & \downarrow i & & \downarrow i & & \downarrow i & & \\
\mathcal{B}: & \cdots \longrightarrow & B_{d+1} & \xrightarrow{\partial_{d+1}} & B_d & \xrightarrow{\partial_d} & B_{d-1} & \xrightarrow{\partial_{d-1}} & \cdots \\
& & \downarrow j & & \downarrow j & & \downarrow j & & \\
\mathcal{C}: & \cdots \longrightarrow & C_{d+1} & \xrightarrow{\partial_{d+1}} & C_d & \xrightarrow{\partial_d} & C_{d-1} & \xrightarrow{\partial_{d-1}} & \cdots \\
& & \downarrow & & \downarrow & & \downarrow & & \\
& & & 0 & & 0 & & 0 &
\end{array}
\tag{2.70}
$$

**Proposition 2.4.5.** *The map $\partial_d : H_d(C) \to H_{d-1}(A)$ is a homomorphism.*

*Proof.* Let $[c_1], [c_2] \in H_d(C)$ be two homology classes, where $c_1, c_2 \in C_d(X, A)$ are cycles. By definition, let $[a_1] = \partial_d([c_1])$ and $[a_2] = \partial_d([c_2])$. This means that $c_1$ and $c_2$ are represented by chains $b_1, b_2 \in C_d(X)$ such that $j(b_1) = c_1$ and $j(b_2) = c_2$. Since $b_1$ and $b_2$ are cycles modulo $A$, we have $\partial_d(b_1) = i(a_1)$ and $\partial_d(b_2) = i(a_2)$. We need to show that $\partial_d([c_1] + [c_2]) = [a_1] + [a_2]$. Consider the sum $c_1 + c_2 \in C_d(X, A)$. By the properties of the quotient map $j$ we yield $j(b_1 + b_2) = j(b_1) + j(b_2) = c_1 + c_2$. Thus, $b_1 + b_2 \in C_d(X)$ is a preimage of $c_1 + c_2$ under the quotient map $j$. Furthermore, the boundary of $b_1 + b_2$ is:

$$
\partial_d(b_1 + b_2) = \partial_d(b_1) + \partial_d(b_2) = i(a_1) + i(a_2) = i(a_1 + a_2). \tag{2.71}
$$

Therefore, the cycle $c_1 + c_2$ maps to the cycle $a_1 + a_2$ under $\partial_d$, which implies on the respective equivalence classes

$$
\partial_d([c_1] + [c_2]) = [a_1 + a_2] = [a_1] + [a_2]. \tag{2.72}
$$

$\square$

Since $i$ and $j$ commute in Eq. 2.70, they induce maps $i_\star$ and $j_\star$ on homology.

**Lemma 2.4.6.** *[18, Theorem 2.16] The given sequence,*

$$
\cdots \to H_d(A) \xrightarrow{i_\star} H_d(B) \xrightarrow{j_\star} H_d(C) \xrightarrow{\partial_d} \cdots
$$
$$
\cdots \xrightarrow{\partial_d} H_{d-1}(A) \xrightarrow{i_\star} H_{d-1}(B) \to \cdots \tag{2.73}
$$

*is exact.*

*Proof.* We need to show that the sequence is exact at each stage:
- We show that $\ker(j_\star) = \operatorname{im}(i_\star)$. Since $j \circ i = 0$, it follows that $j_\star \circ i_\star = 0$, so $\operatorname{im}(i_\star) \subseteq \ker(j_\star)$. Let $[b] \in \ker(j_\star)$, meaning $j_\star([b]) = [c] = 0$ in $H_d(C)$.

This implies that $c$ is a boundary, i.e., there exists a chain $c' \in C_{d+1}$ such that $c = \partial_{d+1}(c')$. Since $j$ is surjective, $c' = j(b')$ for some $b' \in B_{d+1}$. Thus, $j(b) = \partial_{d+1}(c') = \partial_{d+1} \circ j(b')$, which leads to $j(b - \partial_{d+1}(b')) = 0$. Hence, $b - \partial_{d+1}(b') = i(a)$ for an $a \in A_d$. Since $i$ is injective, $a$ is a cycle because:

$$i \circ \partial_d(a) = \partial_d \circ i(a) = \partial_d(b - \partial_{d+1}(b')) = \partial_d(b) = 0, \qquad (2.74)$$

given that $b$ is a cycle. Thus, $i_\star([a]) = [b]$, and $\operatorname{im}(i_\star) = \ker(j_\star)$.

- We need to show that $\ker(\partial_d) = \operatorname{im}(j_\star)$. By definition, if $[b] \in \operatorname{im}(j_\star)$, then $j_\star([b]) = [c] \in H_d(C)$ where $c$ is a cycle. Therefore, $\partial_d([c]) = 0$, and thus $\operatorname{im}(j_\star) \subseteq \ker(\partial_d)$. Let $[c] \in \ker(\partial_d)$. $c$ is a relative cycle such that $\partial_d([c]) = 0$. Hence, there exists $b \in B_d$ such that $j(b) = c$. Therefore, $\ker(\partial_d) \subseteq \operatorname{im}(j_\star)$.
- We show that $\ker(i_\star) = \operatorname{im}(\partial_d)$. If $[c] \in \operatorname{im}(\partial_d)$, then $c$ is a relative boundary. Hence, there exists $b \in B_d$ such that $c = j(b)$. Since $\partial_d(b) = 0$, it follows that $i_\star([c]) = 0$, thus $\operatorname{im}(\partial_d) \subseteq \ker(i_\star)$. Let $[a] \in \ker(i_\star)$, meaning $i_\star([a]) = 0$ in $H_{d-1}(B)$. Thus, $i(a) = \partial_d(b)$ for some $b \in B_d$. Since $\partial_d(j(b)) = j(\partial_d(b)) = j \circ i(a) = 0$, $j(b)$ is a cycle. Therefore, $\partial_d([j(b)]) = [a]$ and $\ker(i_\star) \subseteq \operatorname{im}(\partial_d)$.

$\square$

**Corollary 2.4.7.** *[18, Theorem 2.16] The sequence*

$$\cdots \to H_d(A) \xrightarrow{i_\star} H_d(X) \xrightarrow{j_\star} H_d(X, A) \xrightarrow{\partial_d} \cdots$$
$$\cdots \xrightarrow{\partial_d} H_{d-1}(A) \xrightarrow{i_\star} H_{d-1}(X) \xrightarrow{j_\star} \cdots \to H_0(X, A) \to 0 \qquad (2.75)$$

*is exact.*

*Proof.* The exactness of this sequence is a standard result of algebraic topology, which follows from the long exact sequence of the pair $(X, A)$. The technical argument is analogous to Lemma 2.4.6. The maps in this sequence are:

- $i_\star$ is induced by the inclusion map $i : A \hookrightarrow X$.
- $j_\star$ is induced by the quotient map $j : X \hookrightarrow X/A$, where each chain in $X$ is mapped to its homology class in $X/A$ modulo the image of chains in $A$.
- $\partial_d$ is the boundary homomorphism connecting homology groups, defined by the boundary of a relative cycle in $H_d(X, A)$, which by definition of a chain complex is a cycle in $A$ in one dimension lower.

The exactness at each stage, for example in $H_d(X, A)$, implies that the image of $j_\star$ of $H_d(X)$ is equal to the kernel of $\partial_d$. This shows that any cycle in $H_d(X, A)$ that becomes trivial in $H_{d-1}(A)$ must come from a cycle in $X$ that is not affected by cycles in $A$. The exactness in $H_{d-1}(A)$ further implies that the image of $\partial_d$ is exactly the kernel of $i_\star$ that corresponds to cycles in $A$ that are boundaries in $X$ but not in $A$ itself. The exactness of the entire sequence thus results from the properties of the chain maps and the boundary operators defined in the chain complexes of $A$, $X$ and $X/A$. The additional observation that $\partial_d([c])) = [\partial_d(c)]$ for each relative cycle $c \in H_d(X, A)$ extends the continuity of the exact sequence by the boundary map, since it connects the homology in one dimension in $A$ with the relative homology in $X$ and $A$.

$\square$

## 2.5 Excision

In addition, we refer to the Excision Theorem, a fundamental result in algebraic topology. Simply put, the Excision Theorem allows one to analyse the homology of a space by effectively 'excising' or removing a smaller subspace under certain topological conditions. These conditions usually involve the smaller subspace being 'negligible' in some sense, such as being contained within another subspace. Essentially, this theorem guarantees that the homology of the original space is preserved compared to the homology of the space with the smaller subspace removed. To prove the Excision Theorem, we will need the famous Five Lemma.

**Lemma 2.5.1** (The Five Lemma). *[27, Exercise 1.3.3] Consider a commutative diagram with exact rows:*

$$
\begin{array}{ccccccccc}
A & \xrightarrow{i} & B & \xrightarrow{j} & C & \xrightarrow{k} & D & \xrightarrow{l} & E \\
\downarrow{\alpha} & & \downarrow{\beta} & & \downarrow{\gamma} & & \downarrow{\delta} & & \downarrow{\epsilon} \\
A' & \xrightarrow{i'} & B' & \xrightarrow{j'} & C' & \xrightarrow{k'} & D' & \xrightarrow{l'} & E'.
\end{array}
\tag{2.76}
$$

*If $\alpha, \beta, \delta, \epsilon$ are isomorphisms, then $\gamma$ is also an isomorphism.*

*Proof.* Commutativity of the diagram ensures that $\gamma$ is a homomorphism. We show that $\gamma$ is bijective. Let $c' \in C'$. Since $\delta$ is surjective, there exists $d \in D$ such that $k'(c') = \delta(d)$. Because $\epsilon$ is injective, $\epsilon(l(d)) = l'(\delta(d)) = l'(k'(c')) = 0$, thus $l(d) = 0$. By exactness, $d = k(c)$ for some $c \in C$. Therefore, $k'(c' - \gamma(c)) = 0$, and exactness implies $c' - \gamma(c) = j'(b')$ for some $b' \in B'$. With $\beta$ surjective, $b' = \beta(b)$ for some $b \in B$, hence $\gamma(c + j(b)) = c'$, proving surjectivity of $\gamma$. Assume $\gamma(c) = 0$. Injectivity of $\delta$ and exactness yield $\delta(k(c)) = 0$ and $k(c) = 0$, thus $c = j(b)$ for some $b \in B$. Since $\gamma(j(b)) = j'(\beta(b)) = 0$ and $\beta$ is injective, $b = i(a)$ for some $a \in A$, giving $c = j(i(a)) = 0$. Therefore, $\gamma$ has a trivial kernel and is injective. $\square$

**Theorem 2.5.2** (The Excision Theorem). *[18, Theorem 2.20] Let $Z \subset A \subset X$ such that $\overline{Z} \subseteq \text{int}(A)$, the inclusion $(X \setminus Z, A \setminus Z) \hookrightarrow (X, A)$ induces isomorphisms in homology:*

$$
H_d(X \setminus Z, A \setminus Z) \to H_d(X, A) \quad \text{for all } d.
\tag{2.77}
$$

*Equivalently, for subspaces $A, B \subset X$ whose interiors cover $X$, the inclusion $(B, A \cap B) \hookrightarrow (X, A)$ induces isomorphisms:*

$$
H_d(B, A \cap B) \to H_d(X, A) \quad \text{for all } d.
\tag{2.78}
$$

The equivalence is demonstrated by setting $B = X \setminus Z$ and $Z = X \setminus B$, thereby $A \cap B = A \setminus Z$. The condition $\overline{Z} \subseteq \text{int}(A)$ translates to $X = \text{int}(A) \cup \text{int}(B)$. The proof involves barycentric subdivision, which aids in the computation of homology groups using small singular simplices. In metric spaces, 'smallness' is determined by the diameters of simplices, while in general topological spaces it is defined by their containment within elements of a cover $\mathcal{U}$. Define $\mathcal{U} = \{U_i\}_{i \in I}$, an open cover of $X$ and $I$ an arbitrary index set, and let $C_d^{\mathcal{U}}(X)$ be the subgroup of $C_d(X)$ consisting of chains $\sum_i \alpha_i \tilde{\sigma}_i$ with each $\tilde{\sigma}_i$ contained within one of the $U_i$ and $\alpha_i \in \mathbb{Z}$. The boundary operator $\partial_d : C_d(X) \to C_{d-1}(X)$ preserves this containment, forming a chain complex $C_d^{\mathcal{U}}(X) \to C_{d-1}^{\mathcal{U}}(X)$ with homology groups denoted as $H_d^{\mathcal{U}}(X)$.

**Proposition 2.5.3.** *[18, Proposition 2.21] The inclusion $\iota : C_d^{\mathcal{U}}(X) \hookrightarrow C_d(X)$ is a chain homotopy equivalence. There exists a chain map $\rho : C_d(X) \to C_d^{\mathcal{U}}(X)$ such that the compositions $\iota \circ \rho$ and $\rho \circ \iota$ are chain homotopic to the identity, implying that $\iota$ induces isomorphisms $H_d^{\mathcal{U}}(X) \cong H_d(X)$ for all $d$.*

*Proof. [18, Proof of Proposition 2.21] The proof is divided into four steps, which initially are geometric and become more algebraic over time.*

1. Barycentric subdivision of simplices [18, Proposition 2.21 (1)]: We recall that the points of a $d$-simplex $[v_0, \ldots, v_d]$ are given by linear combinations of the form $\sum_i t_i v_i$ with $\sum_i t_i = 1$ and $t_i \geq 0$ for all $i$.

   **Definition 2.5.4** (Barycenter). *The barycenter of a $d$-simplex $[v_0, \ldots, v_d]$ is the point $b = \sum_i t_i v_i$, where the barycentric coordinates are $t_i = \frac{1}{d+1}$ for all $i$.*

   **Definition 2.5.5** (Barycentric subdivision). *The barycentric subdivision of a $d$-simplex $[v_0, \ldots, v_d]$ is the decomposition into $d$ simplices $[b, w_0, \ldots, w_{d-1}]$. The $(d-1)$-simplex $[w_0, \ldots, w_{d-1}]$ is a face in the barycentric subdivision of $[v_0, \ldots, \hat{v}_i, \ldots, v_d]$.*

   The base case for $d = 0$ is the barycentric subdivision of $[v_0]$, which is simply $[v_0]$ itself. Inductively, it follows that the barycentric subdivision of $[v_0, \ldots, v_d]$ includes the barycenters of all $k$-dimensional faces $[v_{i_0}, \ldots, v_{i_k}]$ of $[v_0, \ldots, v_d]$ for $0 \leq k \leq d$. If $k = 0$, we obtain the original vertices $v_i$, as the barycenter of a $0$-simplex is the $0$-simplex itself. The barycenter of $[v_{i_0}, \ldots, v_{i_k}]$ has coordinates

   $$t_i = \begin{cases} \frac{1}{k+1}, & \text{if } i = i_0, \ldots, i_k, \\ 0, & \text{otherwise.} \end{cases} \tag{2.79}$$

   The $d$-simplices thus form the structure of a simplicial complex. The diameter of a simplex is defined as the maximum distance between its vertices. Consider a point $v$ and a point of the form $\sum_i t_i v_i$ within the simplex $[v_0, \ldots, v_d]$. The distance between these two points satisfies the inequality

   $$\begin{aligned} \left| v - \sum_i t_i v_i \right| &= \left| \sum_i t_i (v - v_i) \right| \\ &\leq \sum_i t_i \left| v - v_i \right| \\ &\leq \sum_i t_i \max_i \left| v - v_i \right| \\ &= \max_i \left| v - v_i \right|. \end{aligned} \tag{2.80}$$

This inequality holds due to the convex combination of the vertices and the properties of the norm.

Let $b_i$ denote the barycenter of the face $[v_0, \ldots, \hat{v}_i, \ldots, v_d]$, where the barycentric coordinates are equal to $\frac{1}{d+1}$ for each vertex except for $v_i$, which has $t_i = 0$. It follows that $b = \frac{1}{d+1} v_i + \frac{d}{d+1} b_i$. The sum of the coefficients is 1, so $b$ lies on the line segment $[v_i, b_i]$, and the distance from $b$ to $v_i$ is $\frac{d}{d+1}$ times the length of $[v_i, b_i]$. Therefore, the distance from $b$ to $v_i$ is bounded by $\frac{d}{d+1}$ times the diameter of $[v_0, \ldots, v_d]$. This bound is independent of the shape of the simplex. If two points $w_j$ and $w_k$ in a simplex are not the barycenter, then these points lie on a proper face of $[v_0, \ldots, v_d]$. Thus, the same bound applies inductively to all 0-simplices. For the repeated $r$-fold barycentric subdivision, we obtain

$$\left( \frac{d}{d+1} \right)^r \xrightarrow{r \to \infty} 0, \tag{2.81}$$

thus, barycentric subdivision can be made with arbitrarily small diameter.

2. Barycentric subdivision of linear chains [18, Proposition 2.21 (2)]: We construct now the subdivision operator $S_d : C_d(X) \to C_d(X)$ and show, that it is chain homotopic to the identity map. For a convex set $Y$ in Euclidean space, the linear maps $\lambda : \tilde{\sigma}^{(d)} \to Y$ generate a subgroup of $C_d(Y)$, which we denote by $LC_d(Y)$, the so-called linear chains. The boundary operator $\partial_d : C_d(Y) \to C_{d-1}(Y)$ maps $LC_d(Y)$ to $LC_{d-1}(Y)$, hence the linear chains form a subcomplex of the singular chain complex of $Y$. We define the unique map $\lambda : \tilde{\sigma}^{(d)} \to Y$ by $[w_0, \ldots, w_d]$, where $w_i$ is the image of the $i$-th vertex of $\tilde{\sigma}^{(d)}$ under $\lambda$. To avoid the constant case distinction regarding the 0-simplices, we complete the chain complex $\mathcal{LC}(Y)$ by defining $LC_{-1}(Y) = \mathbb{Z}$, generated by the equivalence class of the empty simplex $[\varnothing]$, with $\partial_0[w_0] = [\varnothing]$ for all 0-simplices $[w_0]$. Consider the singular $d$-simplex $\tilde{\sigma}^{(d)}$, defined as convex hull of its vertices $\{v_0, \ldots, v_d\}$ in a topological space. A linear map $\lambda : \tilde{\sigma}^{(d)} \to Y$ is determined by the images of the vertices $v_i$, and thus $\lambda$ is entirely defined by the points $w_i = \lambda(v_i)$ in $Y$. The image $\lambda(\tilde{\sigma}^{(d)})$ is a linear simplex in $Y$, and the subgroup $LC_d(Y)$ consists of formal sums of such linear simplices. The boundary of a linear simplex $[w_0, \ldots, w_d]$ is given by the alternating sum

$$\partial_d[w_0, \ldots, w_d] = \sum_{i=0}^{d} (-1)^i [w_0, \ldots, \hat{w}_i, \ldots, w_d], \tag{2.82}$$

where $\hat{w}_i$ denotes the omission of the vertex $w_i$. Since each face of a linear simplex is itself a linear simplex, $\partial_d$ maps $LC_d(Y)$ into $LC_{d-1}(Y)$, ensuring that $\mathcal{LC}(Y)$ is a subcomplex of the singular chain complex $\mathcal{C}(Y)$, for the two discussed chain complexes

$$\mathcal{LC}(Y): \quad \cdots \xrightarrow{\partial_2} LC_1(Y) \xrightarrow{\partial_1} LC_0(Y) \xrightarrow{\partial_0} LC_{-1}(Y) \xrightarrow{\partial_{-1}} 0, \tag{2.83}$$

$$\mathcal{C}(Y): \quad \cdots \xrightarrow{\partial_2} C_1(Y) \xrightarrow{\partial_1} C_0(Y) \xrightarrow{\partial_0} C_{-1}(Y) \xrightarrow{\partial_{-1}} 0. \tag{2.84}$$

The boundary of any 0-simplex $[w_0]$ is defined to be the empty simplex, i.e., $\partial_0[w_0] = [\varnothing]$, for all 0-simplices $[w_0]$. This extension allows us to treat the boundary operator uniformly across all dimensions, including the trivial case where the simplex has no vertices. Each point $b \in Y$ defines a homomorphism $b : LC_d(Y) \to LC_{d+1}(Y)$ on the basis elements, given by $b([w_0, \ldots, w_d]) \mapsto [b, w_0, \ldots, w_d]$. Now, applying $\partial_d$ to $b(\alpha)$, we obtain

$$\partial_{d+1} b([w_0, \ldots, w_d]) = [w_0, \ldots, w_d] - b(\partial_d[w_0, \ldots, w_d]). \tag{2.85}$$

Due to linearity, it follows that

$$\partial_{d+1} b(\alpha) = \alpha - b(\partial_d \alpha) \quad \text{for any } \alpha \in LC_d(Y). \tag{2.86}$$

This shows the equivalence of the equations $\partial_{d+1} b(\alpha) = \alpha - b(\partial_d \alpha)$ and $\partial_{d+1} b + b\partial_d = \text{id}$, meaning that $b$ is a chain homotopy between the identity map and the zero map of the chain complex $\mathcal{LC}(Y)$.

As a consequence, the subdivision homomorphism $S_d : LC_d(Y) \to LC_d(Y)$ can be defined inductively over $d$. Now let $\bar{\sigma}^{(d)}$ be a generator of $LC_d(Y)$, and let $b_\lambda$ be the image of the barycentre of $\bar{\sigma}^{(d)}$ under $\lambda$. We define the value of $S_d(\lambda)$ as $b_\lambda(S_{d-1}\partial_d\lambda)$, such that $b_\lambda : LC_{d-1}(Y) \to LC_d(Y)$ applies. The beginning of the induction shows that $S([\varnothing]) = [\varnothing]$, so that $S$ is the identity on $LC_{-1}(Y)$. Consequently, this mapping is also the identity on $LC_0(Y)$, as the following applies for $d = 0$:

$$S_0([w_0]) = w_0(S_{-1}(\partial_0[w_0])) = w_0(S_{-1}([\varnothing])) = w_0([\varnothing]) = w_0. \tag{2.87}$$

Provided that $\lambda$ is an embedding with $\text{ran}(\lambda) = [w_0, \ldots, w_d]$, $S_d(\lambda)$ can be defined as the sum of the $d$-simplices in the barycentric subdivision of $[w_0, \ldots, w_d]$. Since $S = \mathbb{I}$ applies to $LC_0(Y)$ and $LC_{-1}(Y)$, it follows for all further $d \in \mathbb{N}$ that

$$
\begin{aligned}
\partial_d S_d \lambda &= \partial_d b_\lambda(S_{d-1}\partial_d \lambda) \\
&= S_{d-1}\partial_d \lambda - b_\lambda \partial_{d-1}(S_{d-1}\partial_d \lambda) \quad \text{since } \partial_{d+1} b_\lambda = \mathbb{I} - b_\lambda \partial_d, \\
&= S_{d-1}\partial_d \lambda - b_\lambda S_{d-2}(\partial_d^2 \lambda) \\
&= S_{d-1}\partial_d \lambda. 
\end{aligned} \tag{2.88}
$$

Thus, the chain map $S : \mathcal{LC}(Y) \to \mathcal{LC}(Y)$ satisfies $\partial S = S\partial$.

Next, we define a chain homotopy $T_d : LC_d(Y) \to LC_{d+1}(Y)$ between $S_d$ and the identity, such that the following diagram commutes:

$$
\begin{array}{ccccccccc}
\cdots \to & LC_2(Y) & \to & LC_1(Y) & \to & LC_0(Y) & \to & LC_{-1}(Y) & \to 0 \\
& \downarrow{\scriptstyle S_2} \;{\scriptstyle T_1} & & \downarrow{\scriptstyle S_1} \;{\scriptstyle T_0} & & \downarrow{\scriptstyle S_0 = \mathbb{I}} \;{\scriptstyle T_{-1} = 0} & & \downarrow{\scriptstyle S_{-1} = \mathbb{I}} & \\
\cdots \to & LC_2(Y) & \to & LC_1(Y) & \to & LC_0(Y) & \to & LC_{-1}(Y) & \to 0.
\end{array} \tag{2.89}
$$

We define $T_d$ on $LC_d(Y)$ inductively by setting $T_{-1} = 0$ and $T_d \lambda = b_\lambda(\lambda - T_{d-1}\partial_d\lambda)$ for $d \geq 0$. This formula is an inductive definition of the barycentric subdivision of $\tilde{\sigma}^{(d)} \times \mathbb{I}$ obtained by the union of all simplices in $\tilde{\sigma}^{(d)} \times \{0\} \cup \partial_{d+1}\tilde{\sigma}^{(d+1)} \times \mathbb{I}$ to the barycentre of $\tilde{\sigma}^{(d)} \times \{1\}$. The chain homotopy formula is $\partial_{d+1}T_d + T_{d-1}\partial_d = \mathbb{I} - S_d$. This is trivial for $LC_{-1}(Y)$, where $T_{-1} = 0$ and $S_{-1} = \mathbb{I}$. The formula for $LC_d(Y)$ for all $d \geq 0$ is given by:

$$
\begin{aligned}
\partial_{d+1}T_d\lambda &= \partial_{d+1}b_\lambda(\lambda - T_{d-1}\partial_d\lambda), \\
&= \lambda - T_{d-1}\partial_d\lambda - b_\lambda\partial_d(\lambda - T_{d-1}\partial_d\lambda) \\
&= \lambda - T_{d-1}\partial_d\lambda - b_\lambda[\partial_d\lambda - \partial_d T_{d-1}(\partial_d\lambda)] \\
&= \lambda - T_{d-1}\partial_d\lambda - b_\lambda[S_{d-1}(\partial_d\lambda) + T_{d-2}\partial_d^2(\lambda)] \\
&= \lambda - T_{d-1}\partial_d\lambda - S_d\lambda.
\end{aligned}
\tag{2.90}
$$

Now we can drop the group $LC_{-1}(Y)$ again, and the relation $\partial_{d+1}T_d + T_{d-1}\partial_d = \mathbb{I} - S_d$ still holds since $T_{-1}$ is zero for $LC_{-1}(Y)$.

3. Barycentric subdivision of general chains [18, Proposition 2.21 (3)]: Consider the chain complex $C_d(X)$. Define the operator $S_d : C_d(X) \to C_d(X)$ by $S_d\tilde{\sigma}^{(d)} \coloneqq \sigma_\star S_d\sigma^{(d)}$, where $\tilde{\sigma}^{(d)} : \sigma^{(d)} \to X$ is a singular $d$-simplex, and $S_d\sigma^{(d)}$ denotes the barycentric subdivision of the standard $d$-simplex $\sigma^{(d)}$. The barycentric subdivision $S_d\sigma^{(d)}$ is a signed sum of the $d$-simplices that make up the subdivision of $\sigma^{(d)}$. We will now prove that $S_d$ is a chain map, i.e., that $S_d$ commutes with the boundary operator, $S_{d-1}\partial_d = \partial_d S_d$:

$$
\begin{aligned}
\partial_d S_d\tilde{\sigma}^{(d)} &= \partial_d(\sigma_\star S_d\sigma^{(d)}) = \sigma_\star\partial_d(S_d\sigma^{(d)}) = \sigma_\star S_{d-1}(\partial_d\sigma^{(d)}) \\
&= \sigma_\star S_{d-1}\left(\sum_i (-1)^i \sigma_i^{(d-1)}\right) \\
&= \sum_i (-1)^i \sigma_\star S_{d-1}(\sigma_i^{(d-1)}) \\
&= \sum_i (-1)^i S_{d-1}\left(\tilde{\sigma}^{(d-1)}\big|_{\sigma_i^{(d-1)}}\right) \\
&= S_{d-1}\left(\sum_i (-1)^i \tilde{\sigma}^{(d-1)}\big|_{\sigma_i^{(d-1)}}\right) \\
&= S_{d-1}(\partial_d\tilde{\sigma}^{(d)}).
\end{aligned}
\tag{2.91}
$$

Next, define the operator $T_d$ on $C_d(X)$ by $T_d\tilde{\sigma}^{(d)} \coloneqq \sigma_\star T_d\sigma^{(d)}$, where $\sigma_\star$ denotes the pullback by $T_d\sigma^{(d)}$. We will show that $T_d$ provides a chain homotopy between the subdivision operator $S_d$ and the identity map, satisfying $T_{d-1}\partial_d + \partial_{d+1}T_d = \mathbb{I} - S_d$. The verification is as follows:

$$
\begin{aligned}
\partial_{d+1}T_d\tilde{\sigma}^{(d)} &= \partial_{d+1}(\sigma_\star T_d\sigma^{(d)}) = \sigma_\star(\partial_{d+1}T_d\sigma^{(d)}) \\
&= \sigma_\star(\sigma^{(d)} - S_d\sigma^{(d)} - T_{d-1}\partial_d\sigma^{(d)}) \\
&= \tilde{\sigma}^{(d)} - S_d\tilde{\sigma}^{(d)} - \sigma_\star T_{d-1}\partial_d\sigma^{(d)}
\end{aligned}
$$

$$= \tilde{\sigma}^{(d)} - S_d\tilde{\sigma}^{(d)} - T_{d-1}(\partial_d\tilde{\sigma}^{(d)}). \tag{2.92}$$

4. Iterated barycentric subdivision [18, Proposition 2.21 (4)]: A chain homotopy between the identity map $\mathbb{I}$ and the iterate $S_d^m$ is given by the operator $D_d^m = \sum_{0 \le i < m} T_d S_d^i$. The following calculation verifies this:

$$\begin{aligned}
\partial_{d+1}D_d^m + D_{d-1}^m\partial_d &= \sum_{0 \le i < m}(\partial_{d+1}T_d S_d^i + T_{d-1}S_{d-1}^i\partial_d) \\
&= \sum_{0 \le i < m}(\partial_{d+1}T_d S_d^i + T_{d-1}\partial_d S_d^i) \\
&= \sum_{0 \le i < m}(\partial_{d+1}T_d + T_{d-1}\partial_d)S_d^i \\
&= \sum_{0 \le i < m}(\mathbb{I} - S_d)S_d^i \\
&= \sum_{0 \le i < m}(S_d^i - S_d^{i+1}) \\
&= \mathbb{I} - S_d^m. \tag{2.93}
\end{aligned}$$

In the calculation, we use the fact that $T_{d-1}\partial_d + \partial_{d+1}T_d = \mathbb{I} - S_d$, which shows that $D_d^m$ serves as a homotopy between the identity map $\mathbb{I}$ and $S_d^m$. For each singular $d$-simplex $\tilde{\sigma}^{(d)} : \sigma^{(d)} \to X$, there exists an integer $m$ such that $S_d^m(\tilde{\sigma}^{(d)})$ lies in $C_d^{\mathcal{U}}(X)$, where $\mathcal{U}(X)$ is an open cover of $X$. In this context, $C^{\mathcal{U}}(X)$ denotes the chain complex where simplices are refined to align with the cover $\mathcal{U}$. This holds because, for sufficiently large $m$, the simplices of $S_d^m(\sigma^{(d)})$ have diameters smaller than a Lebesgue number for the cover of $\sigma^{(d)}$ by the open sets $(\tilde{\sigma}^{(d)})^{-1}(\text{int}(U_j))$.[1] We define the operator $D_d : C_d(X) \to C_{d+1}(X)$ by $D_d\tilde{\sigma}^{(d)} := D_d^{m(\tilde{\sigma}^{(d)})}(\tilde{\sigma}^{(d)})$ for each singular $d$-simplex $\tilde{\sigma}^{(d)} : \sigma^{(d)} \to X$. Our objective is to find a chain map $\rho : C_d(X) \to C_d(X)$ whose image lies in $C_d^{\mathcal{U}}(X)$, and which satisfies the chain homotopy equation

$$\partial_{d+1}D_d + D_{d-1}\partial_d = \mathbb{I} - \rho. \tag{2.94}$$

A straightforward approach is to define $\rho$ directly from this equation:

$$\rho = \mathbb{I} - \partial_{d+1}D_d - D_{d-1}\partial_d. \tag{2.95}$$

We can easily verify that $\rho$ is a chain map by calculating:

$$\begin{aligned}
\partial_d\rho(\tilde{\sigma}^{(d)}) &= \partial_d\tilde{\sigma}^{(d)} - \partial_{d+1}^2 D_d\tilde{\sigma}^{(d)} - \partial_d D_{d-1}\partial_d\tilde{\sigma}^{(d)} \\
&= \partial_d\tilde{\sigma}^{(d)} - \partial_d D_{d-1}\partial_d\tilde{\sigma}^{(d)}, \tag{2.96} \\
\rho(\partial_{d+1}\tilde{\sigma}^{(d+1)}) &= \partial_{d+1}\tilde{\sigma}^{(d+1)} - \partial_{d+1}D_d\partial_{d+1}\tilde{\sigma}^{(d+1)} - D_{d-1}\partial_{d+1}^2\tilde{\sigma}^{(d+1)}
\end{aligned}$$

---

[1]The Lebesgue number is a positive real number $\epsilon$ such that any set with diameter less than $\epsilon$ is contained within some element of the cover. The existence of such a number is guaranteed by the compactness of $\sigma^{(d)}$. Since $m$ may vary depending on $\tilde{\sigma}^{(d)}$, we define $m(\tilde{\sigma}^{(d)})$ as the smallest integer such that $S_d^m(\tilde{\sigma}^{(d)})$ lies in $C_d^{\mathcal{U}}(X)$.

$$= \partial_{d+1}\tilde{\sigma}^{(d+1)} - \partial_{d+1}D_d\partial_{d+1}\tilde{\sigma}^{(d+1)}. \tag{2.97}$$

To ensure that $\rho$ maps $C_d(X)$ to $C_d^{\mathcal{U}}(X)$, we compute $\rho(\tilde{\sigma}^{(d)})$:

$$\begin{aligned}
\rho(\tilde{\sigma}^{(d)}) &= \tilde{\sigma}^{(d)} - \partial_{d+1}D_d\tilde{\sigma}^{(d)} - D_{d-1}(\partial_d\tilde{\sigma}^{(d)}) \\
&= \tilde{\sigma}^{(d)} - \partial_{d+1}D_d^{m(\tilde{\sigma}^{(d)})}(\tilde{\sigma}^{(d)}) - D_{d-1}(\partial_d\tilde{\sigma}^{(d)}) \\
&= S_d^{m(\tilde{\sigma}^{(d)})}(\tilde{\sigma}^{(d)}) + D_{d-1}^{m(\tilde{\sigma}^{(d-1)})}(\partial_d\tilde{\sigma}^{(d)}) - D_{d-1}(\partial_d\tilde{\sigma}^{(d)}), \tag{2.98}
\end{aligned}$$

where the last equality follows from $\partial_{d+1}D_d^{m(\tilde{\sigma}^{(d)})} + D_{d-1}^{m(\tilde{\sigma}^{(d)})}\partial_d = \mathbb{I} - S_d^{m(\tilde{\sigma}^{(d)})}$. The term $S_d^{m(\tilde{\sigma}^{(d)})}(\tilde{\sigma}^{(d)})$ lies in $C_d^{\mathcal{U}}(X)$ by the definition of $m(\tilde{\sigma}^{(d)})$. The remaining terms, $D_{d-1}^{m(\tilde{\sigma}^{(d-1)})}(\partial_d\tilde{\sigma}^{(d)}) - D_{d-1}(\partial_d\tilde{\sigma}^{(d)})$, are linear combinations of terms of the form $D_{d-1}^{m(\tilde{\sigma}^{(d-1)})}(\tilde{\sigma}_j^{(d-1)}) - D_{d-1}^{m(\tilde{\sigma}_j^{(d-1)})}(\tilde{\sigma}_j^{(d-1)})$, where $\tilde{\sigma}_j^{(d)}$ denotes the restriction of $\tilde{\sigma}^{(d)}$ to a face of $\sigma^{(d)}$. Since $m(\tilde{\sigma}_j^{(d-1)}) \leq m(\tilde{\sigma}^{(d-1)})$, the difference $D_{d-1}^{m(\tilde{\sigma}^{(d-1)})}(\tilde{\sigma}_j^{(d-1)}) - D_{d-1}^{m(\tilde{\sigma}_j^{(d-1)})}(\tilde{\sigma}_j^{(d-1)})$ consists of terms $T_{d-1}S_{d-1}^i(\tilde{\sigma}_j^{(d-1)})$ with $i \geq m(\tilde{\sigma}_j^{(d-1)})$. These terms lie in $C_d^{\mathcal{U}}(X)$ because the operator $T_{d-1}$ maps $C_{d-1}^{\mathcal{U}}(X)$ into $C_d^{\mathcal{U}}(X)$. Viewing $\rho$ as a chain map from $C_d(X)$ to $C_d^{\mathcal{U}}(X)$, the equation $\partial_{d+1}D_d + D_{d-1}\partial_d = \mathbb{I} - \iota\rho$ holds, where $\iota: C_d^{\mathcal{U}}(X) \hookrightarrow C_d(X)$ is the inclusion. Moreover, $\rho\iota = \mathbb{I}$ since $D_d$ is zero on $C_d^{\mathcal{U}}(X)$. This follows because $m(\tilde{\sigma}^{(d)}) = 0$ for any $\tilde{\sigma}^{(d)} \in C_d^{\mathcal{U}}(X)$, implying that the sum defining $D_d\tilde{\sigma}^{(d)}$ is empty. Consequently, $\rho$ is a chain homotopy inverse to $\iota$.

$\square$

We may now prove Excision, which states that for a space $X = \text{int}(A) \cup \text{int}(B)$, the inclusion map induces an isomorphism $H_d(B, A \cap B) \cong H_d(X, A)$.

*Proof.* Consider the open cover $\mathcal{U} = \{A, B\}$ of $X$, and define the chain group $C_d(A + B)$ as the subgroup of $C_d(X)$ generated by chains in $A$ and chains in $B$, i.e., $C_d(A + B) = C_d(A) + C_d(B)$. The inclusion $C_d(A + B) \subset C_d(X)$ induces a map on the quotient groups $C_d(A + B)/C_d(A) \to C_d(X)/C_d(A)$. Recall that we previously established the identities $\partial_{d+1}D_d + D_{d-1}\partial_d = \mathbb{I} - \iota\rho$ and $\rho\iota = \mathbb{I}$, where $D_d$ is a chain homotopy, $\iota$ denotes the inclusion, and $\rho$ is the associated chain map. These identities persist after passing to the quotient groups, thus the inclusion induces an isomorphism on homology for the chain complex $H_d(\mathcal{C}(A + B)/\mathcal{C}(A)) \cong H_d(\mathcal{C}(X)/\mathcal{C}(A))$. Now, consider the quotient chain group $C_d(B)/C_d(A \cap B)$. It is naturally isomorphic to $C_d(A + B)/C_d(A)$ since both are spanned by singular $d$-simplices in $B$ that are not entirely contained in $A$. Therefore, the inclusion map induces an isomorphism $C_d(B)/C_d(A \cap B) \cong C_d(A + B)/C_d(A)$. Given that these quotient groups are free with bases consisting of singular $d$-simplices in $B$ not contained in $A$, the induced map is indeed an isomorphism. Consequently, the inclusion $(B, A \cap B) \hookrightarrow (X, A)$ induces an isomorphism on homology $H_d(B, A \cap B) \cong H_d(X, A)$.

$\square$

## 2.6   $H_d^\triangle(X) \cong H_d(X)$

We want to establish the equivalence of the homology groups $H_d^\triangle(X)$ and $H_d(X)$. Note that simplicial homology groups $H_d^\triangle(X)$ are defined only for simplicial complexes, whereas singular homology groups $H_d(X)$ can be computed for any topological space. This generality is advantageous, since homeomorphic spaces have isomorphic singular homology groups. To prove the equivalence between $H_d^\triangle(X)$ and $H_d(X)$, we establish an isomorphism between these groups for all dimensions $d$. A homomorphism naturally follows from the map $C_d^\triangle(X) \to C_d(X)$, which sends every $d$-simplex in $X$ to a singular simplex via the mapping $\bar{\sigma}^{(d)} : \sigma^{(d)} \to X$. This induces a homomorphism from $H_d^\triangle(X)$ to $H_d(X)$.

**Theorem 2.6.1** ($H_d^\triangle(X) \cong H_d(X)$). *[18, Theorem 2.27] For each integer $d$, the homomorphisms from the simplicial homology group $H_d^\triangle(X)$ to the singular homology group $H_d(X)$ are isomorphisms. Hence, $H_d^\triangle(X) \cong H_d(X)$ for all $d$.*

*Proof.* Consider a simplicial complex $X$. For each $k$-skeleton $X^k$, the inclusion $X^{k-1} \subset X^k$ leads to the following commutative diagram of exact sequences:

$$H_{d+1}^\triangle(X^k, X^{k-1}) \to H_d^\triangle(X^{k-1}) \to H_d^\triangle(X^k) \to H_d^\triangle(X^k, X^{k-1}) \to H_{d-1}^\triangle(X^{k-1})$$

$$\downarrow \qquad\qquad \downarrow \qquad\qquad \downarrow \qquad\qquad \downarrow \qquad\qquad \downarrow$$

$$H_{d+1}(X^k, X^{k-1}) \longrightarrow H_d(X^{k-1}) \longrightarrow H_d(X^k) \longrightarrow H_d(X^k, X^{k-1}) \longrightarrow H_{d-1}(X^{k-1}).$$

$X^k/X^{k-1}$ contains only simplices of dimension $k$, thus, $C_d(X^k, X^{k-1})$ is trivial for $d \neq k$ and a free abelian group for $d = k$, with a basis formed by the $k$-simplices of $X$. We define a characteristic map $\varsigma : \bigsqcup_i(\sigma_i^{(k)}, \sigma_i^{(k-1)}) \to (X^k, X^{k-1})$, inducing a homeomorphism $\varsigma_* : \bigsqcup_i \sigma_i^{(k)}/\bigsqcup_i \sigma_i^{(k-1)} \to X^k/X^{k-1}$. By Excision 2.5.2 we get

$$H_d^\triangle\left(\bigsqcup_i(\sigma_i^{(k)}, \sigma_i^{(k-1)})\right) \cong H_d(X^k, X^{k-1})$$
$$\cong H_d(X^k/X^{k-1}), \tag{2.99}$$

establishing isomorphisms for $H_d(X^k, X^{k-1})$. Using induction and assuming that the necessary homology isomorphisms hold for dimensions less than $k$, we extend the isomorphism to all dimensions $d$:

$$H_d^\triangle(X^k, X^{k-1}) \cong H_k(X^k, X^{k-1}). \tag{2.100}$$

Up to now, the first and fourth arrow are isomorphisms. By induction we assume the second and fifth arrow to be isomorphisms. The Five Lemma 2.5.1 supports then the isomorphism between $H_k^\triangle(X^k)$ and $H_k(X^k)$, verifying the equivalence of simplicial and singular homology for all $k$. $\qquad\square$

# Chapter 3

# Homological Persistence

This chapter aims to provide an introduction to persistent homology, demonstrating how it can be used to identify multiscale topological features through the application of filtrations to complex data sets represented as point clouds. It discusses barcode isomorphisms and persistence modules that visualise and explain the stability of the detected features within a filtration against perturbations of the data. The Stability Theorem will not be discussed in this thesis [6, Theorem 4.20], as it falls outside the scope of the subject under examination. The chapter concludes with a discussion of persistent chain complexes and the cohomology of these structures, which contribute to a deeper understanding of the evolution of the topology of filtrations on point clouds.

## 3.1 Filtrations of Complexes

The following individual (co)homology groups are defined with coefficients in a fixed field, which we shall denote by $\mathbb{F}$, and form a graded module over the ring of polynomials in one variable $\mathbb{F}[t]$. Replacing the field with a ring leads to substantial problems, but most results still hold for principal ideal domains, see [28, Theorem 2.1, §3.1]. The first step is to define the CW-complex, as every simplicial complex is a CW-complex. In this section, we continue to examine the consequences of absolute and relative homology and cohomology groups for filtrations of CW-complexes, such as the introduced simplicial complexes. In particular, we derive the theory in the context of filtered simplicial complexes upon sets of points, embedded into some metric space. While we can apply the entire theory to simplicial complexes and their (co)homology, CW-complexes provide a much broader context and significantly simplify notation in many instances.

This yields a generalisation of the given situation.

**Definition 3.1.1** (*d*-cell)**.** *A d-cell is a topological space that is homeomorphic to the open d-dimensional ball* $\mathbb{B}^d = \{x \in \mathbb{R}^d \mid \|x\| < 1\}$*, for any norm* $\|\cdot\|$*.*

**Example 3.1.2.**
- *A 0-cell is a single point, homeomorphic to* $\mathbb{B}^0 = \{0\}$*.*
- *A 1-cell is an open line segment, homeomorphic to the open interval* $(-1, 1)$*.*
- *A 2-cell is an open disk, homeomorphic to* $\mathbb{B}^2 = \{(x, y) \mid x^2 + y^2 < 1\}$*.*

- A 3-*cell* is an open ball, homeomorphic to $\mathbb{B}^3 = \{(x, y, z) \mid x^2 + y^2 + z^2 < 1\}$.

**Definition 3.1.3** (Cell complex)**.** *A CW-complex is a Hausdorff topological space $X$ together with a filtration*

$$\mathcal{X}: \quad \varnothing \subset X^0 \subset X^1 \subset X^2 \subset \cdots \subseteq X = \bigcup_{k=0}^{\infty} X^k \coloneqq X^{\infty} \tag{3.1}$$

*such that:*

1. *$X^0$ is a discrete set of points.*
2. *For each $k \geq 1$, $X^k$ is obtained from $X^{k-1}$ by attaching a collection of $k$-cells $\{\sigma_i^{(k)}\}_{i \in I_k}$ via continuous maps $\varphi_i : S^{k-1} \to X^{k-1}$ from the $(k-1)$-sphere, where $I_k$ is an indexing set for the $k$-cells of the CW-complex. This set may be finite or infinite, depending on the number of $k$-cells. Each $i \in I_k$ is an index that uniquely identifies a $k$-cell $\sigma_i^{(k)}$ in the set of all $k$-cells of the complex. Specifically,*

$$X^k = X^{k-1} \cup \bigcup_{i \in I_k} \sigma_i^{(k)}, \tag{3.2}$$

   *where each $\sigma_i^{(k)} \cong \mathbb{B}^k$ and $\sigma_i^{(k)} \cap X^{k-1} = \varphi_i(S^{k-1})$.*

3. *The topology on $X$ is the weak topology with respect to the cells $\{\sigma_i^{(k)}\}$.*

**Remark 3.1.4.** *A CW-complex satisfies the following conditions:*

- *The closure $\overline{\sigma_i^{(m)}}$ intersects finitely many cells $\sigma_j^{(n)}$ for $n \leq m$, and $n, m, i, j \in \mathbb{N}_0$.*
- *A set $A \subseteq X$ is closed if and only if $A \cap X^k$ is closed in $X^k$ for all $k \in \mathbb{N}_0$.*

**Proposition 3.1.5.** *A simplicial complex is a CW-complex.*

*Proof.* Let $K$ be a simplicial complex. We construct $K$ as a CW-complex by

$$\mathcal{K}: \quad \varnothing \subset K^0 \subset K^1 \subset K^2 \subset \cdots \subseteq K = \bigcup_{k=0}^{\infty} K^k, \tag{3.3}$$

where $K^k$ is the $k$-skeleton of $K$. The 0-skeleton $K^0$ consists of the 0-simplices, which form a discrete set. Assume $K^{k-1}$ is constructed. To form $K^k$, we attach each $k$-simplex $\sigma^{(k)} \in K$ to $K^{k-1}$ via a map $\varphi_{\sigma^{(k)}} : \partial_k \sigma^{(k)} \to K^{k-1}$, where $\partial_k \sigma^{(k)}$ is the boundary of $\sigma^{(k)}$. Each $\sigma^{(k)}$ is homeomorphic to an open $k$-ball $\mathbb{B}^k$, and $\partial_k \sigma^{(k)}$ is homeomorphic to the $(k-1)$-sphere $S^{k-1}$. The topology on $K$, the geometric realization of the vertex scheme $\tilde{K}$, is the weak topology with respect to the simplices, thereby satisfying the topology condition of a CW-complex. The closure-finite condition is satisfied because the closure of each simplex intersects only finitely many simplices. The proof is completed by induction over $k$. $\qquad \square$

We investigate the persistent topology of filtered topological spaces, focusing primarily on the prototypical example of a filtered CW-complex $\mathcal{X} : \varnothing \subset X^0 \subset X^1 \subset X^2 \subset \cdots \subseteq X = \bigcup_{k=1}^{\infty} X^k \coloneqq X^{\infty}$, where $X^0$ starts with vertices $\sigma^{(0)}$, and each subsequent complex $X^i$ is constructed by adding a single cell to the previous complex $X^k := X^{k-1} \cup \sigma^{(k)}$. Here, the indexing set is $\{0, 1, 2, \ldots, k, \ldots\}$, denoting the filtration level. Additionally, associated real values $a_i$ are assigned to the indexing sets, satisfying $a_0 \leq a_1 \leq a_2 \leq \cdots \leq a_k$, to indicate for which parameter the particular simplex is created during some filtration process.

**Example 3.1.6.** *[9, §2.2, Example] The following illustrative example is a cellular filtration of the 2-sphere, denoted by $\mathcal{S}^2$. This process constructs a CW-complex and introduces an ordering among the cells to provide a useful framework for organising the cells. In this sequence, the cell $\sigma_0^{(0)}$ stands for the initial point. The set $S^1$ consists of two distinct points. In $S^2$ a path connecting $\sigma_0^{(0)}$ and $\sigma_1^{(0)}$ is added. $S^3$ extends this path with its reversal, from $\sigma_1^{(0)}$ to $\sigma_0^{(0)}$, thus distinguishing it from the path in $S^2$. Finally, $S^4$ and $S^5$ make a distinction between cells representing the upper and lower halves of the sphere, respectively.*

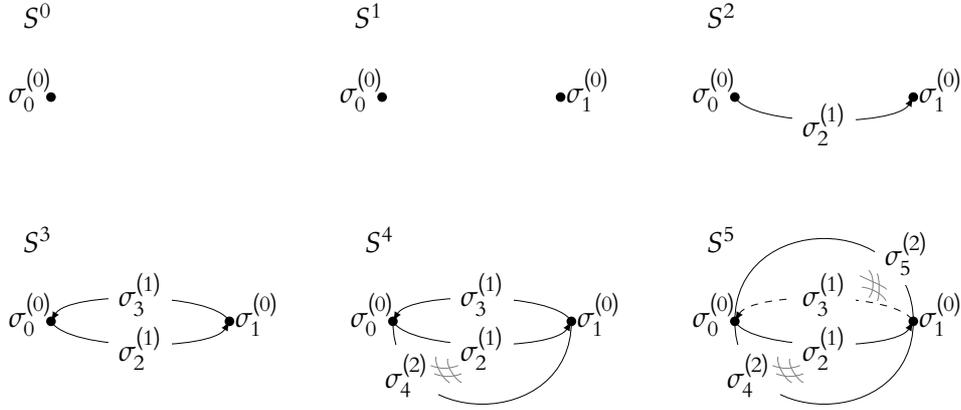

$\mathcal{S}^2 : \varnothing$
$\subset S^0 := \{\sigma_0^{(0)}\}$
$\subset S^1 := \{\sigma_0^{(0)}, \sigma_1^{(0)}\}$
$\subset S^2 := \{\sigma_0^{(0)}, \sigma_1^{(0)}, \sigma_2^{(1)} := (\sigma_0^{(0)}, \sigma_1^{(0)})\}$
$\subset S^3 := \{\sigma_0^{(0)}, \sigma_1^{(0)}, \sigma_2^{(1)}, \sigma_3^{(1)} := (\sigma_1^{(0)}, \sigma_0^{(0)})\}$
$\subset S^4 := \{\sigma_0^{(0)}, \sigma_1^{(0)}, \sigma_2^{(1)}, \sigma_3^{(1)}, \sigma_4^{(2)} := (\sigma_2^{(1)}, \sigma_3^{(1)})\}$
$\subset S^5 := \{\sigma_0^{(0)}, \sigma_1^{(0)}, \sigma_2^{(1)}, \sigma_3^{(1)}, \sigma_4^{(2)}, \sigma_5^{(2)} := (\sigma_3^{(1)}, \sigma_2^{(1)})\}.$

In this example we have six cells $\sigma_0^{(0)}, \sigma_1^{(0)}, \ldots, \sigma_5^{(2)}$ with birth times $a_i = i$ for $i = 0, \ldots, 5$.

## 3.2   Persistence Modules

Persistent homology is derived by applying the homology functor to a filtration of a topological space. This section will examine the concept of persistent modules, which serve to describe persistent homology from an algebraic perspective.

Similarly, as filtrations describe the evolution of a simplicial complex under a one-parameter or multi-parameter family, the persistence module describes the evolution of associated vector spaces. In the case of homology vector spaces or modules, the dimension or rank of the associated homomorphisms are the $k$-th Betti numbers, which indicates the number of $k$-dimensional holes of a topological space in intuitive terms.

**Definition 3.2.1** (Persistence module). *[6, §1.1] A real persistence module $C$ is a collection of vector spaces $\{C^i\}_{i \in I}$ indexed over a partially ordered set $(I, \leq)$ together with linear maps*

$$\phi_i^j : C^i \to C^j, \quad \forall i, j \in I : i \leq j, \tag{3.4}$$

*such that $\phi_i^k = \phi_j^k \circ \phi_i^j$ and $\phi_i^i = \mathbb{I}_{C^i}$ for all $i \leq j \leq k$.*

We rephrase this definition in categorical language, which makes it easier for us to define the normal form of persistence modules. Persistence modules have a decomposition into interval modules, which cannot be further decomposed. On the other hand, an interval module is a persistence module that is indecomposable. Let $\mathrm{Vect}_{\mathbb{F}}$ be the category of finite dimensional vector spaces over $\mathbb{F}$. The reals are a partially ordered category $(\mathbb{R}, \leq)$. Consider the real persistence module:

**Definition 3.2.2** (Persistence module). *[6, §1.3] A real persistence module is a functor $C : (\mathbb{R}, \leq) \to \mathrm{Vect}_{\mathbb{F}}$. It is uniquely defined by a family $\{C^i\}_{i \in \mathbb{R}}$ of finite dimensional vector spaces over $\mathbb{F}$ together with morphisms $\phi_i^j : C^i \to C^j, i \leq j$, such that the following diagram commutes:*

$$\begin{array}{ccc} C^i & \xrightarrow{\phi_i^k} & C^k \\ {\scriptstyle \phi_i^j} \downarrow & \nearrow {\scriptstyle \phi_j^k} & \\ C^j. & & \end{array} \tag{3.5}$$

**Remark 3.2.3.** *A morphism of persistence modules is a natural transformation. This can be written as a morphism $\varphi_i : C^i \to C'^i$, which is compatible with the maps $\phi_i^j$. This defines the category* Pers *of persistence modules. Note that the upper index in this particular case is not a power, but nothing else than an index, induced by the partial order of the domain of the persistence module functor.*

The category of persistence modules is abelian. A notable consequence of the aforementioned results is that a part of the constructions for persistence modules remain valid when the category of finite-dimensional vector spaces over a field is replaced by another abelian category. This includes, for instance, (in)finite dimensional vector spaces over other fields. It is also possible to consider persistence modules over (co)chain complexes in place of vector spaces.

In the following, we will make a number of simplifying assumptions:

**Assumption 3.2.4.**
  *1. For all $i \in \mathbb{R} \setminus J$, where $J$ is a finite subset, there exists a neighbourhood $U_i$ of $i$ such that the map $\phi_i^{i'}$ is an isomorphism for all $i \leq i'$ in $U_i$.*

2. $C^i = 0$ for $i$ sufficiently small.
3. For all $j \in \mathbb{R}$ and for all $i \leq j$ with $j - i$ sufficiently small, the mapping $\phi_i^j$ is an isomorphism.

The assumptions 1, 2 are referred to as finite type conditions. Assumption 1 implies that there is a finite number of 'jumps', while we have that $C^i = C^\infty$ for $i$ large enough. Assumption 3 is called semi-continuity and can be visualised by the following indecomposable interval modules:

**Definition 3.2.5** (Interval module). *Let $I$ be an interval of the form $[a, b)$ or $[a, \infty)$ with $-\infty < a < b \leq \infty$. We define an interval module $\mathbb{F}(I^i)$ over a field $\mathbb{F}$ as follows:*

$$\mathbb{F}(I^i) = \begin{cases} \mathbb{F}, & \text{if } i \in I, \\ 0, & \text{otherwise.} \end{cases} \tag{3.6}$$

*The structure maps $f_i^j : \mathbb{F}(I^i) \to \mathbb{F}(I^j)$ for $i \leq j$ are defined as:*

$$f_i^j = \begin{cases} \mathrm{id}_{\mathbb{F}}, & \text{if } i, j \in I, \\ 0, & \text{if } i \notin I \text{ or } j \notin I. \end{cases} \tag{3.7}$$

### 3.2.1 Decomposition of Persistence Modules

The following theorem describes the decomposition of any persistence module:

**Theorem 3.2.6** (Decomposition). *[6, §1.5] Every persistence module $\mathcal{C}$ over a field $\mathbb{F}$ is isomorphic to a direct sum of interval modules:*

$$\mathcal{C} \cong \bigoplus_{i=0}^{d} \mathbb{F}(I^i)^{m_i}, \tag{3.8}$$

*where $I^i \neq I^j$ for all $i \neq j$, and $m_i$ are multiplicities. This decomposition is unique up to isomorphism and permutation of the summands.*

*Proof.* To establish the existence of such a decomposition, we define a functor from the category Pers of persistence modules to the category of finitely generated graded modules over the polynomial ring $\mathbb{F}[i]$. For a persistence module $\mathcal{C} = \{C^i, f_i^j\}_{i \leq j}$, we map $\mathcal{C}$ to the graded $\mathbb{F}[i]$-module $M = \bigoplus_{i \in \mathbb{R}} C^i$, where the action of $i$ on $M$ is defined by $i \cdot (v_0, v_1, \ldots) \coloneqq (0, \ldots, f_i^{i+1}(v_i), f_{i+1}^{i+2}(v_{i+1}), \ldots)$. The conditions ensuring the persistence module structure (finite generation and compatible transition maps) imply that $M$ is finitely generated as an $\mathbb{F}[i]$-module. This functor is an equivalence of categories: the inverse functor takes a graded $\mathbb{F}[i]$-module and interprets the grading as a persistence parameter over $\mathbb{R}$. Since $\mathbb{F}[i]$ is a principal ideal domain (PID), we can apply the structure theorem for finitely generated modules over a PID [28, §2.1], which yields:

$$M = \bigoplus_{i \in \mathbb{R}} C^i \cong \bigoplus_{k=0}^{d} T^{a_k} \mathbb{F}[i] \oplus \bigoplus_{l} T^{b_l} \mathbb{F}[i]/(i^{d_i}). \tag{3.9}$$

Here, $T^{a_k}\mathbb{F}[i]$ denotes the graded free module over $\mathbb{F}[i]$ shifted by $a_k$, meaning $(T^{a_k}\mathbb{F}[i])_d = \mathbb{F}[i]_{d-a_k}$. This corresponds to interval modules of the form $[a_k, \infty)$. Similarly, $T^{b_l}\mathbb{F}[i]/(i^{d_l})$ represents a torsion module over $\mathbb{F}[i]$ that is shifted by $b_l$, defined by $(T^{b_l}\mathbb{F}[i]/(i^{d_l}))_d = \mathbb{F}[i]_{d-b_l}/(i^{d_l})$. This corresponds to interval modules of the form $[b_l, b_l + d_l)$.

For uniqueness, consider the endomorphism ring $\text{End}(\mathbb{F}(I))$. We show that $\text{End}(\mathbb{F}(I)) \cong \mathbb{F}$. Any endomorphism of interval modules must respect the persistence module structure, meaning $\phi_i : \mathbb{F}(I^i) \to \mathbb{F}(I^j)$ must be a linear map. Since $\mathbb{F}(I^i) = \mathbb{F}$ for $i \in I$ and $0$ otherwise, $\phi_i$ must be multiplication by a scalar $\lambda_i \in \mathbb{F}$. The critical technical reason for this is that $\mathbb{F}(I^i)$ is a simple module (i.e., it has no non-trivial submodules) when $i \in I$. Furthermore, for the endomorphism to be compatible with the structure maps $f_i^j$, we must have $\lambda_i = \lambda_j$ for all $i, j \in I$, meaning $\phi$ acts as multiplication by a constant $\lambda \in \mathbb{F}$ uniformly across the entire interval $I$. Hence, $\text{End}(\mathbb{F}(I)) \cong \mathbb{F}$.

Now, suppose we have two decompositions:

$$\bigoplus_i \mathbb{F}(I^i) \cong \mathcal{C} \cong \mathcal{C}' \cong \bigoplus_j \mathbb{F}(J^j). \tag{3.10}$$

Consider the composition of morphisms:

$$f_i^j : \mathbb{F}(I^i) \hookrightarrow \mathcal{C} \cong \mathcal{C}' \twoheadrightarrow \mathbb{F}(J^j), \tag{3.11}$$

$$g_j^i : \mathbb{F}(J^j) \hookrightarrow \mathcal{C}' \cong \mathcal{C} \twoheadrightarrow \mathbb{F}(I^i). \tag{3.12}$$

The equation $\sum_{i,j} g_j^i f_i^j = 1$ implies that at least one of the compositions $f_i^j g_j^i \neq 0$, yielding an isomorphism between $\mathbb{F}(I^i)$ and $\mathbb{F}(J^j)$ for some $i, j$. Since the interval modules are indecomposable and distinct intervals $I^i$ and $J^j$ cannot yield non-trivial morphisms between each other, this forces $I^i = J^j$ and $m_i = m_j$.

By induction on the number of summands $i, j$, we conclude that the multiplicities $m_i$ and the intervals $I_i$ in both decompositions must coincide. $\square$

**Definition 3.2.7** (Barcode). *The barcode associated with a persistence module $\mathcal{C}$ is a collection of tuples consisting of interval modules and their corresponding multiplicities, given by the decomposition in Theorem 3.2.6:*

$$\mathcal{B}(\mathcal{C}) \coloneqq \{(I^i, m_i)\}_{i=0}^d. \tag{3.13}$$

*Here, each $I^i$ is an interval, and $m_i$ are the multiplicities of the interval modules $\mathbb{F}(I^i)$ in the decomposition of $\mathcal{C}$, respectively.*

At this point, we can establish a connection to spectral theory. In particular, for real-valued filtrations, the spectrum of a persistence module $\mathcal{C}$ can be represented by specific intervals arising from the decomposition of persistence modules.

**Definition 3.2.8** (Spectrum). *A point $k \in \mathbb{R}$ is called a spectral point of a persistence module $\mathcal{C}$ if, for every neighbourhood $U_k$ of $k$, there exist $i < j$ in $U_k$ such that the map $f_i^j$* *is not an isomorphism. The (finite) spectrum of $\mathcal{C}$ is defined as:*

$$Spec(\mathcal{C}) := \{spectral\ points\} \cup \{\infty\}. \tag{3.14}$$

**Remark 3.2.9.**
- *For each interval $I = [a, b)$ or $[a, \infty)$, and for any two points $i < j$ within this interval, the corresponding structure map $f_i^j$ of the persistence module $\mathcal{C}$ is an isomorphism. A point $k$ becomes a spectral point if, within any neighbourhood $U_k$ around $k$, there exists a pair $(i, j)$ with $i, j \in U_k$ such that the map $f_i^j$ is not an isomorphism. This identifies a change in the module's behavior, such as the birth of a new feature or the death of an existing one.*
- *The barcode $\mathcal{B}(\mathcal{C})$ provides a visual and combinatorial representation of the intervals and their multiplicities. Each interval in the barcode corresponds to a range of points $k$ on the real line where the persistence module maintains consistent behavior. A spectral point indicates a transition between these behaviors, which is encoded withi the endpoints of the intervals in the barcode.*
- *The spectrum $Spec(\mathcal{C})$ summarizes the critical values on the real line where the module's structure changes. This concept is similar to the spectrum in functional analysis, where one considers the set of scalars (like eigenvalues) that describe how a linear operator acts on a vector space.*
- *Since the persistence module $\mathcal{C}$ is finitely generated, the set of intervals in its decomposition is finite. This implies that the set of spectral points, which correspond to transitions between intervals, is also finite. The inclusion of $\infty$ in $Spec(\mathcal{C})$ represents the asymptotic behavior of $\mathcal{C}$, capturing the intervals of the form $[a, \infty)$.*

## 3.2.2 Persistent Homology on Complexes

In the context of algebraic topology, applying the homology functor $H(-)$ to a finite filtration of a topological space $\mathcal{X}$ yields a sequence of algebraic structures:

$$H(\mathcal{X}): \quad H(X^0) \to H(X^1) \to H(X^2) \to \cdots \to H(X^d) \to H(X), \tag{3.15}$$

where $H(-)$ generally represents either the $k$-th dimensional homology, denoted by $H_k(-; \mathbb{F})$, or the absolute homology, expressed as $H_\bullet(-; \mathbb{F})$. Here, this diagram characterizes a sequence of abelian groups or finite-dimensional vector spaces, connected through linear maps, forming a persistence module. Persistence modules are central to understanding how the features of a space evolve over time. According to Theorem 3.2.6, they can be decomposed into interval modules. Each interval module is associated with an ordered pair of integers $(p, q)$ where $0 \leq p \leq q \leq d$, within a finite filtration. These pairs $(p, q)$ signify topological features that persist over an index set $I := \{p, \ldots, q\}$, where $\inf\{I\} = p$ and $\sup\{I\} = q$. Conventionally, these tuples are interpreted as half-open intervals $[a_p, a_{q+1})$, with $a_{d+1} = \infty$ being a customary notation when the sequence extends beyond the largest indexed space. The decomposition of a persistence module into its constituent interval modules is represented in a persistence diagram or a barcode. This barcode is a multiset of ordered tuples $(p, q)$ or, alternatively, a multiset of half-open intervals $[a_p, a_{q+1})$.

$$H_\bullet(\mathcal{X}): \quad H_\bullet(X^0) \to H_\bullet(X^1) \to \cdots \to H_\bullet(X^{d-1}) \to H_\bullet(X^d),$$

$$H^\bullet(\mathcal{X}): \quad H^\bullet(X^0) \leftarrow H^\bullet(X^1) \leftarrow \cdots \leftarrow H^\bullet(X^{d-1}) \leftarrow H^\bullet(X^d),$$

$$H_\bullet(X^\infty, \mathcal{X}): \quad H_\bullet(X^d) \to H_\bullet(X^d, X^0) \to H_\bullet(X^d, X^1) \to \cdots \to H_\bullet(X^d, X^{d-1}),$$

$$H^\bullet(X^\infty, \mathcal{X}): \quad H^\bullet(X^d) \leftarrow H^\bullet(X^d, X^0) \leftarrow H^\bullet(X^d, X^1) \leftarrow \cdots \leftarrow H^\bullet(X^d, X^{d-1}).$$

Figure 3.2: The four standard persistence modules are from top to bottom: persistent absolute homology, persistent absolute cohomology, persistent relative homology and persistent relative cohomology. The ($\bullet$) indicates the direct sum of (co)homology over dimensions, in other words the absolute (co)homology. We can also write $H^\bullet(X^d, \mathcal{X})$, $H_\bullet(X^d, \mathcal{X})$ or $H^\bullet(X, \mathcal{X})$, $H_\bullet(X, \mathcal{X})$ instead of $H^\bullet(X^\infty, \mathcal{X})$, $H_\bullet(X^\infty, \mathcal{X})$, as the $k$-th (co)homology groups for $k > d$ become trivial.

This collection is formally expressed through the forgetful functor $\mathrm{Pers}(-)$, called the persistence diagram:

$$\mathrm{Pers}(H_k(\mathcal{X}; \mathbb{F})) = \{(p_1, q_1), \ldots, (p_m, q_m)\} \tag{3.16}$$

$$\cong \{[a_{p_1}, a_{q_1+1}], \ldots, [a_{p_m}, a_{q_m+1}]\}. \tag{3.17}$$

In practical applications, intervals where $a_p = a_{q+1}$ are usually omitted, as they represent ephemeral topological features.

**Example 3.2.10.** *[9, §2.3, Example] We consider the topological subspaces $S^0$, $S^2$, $S^4$, all of which are contractible. Meanwhile, $S^1$, $S^3$, $S^5$ are homeomorphic to the $0$-sphere, $1$-sphere, and $2$-sphere, respectively. This structural distinction leads to four distinct intervals in the persistence diagram of the absolute homology of a sphere, specifically $\mathcal{S}^2$:*

$$\mathrm{Pers}(H_\bullet(\mathcal{S}^2)) = \{(0, 5)_0, (1, 1)_0, (3, 3)_1, (5, 5)_2\}$$
$$\cong \{[0, \infty)_0, [1, 2)_0, [3, 4)_1, [5, \infty)_2\}. \tag{3.18}$$

*Here, the subscript $k$ in $(p, q)_k$ or $[a_p, a_{q+1}]_k$ denotes a topological feature in the $k$-dimensional homology.*

### 3.2.3 The Four Standard Persistence Modules

The standard module of absolute persistent homology, $H_\bullet(\mathcal{X})$, illustrates how the absolute homology groups $H_\bullet(X^i)$ relate to each other as the index $i$ changes. Similar observations can be made by considering the absolute cohomology groups $H^\bullet(X^i)$, the relative homology groups $H_\bullet(X^d, X^i)$, and the relative cohomology groups $H^\bullet(X^d, X^i)$ [9, §2.4]. The persistence diagram for absolute cohomology is represented as a multiset of integer pairs $(p, q)$, where $0 \le p \le q \le d$ for a finite filtration. For relative homology and cohomology, the persistence diagram consists of a multiset of tuples $(p, q)$ where $0 \le p \le q \le d - 1$ for a finite filtration. In each case, we interpret $(p, q)$ as a half-open interval $[a_p, a_{q+1})$ with the convention that $a_0 = -\infty$ and $a_{d+1} = \infty$.

**Example 3.2.11.** *[9, §2.4, Example] For $\mathcal{S}^2$ we yield*

$$\mathrm{Pers}(H_\bullet(S^5, \mathcal{S}^2)) = \{(0,0)_0, (1,1)_1, (3,3)_2, (0,4)_2\} \tag{3.19}$$

$$\cong \{[-\infty, 0)_0, [1, 2)_1, [3, 4)_2, [-\infty, 5)_2\}. \tag{3.20}$$

*At index 1 there is a nontrivial element of $H_1(S^5, S^1)$ represented by any arc connecting the two points of $S^1$ – the homology class is $[\sigma_2^{(1)}] = [\sigma_3^{(1)}]$. This class vanishes in $H_1(S^5, \mathcal{S}^2)$, thus we yield the interval $[1, 2)$.*

## 3.3 Persistent (Co)homology

Inverse problems primarily involve inferring geometric shapes from measurements like path integrals. Classical methods such as Fourier transforms provide extensive information but struggle with nonlinearity and ill-posed conditions, requiring substantial regularisation. Topology, particularly through persistent homology, offers alternative methods for deducing topological rather than geometric information. This approach is especially useful in high-dimensional, discrete sets of points, exemplified in the finite case by geological sonar to detect subterranean features based on density variations [9, §1]. Persistent homology identifies topological features represented as intervals in a barcode or persistence diagram, crucial for understanding the presence and persistence of features such as holes or voids in topological spaces. This method is statistically robust and can provide both qualitative and quantitative insights into point sets, which we suspect to lie on some compact topological object [5, 7]. We follow the results of de Silva, Morozov, and Vejdemo-Johansson for our explanations [9, §1]. In particular, there are at least four naturally arising persistent objects that can be extracted from a filtration of any topological space. They are

$$\text{persistent} \left\{ \begin{array}{c} \text{absolute} \\ \text{relative} \end{array} \right\} \left\{ \begin{array}{c} \text{homology} \\ \text{cohomology} \end{array} \right\}.$$

We address the computation of barcodes for all four types of persistent objects. We demonstrate that both absolute and relative (co)homologies yield identical barcodes and that transitions between these states are facilitated by established duality principles. The duality between homology and cohomology is akin to the duality in vector spaces, whereas a global duality specific to persistent topology allows for a unique interchange:

$$\text{Absolute homology} \rightleftarrows \text{relative cohomology.}$$
$$\text{Absolute cohomology} \rightleftarrows \text{relative homology.}$$

The main results suggest that a single calculation is sufficient to compute all four persistent objects due to the commutative nature of the global duality.

### 3.3.1 Barcode Isomorphisms

We characterise the multisets for persistence modules that are decomposable into interval modules. The persistence diagram partitions into $\text{Pers}_0$, comprising finite intervals $[a, b)$ as per [9, §2.3], and $\text{Pers}_\infty$, consisting of intervals $[a, \infty)$. This leads to the decomposition $\text{Pers} = \text{Pers}_0 \cup \text{Pers}_\infty$.

In this chapter, we establish that persistent homology and cohomology yield the same intervals, or barcodes, for both absolute and relative (co-)homology frameworks. This equivalence necessitates invoking the Universal Coefficients Theorem from algebraic topology, which we will prove beforehand. The Universal Coefficient Theorem for Cohomology elegantly ties together the cohomology of a space with coefficients in any abelian group $G$ to the homology of the space with integer coefficients. Specifically, as notation for this proof, $H_d(X; \mathbb{Z})$ and $H^d(X; \mathbb{Z})$ denote the $d$-th singular homology and cohomology groups with coefficients in $\mathbb{Z}$, respectively. We further involve $\text{Hom}(A, G)$, denoting group homomorphisms from an abelian group $A$ to another abelian group $G$, and $\text{Ext}^1(A, G)$ [18, §3.1, p.195], which measures obstructions in the splitting of a short exact sequence of abelian groups. Further, we use the properties of derived functors [27, §2.7].

**Theorem 3.3.1** (Universal coefficients for cohomology). *[18, Theorem 3.2] Let X be a topological space and G an abelian group. For any integer $d \geq 0$, there is a short exact sequence:*

$$0 \to \text{Ext}^1(H_{d-1}(X; \mathbb{Z}), G) \to H^d(X; G) \to \cdots$$
$$\cdots \to \text{Hom}(H_d(X; \mathbb{Z}), G) \to 0, \tag{3.21}$$

*which splits, though not canonically.*

*Proof.* Let $(\mathcal{C}(X), \partial)$ be the singular chain complex of a topological space $X$ with integer coefficients. The homology groups $H_d(X; \mathbb{Z})$ are defined as:

$$H_d(X; \mathbb{Z}) := \ker(\partial_d)/\text{im}(\partial_{d+1}), \tag{3.22}$$

where $\partial_d$ are the boundary maps in $\mathcal{C}(X)$. The chain group $C_d$ consists of formal sums of singular $d$-simplices in $X$ with integer coefficients. The boundary maps $\partial_d : C_d \to C_{d-1}$ are defined by

$$\partial_d(\tilde{\sigma}^{(d)}) = \sum_{i=0}^{d}(-1)^i \tilde{\sigma}^{(d)}|_{[v_0,\ldots,\hat{v}_i,\ldots,v_d]}, \tag{3.23}$$

where $\tilde{\sigma}^{(d)} : \tilde{\sigma}^d \to X$ is a singular simplex, and $\tilde{\sigma}|_{[v_0,\ldots,\hat{v}_i,\ldots,v_d]}$ denotes the restriction of $\tilde{\sigma}$ to the $i$-th face of the simplex, omitting the $i$-th vertex. The functor $\text{Hom}(-, G)$ for an abelian group $G$ applied to $C_d$ yields a group $\text{Hom}(C_d, G)$, and the coboundary $\delta^d$ for the cochain complex $\text{Hom}(\mathcal{C}^\flat(X), G)$ — see Eq. 3.3.2 — is $\delta^d(f) = f \circ \partial_{d+1}$ for $f \in \text{Hom}(C_{d+1}(X), G)$. This leads to the cohomology groups $H^d(X; G) := \ker(\delta^d)/\text{im}(\delta^{d-1})$. We consider the projective resolution of $\mathcal{C}^\flat(X)$ to effectively apply the $\text{Ext}^1$–functor. For an abelian group $A$, $\text{Ext}^1(A, G) = R^1 \text{Hom}(A, G)$,

where $R^1$ denotes the first right derived functor of Hom [27, Vista 3.4.6, §3.5]. To show this equality, we begin by taking a projective resolution of $A$ [27, Definition 2.2.4]: $\cdots \to P_2 \to P_1 \to P_0 \to A \to 0$, where each $P_i$ is a projective abelian group. Then we can apply the functor $\mathrm{Hom}(-, G)$ to the projective resolution:

$$0 \to \mathrm{Hom}(P_0, G) \to \mathrm{Hom}(P_1, G) \to \mathrm{Hom}(P_2, G) \to \cdots. \qquad (3.24)$$

This sequence is exact on the left because $\mathrm{Hom}(-, G)$ is left exact and each $P_i$ is projective. The first right derived functor $R^1 \mathrm{Hom}(A, G)$ is defined as the cohomology of this sequence at $P_1$:

$$R^1 \mathrm{Hom}(A, G) = \frac{\ker(\mathrm{Hom}(P_1, G) \to \mathrm{Hom}(P_2, G))}{\mathrm{im}(\mathrm{Hom}(P_0, G) \to \mathrm{Hom}(P_1, G))}. \qquad (3.25)$$

The group $\mathrm{Ext}^1(A, G)$ classifies extensions of $A$ by $G$, equivalent to the kernel / image calculation in the cohomology of the Hom–sequence $R^1 \mathrm{Hom}(A, G)$. Given $H_d(X; \mathbb{Z})$, consider the short exact sequence obtained from the projective resolution of $\mathbb{Z}$: $0 \to \mathbb{Z} \to F \to H_d(X; \mathbb{Z}) \to 0$, where $F$ is free. Applying $\mathrm{Hom}(-, G)$ gives

$$0 \to \mathrm{Hom}(H_d(X; \mathbb{Z}), G) \to \mathrm{Hom}(F, G) \to \cdots$$
$$\cdots \to \mathrm{Hom}(\mathbb{Z}, G) \to \mathrm{Ext}^1(H_d(X; \mathbb{Z}), G) \to 0. \qquad (3.26)$$

We apply $\mathrm{Ext}^1$, which gives rise to the long exact sequence of $\mathrm{Ext}^1$-groups:

$$0 \to \mathrm{Hom}(H_d(X; \mathbb{Z}), G) \to H^d(X; G) \to \cdots$$
$$\cdots \to \mathrm{Ext}^1(H_{d-1}(X; \mathbb{Z}), G) \to 0. \qquad (3.27)$$

Finally, we verify exactness and splitting. The term $\mathrm{Ext}^1(H_{d-1}(X; \mathbb{Z}), G)$ measures the non-trivial extensions of $G$ by $H_{d-1}(X; \mathbb{Z})$, which corresponds to the obstructions to lifting $H_{d-1}(X; \mathbb{Z})$ linearly over $G$. The term $\mathrm{Hom}(H_d(X; \mathbb{Z}), G)$ represents the group homomorphisms from $H_d(X; \mathbb{Z})$ to $G$, which naturally includes in $H^d(X; G)$. The sequence is exact at each stage by the properties of derived functors and their application to the singular chain complex [27, Horseshoe Lemma 2.2.8]. The sequence ends with $0$ because $\mathrm{Ext}^1$ of a projective (or free) module vanishes, and $\mathbb{Z}$ is free. The sequence splits because the functor $\mathrm{Hom}(-, G)$ preserves products and coproducts. However, the way it splits is not canonical and depends on the choice of a splitting homomorphism, which is not unique. For cohomology the same identity holds with reversed arrow, according to isomorphisms between homology and cohomology, which gives the desired result.

In the generalization of the Universal Coefficient Theorem 3.3.1 to the case of modules over a PID, the $\mathrm{Ext}^1$ terms vanish since $\mathbb{F}$ is a field, thus we obtain $H^d(X; \mathbb{F}) \cong \mathrm{Hom}(H_d(X; \mathbb{F}), \mathbb{F})$ [18, §3.3.1, p.198]. $\qquad \square$

**Theorem 3.3.2.** *[9, Proposition 2.3] For all integers $d \geq 0$, it holds that:*

$$\mathrm{Pers}(H_d(\mathcal{X}; \mathbb{F})) = \mathrm{Pers}(H^d(\mathcal{X}; \mathbb{F})), \qquad (3.28)$$

$$\mathrm{Pers}(H_d((X^\infty, \mathcal{X}); \mathbb{F})) = \mathrm{Pers}(H^d((X^\infty, \mathcal{X}); \mathbb{F})). \qquad (3.29)$$

*Proof.* When considering coefficients in a field $\mathbb{F}$ rather than in a ring, the Universal Coefficient Theorem 3.3.1 assures us of a natural isomorphism between the $d$-th cohomology group and the homomorphisms from the $d$-th homology group to the base field: $H^d(X; \mathbb{F}) \cong \operatorname{Hom}(H_d(X; \mathbb{F}), \mathbb{F})$. Therefore, the associated maps

$$H_d(X^i; \mathbb{F}) \to H_d(X^j; \mathbb{F}) \quad \text{and} \quad H^d(X^i; \mathbb{F}) \leftarrow H^d(X^j; \mathbb{F}) \tag{3.30}$$

are adjoint and hence possess the same rank. Since the persistence intervals over a field are uniquely determined by the dimension and rank of the homology vector spaces, it follows that this holds for both homology and cohomology. Consequently, they share the same barcode. $\qquad\square$

Consider a finite filtration $X^0 \subset X^1 \subset X^2 \subset \cdots \subset X^d \subseteq X^\infty$ of a topological space $X$. The homology groups $H_k(X^d)$ for some fixed dimension $k \leq d$ serve as the initial terms for the relative homology groups $H_k(X^\infty, X)$. Since $H_k(X)$ is consistent with $H_k(X^d)$ as $d$ approaches infinity, these sequences can be unified into a single sequence: $H_k(X) \to H_k(X^\infty, X)$. For this concatenated sequence, the indices are

$$\{0, 1, 2, \dots, d = \bar{0}, \bar{1}, \bar{2}, \dots, \overline{(d-1)}\}, \tag{3.31}$$

using the barred numbers to indicate relative homology part of the sequence. This structure allows us to discuss the persistence diagram associated with this homological configuration. The persistence intervals in this diagram can generally be categorized into three types:

1. $(p, q)$ where $0 \leq p \leq q < d$, denoted as $[p, q+1)$ or $[a_p, a_{q+1})$.
2. $(\bar{p}, \bar{q})$ where $0 < p \leq q \leq d-1$, denoted as $[\bar{p}, \bar{q}+1)$ or $[a_{\bar{p}}, a_{\bar{q}+1})$.
3. $(p, \bar{q})$ where $0 \leq p \leq d, 0 \leq q \leq d-1$, denoted as $[p, \bar{q}+1)$ or $[a_p, a_{\bar{q}+1})$.

**Corollary 3.3.3.** *[9, Proposition 2.5] The barcode* $\operatorname{Pers}(H_k(\mathcal{X}) \to H_k(X^\infty, \mathcal{X}))$ *comprises the following collections of intervals, which match the previously discussed types of intervals bijectively:*

1. *An interval $[a, b)$ for every interval $[a, b)$ in $\operatorname{Pers}_0(H_k(\mathcal{X}))$.*
2. *An interval $[\bar{a}, \bar{b})$ for every interval $[a, b)$ in $\operatorname{Pers}_0(H_{k-1}(\mathcal{X}))$.*
3. *An interval $[a, \bar{a})$ for every interval $[a, \infty)$ in $\operatorname{Pers}_\infty(H_k(\mathcal{X}))$.*

*Proof.* We begin by analyzing the first two types of intervals in the persistence diagram $\operatorname{Pers}(H_k(X) \to H_k(X^\infty, X))$. These intervals either do not intersect the intermediate term $H_k(X^d)$ or terminate before it. Consequently, they correspond precisely to the finite intervals in $\operatorname{Pers}(H_k(X))$ and $\operatorname{Pers}(H_k(X^\infty, X))$, clarifying the first two cases. The correspondence $\operatorname{Pers}_0(H_k(X^\infty, X)) = \operatorname{Pers}_0(H_{k-1}(X))$ helps in mapping these relationships.

The third case requires examining intervals of the form $[a, \bar{b})$ and proving that they are invariably of the form $[a, \bar{a})$, which means that the paired intervals $[a, \infty)$ and $[-\infty, a)$ in $\operatorname{Pers}_\infty(H_k(\mathcal{X}))$ and $\operatorname{Pers}_\infty(H_k(X^\infty, \mathcal{X}))$ are restrictions of a single interval $[a, \bar{a})$ in the concatenated sequence [9, Proof of Proposition 2.5]. To establish this, we compare the ascending filtration defined by the images of $H_k(X^i)$ in $H_k(X^d)$ for $i = 0, 1, 2, \dots, d-1$ (denoted as $\operatorname{im}(H_k(X^i) \to H_k(X^d))$) with the descending filtration defined by the kernels of $H_k(X^d)$ in $H_k(X^\infty, X)$ for $i = 0, 1, 2, \ldots, d-1$ (denoted as $\ker(H_k(X^d) \to H_k(X^d, X^i))$). This examination revolves around the fundamental properties of the homology groups in a filtration setting. For each index $i$, the image and kernel correspond to the same subspace of $H_k(X^d)$. This equivalence is guaranteed by the homology long exact sequence associated with the pair $(X^d, X^i)$, which links the relative and absolute homology groups. Specifically, the exact sequence implies that any cycle in $\mathrm{im}(H_k(X^i) \to H_k(X^d))$ that becomes a boundary in $H_k(X^d, X^i)$ must vanish, thus equating the image and kernel. As a result, both filtrations align perfectly, establishing that the third type of interval indeed maps to self-closing intervals of the form $[a, \bar{a}]$. $\qquad \square$

### 3.3.2 Persistent Chain Complexes

To enhance readability, and because we mostly move in the horizontal direction of a persistent chain complex, see Fig. 3.3, we omit in this chapter the dimension $\sigma_i^{(k)} \in C_k^i$ and write $\sigma_i \in C^i$ instead as well as for the boundary operator, as the results hold for all dimensions. Further, we apply a slight change in notation. The filtration runs not anymore from $0, 1, \ldots, d$ but from $0, 1, \ldots, n$. We will need the variable $d$ to indicate the endpoint of a persistence interval or death of a homological feature. Following [9, §2.6], the standard persistence module of a topological space $X$ can be described via a filtered CW-complex $\mathcal{X} := \sigma_0 \cup \sigma_1 \cup \cdots \cup \sigma_n$, adding a single cell for each filtration step, so that the persistence module is the sequence

$$\mathcal{C}(X): \quad C^0 \to C^1 \to C^2 \to \cdots \to C^n,$$

for a finite filtration, where each $C^i := \langle \sigma_0, \sigma_1, \ldots, \sigma_i \rangle$ is a vector space (or an abelian group) over a field $\mathbb{F}$, generated by the simplices $\sigma_j$. The boundary operator is defined by $\partial \sigma_j = \sum_{i<j} \lambda_{ij} \sigma_i$ for $\lambda_{ij} \in \mathbb{F}$ for all $i, j \in \{0, 1, \ldots, n\}$. Then $\mathcal{C}(X) := \{C^i, \partial\}_{i=0}^n$ is the chain complex of $X$, and $\mathcal{C}(\mathcal{X})$ represents the persistent version for $\mathcal{X}$, see Fig. 3.3. We define persistent absolute homology of $\mathcal{X}$ as

$$H_\bullet(\mathcal{X}) = H(\mathcal{C}(\mathcal{X}), \partial): \quad \frac{\ker(\partial_0)}{\mathrm{im}(\partial_1)} \to \frac{\ker(\partial_1)}{\mathrm{im}(\partial_2)} \to \frac{\ker(\partial_2)}{\mathrm{im}(\partial_3)} \to \cdots \to \frac{\ker(\partial_d)}{\mathrm{im}(\partial_{d+1})}.$$

For the persistent absolute cohomology $H^\bullet(\mathcal{X})$, we define

$$\mathcal{C}^\flat(X): \quad C^{\flat 0} \leftarrow C^{\flat 1} \leftarrow C^{\flat 2} \leftarrow \cdots \leftarrow C^{\flat n},$$

where $C^{\flat i} = \mathrm{Hom}(C^i, \mathbb{F}) = \langle \sigma_0^\flat, \sigma_1^\flat, \ldots, \sigma_i^\flat \rangle$, with the elements $\sigma_i^\flat$ as the dual basis of $\sigma_i$. The coboundary $\delta = \partial^\flat$ is defined as the adjoint of $\partial$. Considering $(\mathcal{C}^\flat(X), \delta)$, its homology forms a module structure $H^\bullet(\mathcal{X}) = H(\mathcal{C}^\flat(\mathcal{X}), \delta) = \ker(\delta)/\mathrm{im}(\delta)$. This is a persistence module, with morphisms in the reverse direction.

**Example 3.3.4.** *[9, §2.6, Example] For $\mathcal{S}^2$, the boundary operator is given by*

$$\partial_0 \sigma_0^{(0)} = \partial_0 \sigma_1^{(0)} = 0,$$
$$\partial_1 \sigma_2^{(1)} = \partial_1 \sigma_3^{(1)} = \sigma_0^{(0)} - \sigma_1^{(0)},$$

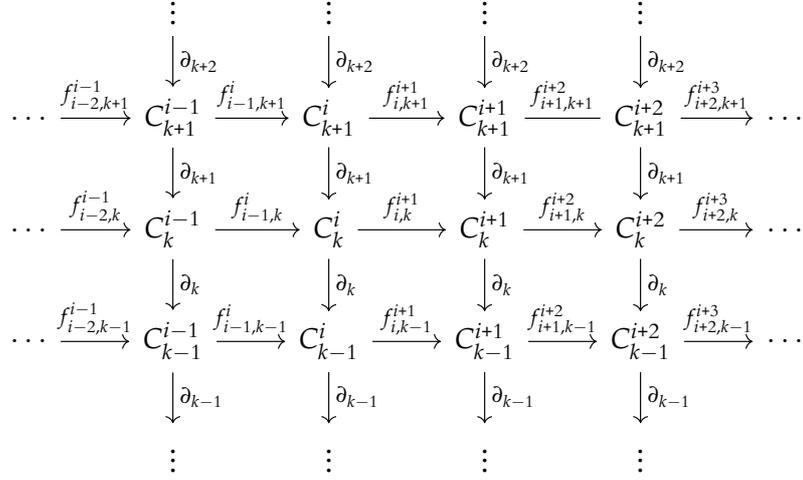

Figure 3.3: A persistent chain complex, where moving to the right increases the filtration index, while moving downwards decreases the dimension. The maps $f_{i,k}^{i+1}$ correspond to the structure map from the $i$-th to the $(i+1)$-st chain group over $\mathbb{Z}$ – or vector space over $\mathbb{F}$ – in dimension $k$. The boundary operator is defined for the chain complex, restricting it for each element to the respective chain group.

$$\partial_2 \sigma_4^{(2)} = \partial_2 \sigma_5^{(2)} = \sigma_2^{(1)} - \sigma_3^{(1)}. \qquad (3.32)$$

*Similarly, the coboundary operator is given by*

$$\begin{aligned}
\delta^0 \sigma_0^{(0)\flat} &= -\delta^0 \sigma_1^{(0)\flat} = \sigma_2^{(1)\flat} + \sigma_3^{(1)\flat}, \\
\delta^1 \sigma_2^{(1)\flat} &= -\delta^1 \sigma_3^{(1)\flat} = \sigma_4^{(2)\flat} + \sigma_5^{(2)\flat}, \\
\delta^2 \sigma_4^{(2)\flat} &= -\delta^2 \sigma_5^{(2)\flat} = 0.
\end{aligned} \qquad (3.33)$$

The data in Eq. 3.32 and in Eq. 3.33 can be summarized in a matrix and its transposed, respectively. The relative homology and cohomology persistence modules are then defined as the homology of the persistence modules:

$$(C^n/\mathcal{C}(X)): \quad C^n \to (C^n/C^0) \to (C^n/C^1) \to (C^n/C^2) \to \cdots \to (C^n/C^{n-1}),$$
$$(C^n/\mathcal{C}(X))^\flat: \quad C^{\flat n} \leftarrow (C^n/C^0)^\flat \leftarrow (C^n/C^1)^\flat \leftarrow (C^n/C^2)^\flat \leftarrow \cdots \leftarrow (C^n/C^{n-1})^\flat.$$

We note that the mappings $\to$ from $\mathcal{C}(X)$ and the mappings $\leftarrow$ from $(C^n/\mathcal{C}(X))^\flat$ are injective, while the mappings $\leftarrow$ from $\mathcal{C}^\flat(X)$ and the mappings $\to$ from $(C^n/\mathcal{C}(X))$ are surjective. Thus, absolute and relative cohomology are structurally similar but qualitatively different from absolute homology and relative homology.

**Theorem 3.3.5** (Partition of persistence modules). *[9, Theorem 2.6] Let $\mathcal{C}(X)$ be given and $\partial$ as above, then there exists a partition $\{0, 1, 2, \ldots, n\} = B^\ltimes \sqcup B \sqcup D$ with a bijective pairing $B \leftrightarrow D$, such that $b \in B$ is paired with $d \in D$ if and only if $[b, d) \in \mathrm{Pairs}(\mathcal{C}(X), \partial)$. Furthermore, there is a basis $\hat{\sigma}_0, \hat{\sigma}_1, \hat{\sigma}_2, \ldots, \hat{\sigma}_n$ of $C^n$ such that*

- $C^i = \langle \hat{\sigma}_0, \hat{\sigma}_1, \ldots, \hat{\sigma}_i \rangle$ for each $i$.
- $\partial \hat{\sigma}_b = 0$ for all $b \in B^\times$.
- $\partial \hat{\sigma}_d = \hat{\sigma}_b$, and thus $\partial \hat{\sigma}_b = 0$, for all $[b, d] \in$ Pairs.

It follows that the persistence diagram $\mathrm{Pers}(H(\mathcal{C}(X), \partial))$ consists of $[a_b, \infty)$ for $b \in B^\times$ together with intervals $[a_b, a_d)$ for $[b, d] \in$ Pairs.

*Proof.* Let $\mathcal{C}(X)$ be a finite filtered chain complex associated with a topological space, structured over a field $\mathbb{F}$. The filtration indices $\{0, 1, 2, \ldots, n\}$ classify the stages at which elements are born or die in homology. The partition of $\{0, 1, 2, \ldots, n\}$ into $B^\times \sqcup B \sqcup D$ is then constructed such that $B^\times$ contains indices corresponding to elements that persist indefinitely, i.e., those for which $\partial \sigma_i = 0$ and $\sigma_i$ does not become a boundary at any higher index. $B$ and $D$ are paired bijectively, where each $b \in B$ (births) corresponds uniquely to a $d \in D$ (deaths), signifying the termination of the homological feature introduced at filtration step $d$. This pairing is determined through the matrix reduction of $\partial$, ensuring that every cycle born at stage $b$ becomes a boundary at stage $d$. We then construct a basis $\hat{\sigma}_i$ such that for each $i$, $\hat{\sigma}_i$ is a generator of $C^i$. Specifically, if $i \in B^\times$, then $\partial \hat{\sigma}_i = 0$, indicating these elements are cycles that persist indefinitely. For pairs $[b, d] \in \mathrm{Pairs}(\mathcal{C}(X), \partial)$ with $b \in B$ and $d \in D$, set $\partial \hat{\sigma}_d = \hat{\sigma}_b$ and $\partial \hat{\sigma}_b = 0$, reflecting the fact that the cycle $\hat{\sigma}_b$ born at $b$ becomes a boundary at $d$ and thus ceases to contribute to homology past $d$. Hence, the persistence diagram $\mathrm{Pers}(H(\mathcal{C}(X), \partial))$ contains intervals $[a_b, \infty)$ for each $b \in B^\times$, indicating persistent homological features and intervals $[a_b, a_d)$ for each $[b, d] \in$ Pairs, describing features with finite lifespans, starting as a cycle at $b$ and terminating as a boundary at $d$. $\qquad \square$

**Remark 3.3.6.** *[12, §2, p.516] In Theorem 3.3.5 the index sets can be reexpressed in the original language of positive and negative simplices proposed by Edelsbrunner et al. [9, p.8], where a simplex is called positive if it belongs to a cycle and negative if it doesn't, respectively:*
- *$B$ identifies the positive simplices that remain unpaired,*
- *$B^\times$ identifies the positive simplices that do become paired,*
- *$D$ identifies the negative simplices.*

*The simplices $\hat{\sigma}_b$ and $\hat{\sigma}_d$ are cycles characterized by their leading terms $\sigma_b$ and $\sigma_d$ respectively, while the simplex $\hat{\sigma}_d$ is a chain with leading term $\sigma_d$. This chain 'kills' the homology class of its paired cycle $\hat{\sigma}_b$ through the boundary relation $\partial \hat{\sigma}_d = \hat{\sigma}_b$.*

**Theorem 3.3.7.** *[9, Proposition 2.4] For all integers $d \geq 0$, it holds that:*

$$\mathrm{Pers}_0(H_d(\mathcal{X})) \cong \mathrm{Pers}_0(H_{d+1}(X^\infty, \mathcal{X})),$$
$$\mathrm{Pers}_\infty(H_d(X^\infty, \mathcal{X})) \cong \mathrm{Pers}_\infty(H_d(X^\infty, \mathcal{X})).$$

*Proof. [9, Proof of Proposition 2.4]* The decomposition $\{0, 1, 2, \ldots, n\} = B^\times \sqcup B \sqcup D$ and the new basis $\hat{\sigma}_0, \hat{\sigma}_1, \ldots, \hat{\sigma}_n$ allow us to write the persistent chain complex $\mathcal{C}(\mathcal{X})$ as a direct sum of persistent chain complexes, such that

$$\mathcal{C}(X) = \bigoplus_{b^\times \in B^\times} \mathcal{C}(X)^{b^\times} \oplus \bigoplus_{[b, d] \in \mathrm{Pairs}} \mathcal{C}(X)^{bd}$$

$$= \bigoplus_{b^{\times} \in B^{\times}} \langle \hat{\sigma}_{b^{\times}} \rangle \oplus \bigoplus_{[b,d) \in \text{Pairs}} \langle \hat{\sigma}_b, \hat{\sigma}_d \rangle. \tag{3.34}$$

The boundary operator is compatible with this decomposition, where $\partial$ is defined for each summand as the corresponding boundary operator of the chain complex. Therefore, we can compute the homology $H(C^i/\mathcal{C}(X), \partial)$ separately for each summand. For summands of the form $\mathcal{C}(X)^{b^{\times}}$, the persistence modules are constant over two phases, with indices in $\{0, \dots, b^{\times} - 1\}$ and $\{b^{\times}, \dots, n-1\}$. Thus,

$$
\begin{aligned}
(C^{b^{\times}}/\mathcal{C}(X)^{b^{\times}}): & \quad \langle \hat{\sigma}_{b^{\times}} \rangle \to 0, \\
\ker(\partial): & \quad \langle \hat{\sigma}_{b^{\times}} \rangle \to 0, \\
\text{im}(\partial): & \quad 0 \to 0, \\
H := \ker(\partial)/\text{im}(\partial): & \quad \langle \hat{\sigma}_{b^{\times}} \rangle \to 0.
\end{aligned} \tag{3.35}
$$

It follows that $H(C^{b^{\times}}/\mathcal{C}(X)^{b^{\times}})$ contributes an interval of the form $[-\infty, a_{b^{\times}})$, generated by the equivalence class of $[\hat{\sigma}_{b^{\times}}] \in \text{Pers}(H(\mathcal{C}(X)))$, and thus has the same homological dimension as $[a_{b^{\times}}, \infty) \in \text{Pers}(H(\mathcal{C}(X)))$. For summands of the type $\mathcal{C}(X)^{bd}$, the persistence modules can be divided into constants over three phases, with indices in the ranges $\{0, \dots, b-1\}$, $\{b, \dots, b^{\times} - 1\}$, and $\{b^{\times}, \dots, n-1\}$, or explicitly:

$$
\begin{aligned}
(C^{bd}/\mathcal{C}(X)^{bd}): & \quad \langle \hat{\sigma}_b, \hat{\sigma}_d \rangle \to \langle \hat{\sigma}_d \rangle \to 0, \\
\ker(\partial): & \quad \langle \hat{\sigma}_b \rangle \to \langle \hat{\sigma}_d \rangle \to 0, \\
\text{im}(\partial): & \quad \langle \hat{\sigma}_b \rangle \to 0 \to 0, \\
H := \ker(\partial)/\text{im}(\partial): & \quad 0 \to \langle \hat{\sigma}_b \rangle \to 0.
\end{aligned} \tag{3.36}
$$

It follows that $H(C^{bd}/\mathcal{C}(X)^{bd})$ contributes a single interval of the form $[a_b, a_d)$, generated by the equivalence class $[\hat{\sigma}_b]$, and thus has one dimension more than $[a_b, a_d) \in \text{Pers}(H(\mathcal{C}(X)))$, which is generated by $[\hat{\sigma}_d]$. $\qquad \square$

**Remark 3.3.8.** *[9, Proposition 2.4] In this case, we get an isomorphism of multisets. This is due to the identification of the intervals $[a, \infty) \leftrightarrow [-\infty, a)$ for $\text{Pers}_{\infty}$. Thus, persistent homology and relative homology barcodes carry the same information, with a dimension shift for the finite intervals.*

In this closing argument, we will prove the four statements in Morozov's paper that were left open [9, §2.8].

Consider a finite filtration of a topological space $X^0 \subset X^1 \subset \cdots \subset X^n := \sigma_0 \cup \sigma_1 \cup \cdots \cup \sigma_n = X$, with chain complexes $\mathcal{C}(X) := (C^i, \partial)_{i \in \{0, \dots, n\}}$ over $\mathbb{F}$. We regard $\mathcal{C}(X)$ as a free module over $\mathbb{F}[i]$ with $(n+1)$ generators $\sigma_0, \sigma_1, \dots, \sigma_n$, where $\sigma_i$ has degree $i$. The associated chain complexes $\mathcal{C}(X)$ over $\mathbb{F}$ have inclusions $X^i \subseteq X^{i+1}$ that induce structure maps $f_i^{i+1} : C^i \to C^{i+1}$. The boundary operator $\partial : \mathcal{C}(X) \to \mathcal{C}(X)$ is a homomorphism of graded modules. We define the pointwise dual as the dual to the ground field $\mathbb{F}$, such that

$$\mathcal{C}^{\dagger}(X) = \text{Hom}_{\mathbb{F}}(\mathcal{C}(X), \mathbb{F}), \tag{3.37}$$

$$C^{\dagger n} : C^{-n} \to \mathbb{F} \text{ is a linear map.} \tag{3.38}$$

Analogously, we define the global dual as dual to the polynomial ring $\mathbb{F}[i]$:

$$\mathcal{C}^{\circ}(X) = \mathrm{Hom}_{\mathbb{F}[i]}(\mathcal{C}(X), \mathbb{F}[i]), \tag{3.39}$$

$$C^{\circ n} : \mathcal{C}(X) \to \mathbb{F}[i] \text{ graded module homomorphism of degree } n. \tag{3.40}$$

These are graded $\mathbb{F}[i]$-modules themselves. The operations $-^{\circ}$ and $-^{\dagger}$ are contravariant functors, thus the boundary map on $\mathcal{C}(X)$ induces boundary maps on the two new modules [9, §2.8].

**Theorem 3.3.9.**
  1. $H_{\bullet}(\mathcal{C}(X), \partial) \cong$ *persistent absolute homology of* $\mathcal{X}$.
  2. $H^{\bullet}(\mathcal{C}^{\dagger}(X), \partial^{\dagger}) \cong$ *persistent absolute cohomology of* $\mathcal{X}$.
  3. $H^{\bullet}(\mathcal{C}^{\circ}(X), \partial^{\circ}) \cong$ *persistent relative cohomology of* $\mathcal{X}$.
  4. $H_{\bullet}(\mathcal{C}^{\circ\dagger}(X), \partial^{\circ\dagger}) \cong$ *persistent relative homology of* $\mathcal{X}$.

*Proof.*

  1. The persistent homology groups associated with a filtration $(X^i)_{0 \leq i \leq d}$ are defined as

$$H_k(X^i \subseteq X^j) = \mathrm{im}\left(H_k(X^i) \xrightarrow{f_i^j} H_k(X^j)\right),$$

where $H_k(X^i)$ denotes the $k$-th homology group of $X^i$, and $f_i^j$ is the map induced by the inclusion $i : X^i \hookrightarrow X^j$. To formalize the algebraic structure of this filtration, we consider the free graded module $\mathcal{C}(X)$ over the polynomial ring $\mathbb{F}[i]$, generated by $\{\sigma_0, \sigma_1, \ldots, \sigma_n\}$, which correspond to the basis elements of the chain complexes $C^i$. The action of $i$ on $\mathcal{C}(X)$ encodes the inclusions $f_i^{i+1}$ via index shifting; specifically, multiplication by $i$ transitions from $C^i$ to $C^{i+1}$. The boundary operator $\partial : \mathcal{C}(X) \to \mathcal{C}(X)$ is defined by $\partial\sigma_i \in C^{i-1}$ for each generator $\sigma_i$. This guarantees that $\partial^2 = 0$, making $\partial$ an operator that respects the grading. Additionally, $\partial$ commutes with the action of $i$. The homology of the complex is given by $H_k(\mathcal{C}(X), \partial) = \ker(\partial)/\mathrm{im}(\partial)$, with $\partial : C^i \to C^{i+1}$. Thus, $H_k(\mathcal{C}(X), \partial)$ covers the persistent homology groups $H_k(X^i \subseteq X^j)$ for all $i \leq j$. We conclude with the isomorphism:

$$H(\mathcal{C}(X), \partial) \cong \bigoplus_{0 \leq i \leq j \leq n} H_k(X^i \subseteq X^j).$$

  2. The persistent cohomology groups associated with a filtration $(X^i)_{0 \leq i \leq n}$ are defined as

$$H^k(X^i \subseteq X^j) = \mathrm{im}\left(H^k(X^j) \xrightarrow{f_i^j} H^k(X^i)\right),$$

where $H^k(X^i)$ denotes the $k$-th cohomology group of $X^i$, and the map $f_i^j$ is induced by the inclusion $i : X^i \hookrightarrow X^j$. To define the dual chain complex, we consider the dual module $\mathcal{C}^{\dagger}(X) = \mathrm{Hom}_{\mathbb{F}}(\mathcal{C}(X), \mathbb{F})$. For the free module $\mathcal{C}(X)$ with generators $\{\sigma_0, \sigma_1, \ldots, \sigma_n\}$ corresponding to the chain complex $C^i$, the dual basis elements $\sigma_i^{\dagger} \in \mathcal{C}^{\dagger}(X)$ are defined by $\sigma_i^{\dagger}(\sigma_j) = \delta_{ij}$, where $\delta_{ij}$ is the Kronecker delta. The dual boundary operator $\partial^{\dagger} : \mathcal{C}^{\dagger}(X) \to \mathcal{C}^{\dagger}(X)$ is defined by $\partial^\dagger(\phi) = \phi \circ \partial$, for any $\phi \in \mathcal{C}^\dagger(X)$. This operator is a graded map that increases the degree by one and satisfies $(\partial^\dagger)^2 = 0$. We compute the cohomology of the dual complex $H^k(\mathcal{C}^\dagger(X), \partial^\dagger) = \ker(\partial^\dagger)/\mathrm{im}(\partial^\dagger)$, where $\partial^\dagger : C^i \to C^{i+1}$. To establish an isomorphism with persistent cohomology, we note that $\mathcal{C}^\dagger(X)$ cover the algebraic structure of cochains under the action of inclusions $f_i^{i+1}$ in the dual setting. The multiplication by $i$ in $\mathbb{F}[i]$ reflects these inclusions within the dual module, and $\partial^\dagger$ acts compatibly. Thus, $H^k(\mathcal{C}^\dagger(X), \partial^\dagger)$ encodes the persistent cohomology groups $H^k(X^i \subseteq X^j)$ for all $0 \le i \le j \le n$. Consequently, we conclude that

$$H(\mathcal{C}^\dagger(X), \partial^\dagger) \cong \bigoplus_{0 \le i \le j \le n} H^k(X^i \subseteq X^j).$$

3. The persistent relative cohomology groups associated with a filtration $(X^i)_{0 \le i \le n}$ are defined as

$$H^k(X^i, X^j) = \mathrm{coker}\left(H^k(X^j) \xrightarrow{f_i^j} H^k(X^i)\right),$$

where $H^k(X^i)$ denotes the $k$-th cohomology group of $X^i$, and the map $f_i^j$ is induced by the inclusion $i : X^i \hookrightarrow X^j$. To construct the dual boundary operator, we define $\partial^\circ : \mathcal{C}^\circ(X) \to \mathcal{C}^\circ(X)$ by $\partial^\circ(\psi) = \psi \circ \partial$, for any $\psi \in \mathcal{C}^\circ(X)$. This operator is a graded module homomorphism of degree 1, ensuring that $\partial^\circ$ satisfies $(\partial^\circ)^2 = 0$, due to the fact that $\partial^2 = 0$ in the original chain complex. We compute the cohomology of the dual complex $H^k(\mathcal{C}^\circ(X), \partial^\circ) = \ker(\partial^\circ)/\mathrm{im}(\partial^\circ)$, where $\partial^\circ : C^i \to C^{i+1}$. To establish an isomorphism with persistent relative cohomology, we observe that the module $\mathcal{C}^\circ(X)$ encodes the algebraic structure of cochains with coefficients in $\mathbb{F}[i]$, reflecting the inclusions of the filtration in a dual setting. The action of $i$ in $\mathbb{F}[i]$ represents the structure maps $f_i^{i+1}$ from $C^i$ to $C^{i+1}$, while $\partial^\circ$ acts compatibly as a graded map. Thus, $H^k(\mathcal{C}^\circ(X), \partial^\circ)$ captures the persistent relative cohomology groups $H^k(X^i, X^j)$ for all $0 \le i \le j \le d$. Therefore, we conclude that

$$H(\mathcal{C}^\circ(X), \partial^\circ) \cong \bigoplus_{0 \le i \le j \le n} H^k(X^i, X^j).$$

4. We define the double dual complex as

$$\mathcal{C}^{\circ\dagger}(X) = \mathrm{Hom}_\mathbb{F}(\mathcal{C}^\circ(X), \mathbb{F}),$$

where $\mathcal{C}^\circ(X) = \mathrm{Hom}_{\mathbb{F}[i]}(\mathcal{C}(X), \mathbb{F}[i])$ is the global dual complex. Thus, $\mathcal{C}^{\circ\dagger}(X)$ serves as the dual of $\mathcal{C}^\circ(X)$ over the field $\mathbb{F}$, capturing cochains with coefficients in $\mathbb{F}$. The dual boundary operator $\partial^{\circ\dagger} : \mathcal{C}^{\circ\dagger}(X) \to \mathcal{C}^{\circ\dagger}(X)$ is defined by $\partial^{\circ\dagger}(\phi) = \phi \circ \partial^\circ$, for any $\phi \in \mathcal{C}^{\circ\dagger}(X)$. This definition guarantees that $(\partial^{\circ\dagger})^2 = 0$, since $(\partial^\circ)^2 = 0$ holds in the dual complex $\mathcal{C}^\circ(X)$. We compute the cohomology of the double dual complex $H_k(\mathcal{C}^{\circ\dagger}(X), \partial^{\circ\dagger}) = \ker(\partial^{\circ\dagger})/\mathrm{im}(\partial^{\circ\dagger})$, where $\partial^{\circ\dagger} : C^i \to C^{i+1}$. To establish an isomorphism with persistent relative homology, we observe that the module $\mathcal{C}^{\circ\dagger}(X)$ reflects the algebraic structure of chains in the dual–dual setting, with coefficients in $\mathbb{F}$. The dual action captures the inverse of the structure maps $f_i^{i+1}$ from $C^i$ to $C^{i+1}$, while $\partial^{\circ\dagger}$ acts compatibly as a graded map. Thus, $H_k(\mathcal{C}^{\circ\dagger}(X), \partial^{\circ\dagger})$ encodes the persistent relative homology groups $H_k(X^i, X^j)$ for all $0 \leq i \leq j \leq n$. Finally, we conclude that

$$H(\mathcal{C}^{\circ\dagger}(X), \partial^{\circ\dagger}) \cong \bigoplus_{0 \leq i \leq j \leq n} H_k(X^i, X^j).$$

$\square$

### 3.3.3 Cohomology of Chain Complexes

Persistent relative cohomology is structurally similar to persistent absolute homology. We adopt Morozov's notation for the reversed sequence [9, §2.7]:

$$(C^n/\mathcal{C}(X))^\sharp : \quad C^{\sharp n} \leftarrow (C^n/C^0)^\sharp \leftarrow (C^n/C^1)^\sharp \leftarrow \cdots \leftarrow (C^n/C^{n-1})^\sharp, \quad (3.41)$$

$$\mathcal{C}^\perp : \quad C^{\perp 0} \rightarrow C^{\perp 1} \rightarrow C^{\perp 2} \rightarrow \cdots \rightarrow C^{\perp n}, \quad (3.42)$$

where we set $C^{\perp i} := (C^n/C^{n-i})^\sharp$. It is already known that $\sigma_0^\sharp, \sigma_1^\sharp \ldots, \sigma_n^\sharp$ form a basis of $C^{\sharp n}$, dual to the generators $\sigma_0, \sigma_1, \ldots, \sigma_n$ of $C^n$. We define $\tau_i = \sigma_{n+1-i}^\sharp$, leading to

$$C^{\perp i} = \langle \sigma_n^\sharp, \sigma_{n-1}^\sharp, \ldots, \sigma_{n+1-i}^\sharp \rangle = \langle \tau_1, \tau_2, \ldots, \tau_i \rangle \quad (3.43)$$

through elementary linear algebra, as can be seen in Example 3.3.12.

**Proposition 3.3.10.** *Let $\partial$ be the boundary operator in homology. If $\partial\sigma_j = \sum_{i<j} D_{ij}\sigma_i$, then the corresponding coboundary operator $\delta$ in cohomology is given by $\delta\tau_j = \sum_{i<j} D_{ij}^\perp \tau_i$, where $D_{ij}^\perp = D_{(n-j),(n-i)}$ is the transposed and reversed version of the matrix $D_{ij}$.*

*Proof.* Recall the duality between homology and cohomology. The boundary operator $\partial$ in homology is defined by its action on the basis elements $\sigma_j \in C^n$ as $\partial\sigma_j = \sum_{i<j} D_{ij}\sigma_i$, where $D_{ij}$ are the coefficients of the boundary map. The corresponding coboundary operator $\delta$ in cohomology acts on the dual space $C^{\sharp n}$, which is spanned by the dual basis $\{\sigma_0^\sharp, \sigma_1^\sharp, \ldots, \sigma_n^\sharp\}$, where $\sigma_i^\sharp(\sigma_k) = \delta_{ik}$. Define $\tau_i = \sigma_{n+1-i}^\sharp$, which reverses the index order in the dual space. The coboundary operator $\delta$ is the adjoint of the boundary operator $\partial$, meaning that for any $\tau_j \in C^{\sharp n}$, we have $\delta\tau_j(\sigma_k) = \tau_j(\partial\sigma_k)$. Substituting the expression for $\partial\sigma_j$, we get $\partial\sigma_j = \sum_{i<j} D_{ij}\sigma_i$ and

$$\delta\tau_j(\sigma_k) = \tau_j \left( \sum_{i<j} D_{ij}\sigma_i \right) = \sum_{i<j} D_{ij}\tau_j(\sigma_i). \quad (3.44)$$

By definition of the dual basis, $\tau_j(\sigma_i) = \delta_{n+1-j,i}$, which means the action of $\delta$ on $\tau_j$ depends on the corresponding entries of the boundary matrix $D_{ij}$. The coboundary operator $\delta$ is thus governed by the transposed and index-reversed matrix $D^\perp$, where $D_{ij}^\perp = D_{(n-j),(n-i)}$. Therefore, the action of $\delta$ on $\tau_j$ is given by $\delta\tau_j = \sum_{i<j} D_{ij}^\perp \tau_i$. $\square$

**Remark 3.3.11.** *[9, §2.7] Due to the structural similarity of $\mathcal{C}^{\perp}(X)$ and $\mathcal{C}(X)$, persistent absolute cohomology is interpreted as 'relative persistent relative cohomology'.*

**Example 3.3.12.** *Consider the construction explicitly for $C^2$ with basis $\{\sigma_0, \sigma_1, \sigma_2\}$. The dual basis vectors in $C^{\sharp 2}$ are $\{\sigma^{\sharp 0}, \sigma^{\sharp 1}, \sigma^{\sharp 2}\}$, and we set $\tau_0 = \sigma^{\sharp 2}$, $\tau_1 = \sigma^{\sharp 1}$, and $\tau_2 = \sigma^{\sharp 0}$. For $C^{\perp 1}$, we obtain:*

$$C^{\perp 1} = \langle \sigma^{\sharp 2}, \sigma^{\sharp 1} \rangle = \langle \tau_0, \tau_1 \rangle.$$

*If the boundary map for $\partial \sigma_2$ is given by:*

$$\partial \sigma_2 = D_{20} \sigma_0 + D_{21} \sigma_1,$$

*then the coboundary operator $\delta$ acts on $\tau_2$ as:*

$$\delta \tau_2 = D_{21}^{\perp} \tau_1 + D_{20}^{\perp} \tau_0,$$

*where $D_{21}^{\perp} = D_{12}$ and $D_{20}^{\perp} = D_{02}$ are the transposed entries.*

**Proposition 3.3.13.** *[9, Proposition 2.8] The persistence module $(C^{\perp n}/C^{\perp}(X))$ is the reverse of $\mathcal{C}^{\sharp}(X)$, and the respective coboundary maps agree.*

*Proof.* Each submodule $C^{\perp i}$ is defined as the dual of the quotient $C^n/C^{n-i}$, i.e., $C^{\perp i} = (C^n/C^{n-i})^{\sharp}$. Thus, the persistence module $(C^{\perp n}/C^{\perp}(X))$ is constructed by taking the quotient of $C^{\perp n}$ by $\mathcal{C}^{\perp}(X)$, which gives the reversed filtration:

$$C^{\sharp n} \leftarrow (C^n/C^0)^{\sharp} \leftarrow (C^n/C^1)^{\sharp} \leftarrow \cdots \leftarrow (C^n/C^{n-1})^{\sharp}. \qquad (3.45)$$

This is the filtration of $\mathcal{C}^{\sharp}(X)$, but in reverse order. Therefore, the persistence module $(C^{\perp n}/C^{\perp}(X))$ is the reverse of $\mathcal{C}^{\sharp}(X)$.

Consider the boundary operator $\partial$ on the chain complex $C^n$, given by $\partial \sigma_j = \sum_{i<j} D_{ij}\sigma_i$, where $D_{ij}$ are the coefficients of the boundary operator. In the dual space $C^{\sharp n}$, the coboundary operator $\delta$ acts on the dual basis elements $\tau_j = \sigma_{n+1-j}^{\sharp}$. The coboundary map is induced by the transposed and reversed matrix of the boundary operator $\delta \tau_j = \sum_{i<j} D_{ij}^{\perp} \tau_i$, where $D_{ij}^{\perp} = D_{(n-j),(n-i)}$, the transposed and index-reversed version of the boundary matrix $D_{ij}$. Since the filtration in $(C^{\perp n}/C^{\perp}(X))$ corresponds to the reverse of $\mathcal{C}^{\sharp}(X)$, and the coboundary operators are induced by the same transposed matrix $D^{\perp}$. $\qquad \square$

**Corollary 3.3.14.** *[9, Proposition 2.10] Let $\{\hat{\sigma}_0^{\sharp}, \ldots, \hat{\sigma}_n^{\sharp}\}$ be the dual of $\{\hat{\sigma}_0, \ldots, \hat{\sigma}_n\}$, where $C^i = \langle \hat{\sigma}_0, \ldots, \hat{\sigma}_i \rangle$ for $0 \leq i \leq n$. Then $\hat{\tau}_i = \hat{\sigma}_{n+1-i}^{\sharp}$ for all i.*

*Proof.* We use Theorem 3.3.5 and obtain a partition of indices $\{0, 1, 2, \ldots, n\} = B^{\times} \sqcup B \sqcup D$ for $\mathcal{C}^{\perp}(X)$ and new generators $\hat{\tau}_i$, where:

$$\delta \hat{\tau}_i = 0 \quad \text{for } i \in B \sqcup D, \qquad \delta \hat{\tau}_d = \hat{\tau}_b \quad \text{for pairs } (b, d) \in \text{Pairs}(\mathcal{C}^{\perp}(X), \delta). \quad (3.46)$$

Now, let $\{\hat{\sigma}_0^{\sharp}, \ldots, \hat{\sigma}_n^{\sharp}\}$ denote the dual basis of $\{\hat{\sigma}_0, \ldots, \hat{\sigma}_n\}$, such that $\hat{\sigma}_i^{\sharp}(\hat{\sigma}_j) = \delta_{ij}$ for all $i, j$. Since $\mathcal{C}^{\perp}(X)$ is the dual filtration of $\mathcal{C}(X)$, we observe that the filtration structure in $\mathcal{C}^{\perp}(X)$ is reversed. Specifically, the duality induces the relationship between the generators $\hat{\tau}_i$ and the elements of the dual basis $\hat{\tau}_i = \hat{\sigma}_{n+1-i}^\sharp$. This indexing reversal arises naturally from the dual filtration structure. The coboundary operator $\delta$ respects this duality, so that for any pair $(b, d) \in \text{Pairs}(\mathcal{C}^\perp(X), \delta)$, we have $\delta \hat{\tau}_d = \hat{\tau}_b$, preserving the dual pairing structure of the persistence module. Thus, for each $i$, the relation $\hat{\tau}_i = \hat{\sigma}_{n+1-i}^\sharp$ holds, completing the proof. $\qquad\square$

# Index